\documentclass{article}
\usepackage{dirtytalk}

\usepackage{svg}
\usepackage{amsmath}
\usepackage{amssymb}
\usepackage{amsthm}
\usepackage{slashed}
\usepackage{cite}
\usepackage{color} 
\usepackage{caption}
\usepackage{euscript}
\usepackage{comment}
\usepackage[arrow,curve,matrix,arc,2cell]{xy}
\usepackage{graphicx}    
\usepackage[utf8]{inputenc}
\usepackage[unicode]{hyperref}
\usepackage[a4paper, margin=1.5in]{geometry}
\usepackage{pifont}
\UseAllTwocells
\DeclareFontFamily{U}{rsfs}{} \DeclareFontShape{U}{rsfs}{n}{it}{<->
rsfs10}{} \DeclareSymbolFont{mscr}{U}{rsfs}{n}{it}
\DeclareSymbolFontAlphabet{\scr}{mscr}
\def\mathscr{\scr}
\begin{document}
%%%%%%%%%%%%%%%%%%%%%%%%%%%%%%%%%%%%%%%%%%%%%%%%%%%%%%%%%%%%%%%%%%%%%%%%
%%%%%%%%%%%%%%%%%%%%%%%%%%     Macros      %%%%%%%%%%%%%%%%%%%%%%%%%%%%%
%%%%%%%%%%%%%%%%%%%%%%%%%%%%%%%%%%%%%%%%%%%%%%%%%%%%%%%%%%%%%%%%%%%%%%%%

\newcommand{\red}[1]{\textcolor{red}{#1}}
\def\e#1\e{\begin{equation}#1\end{equation}}
\def\ea#1\ea{\begin{align}#1\end{align}}
\def\eas#1\eas{\begin{align*}#1\end{align*}}
\def\eq#1{{\rm(\ref{#1})}}
\newenvironment{claim}[1]{\par\noindent\underline{Claim:}\space#1}{}
\newenvironment{claimproof}[1]{\par\noindent\underline{Proof:}\space#1}{\hfill $\blacksquare$}
\theoremstyle{plain}% default
\newtheorem{thm}{Theorem}[section]
\newtheorem{prop}[thm]{Proposition}
\newtheorem{lem}[thm]{Lemma}
\newtheorem{cor}[thm]{Corollary}
\newtheorem{quest}[thm]{Question}
\newtheorem{conj}[thm]{Conjecture}
\theoremstyle{definition}
\newtheorem{dfn}[thm]{Definition}
\newtheorem{ex}[thm]{Example}
\newtheorem{rem}[thm]{Remark}
\numberwithin{equation}{section}

\def\Amp{\mathop{\rm Amp}\nolimits}
\def\loc{{\mathop{\rm loc}\nolimits}}
\def\ope{{\mathop{\rm op}\nolimits}}

\def\cyl{{\mathop{\rm cyl}\nolimits}}
\def\con{{\mathop{\rm con}\nolimits}}

\def\cay{\mathop{\rm cay}\nolimits}
\def\Cay{\mathop{\rm Cay}\nolimits}
\def\ACay{\mathop{\rm ACay}\nolimits}

\def\glu{\mathop{\,\sharp\,}}
\def\glue#1{\mathop{\,\sharp_{#1}\,}}
\def\ob{{\mathop{\rm ob}\nolimits}}

\def\CayD{\slashed{D}}
\def\delbar{\bar{\partial}}

\def\ACyl{{\mathop{\rm ACyl}\nolimits}}
\def\AC{{\mathop{\rm AC}\nolimits}}
\def\CS{{\mathop{\rm CS}\nolimits}}
\def\supp{\mathop{\rm supp}\nolimits}
\def\dist{\mathop{\rm dist}\nolimits}
\def\sgn{\mathop{\rm sgn}\nolimits}
\def\dim{\mathop{\rm dim}\nolimits}
\def\Ker{\mathop{\rm Ker}}
\def\Coker{\mathop{\rm Coker}}
\def\rank{\mathop{\rm rank}}
\def\Ho{H}
\def\sign{\mathop{\rm sign}\nolimits}
\def\id{\mathop{\rm id}\nolimits}
\def\dvol{\mathop{\rm dvol}\nolimits}
\def\spn{\mathop{\rm span}\nolimits}
\def\SO{\mathop{\rm SO}\nolimits}
\def\Orth{\mathop{\rm O}\nolimits}
\def\Fr{\mathop{\rm Fr}\nolimits}
\def\Gr{\mathop{\rm Gr}\nolimits}
\def\cay{\mathop{\rm cay}\nolimits}
\def\fix{\mathop{\rm fix}\nolimits}
\def\inj{\mathop{\rm inj}\nolimits}
\def\SF{\mathop{\rm SF}\nolimits}
\def\Or{\mathop{\rm Or}\nolimits}
\def\ad{\mathop{\rm ad}\nolimits}
\def\Hom{\mathop{\rm Hom}\nolimits}
\def\Map{\mathop{\rm Map}\nolimits}
\def\Crit{\mathop{\rm Crit}\nolimits}
\def\ev{\mathop{\rm ev}\nolimits}
\def\Univ{\mathop{\rm Univ}\nolimits}
\def\Fix{\mathop{\rm Fix}\nolimits}
\def\Hol{\mathop{\rm Hol}\nolimits}
\def\Iso{\mathop{\rm Iso}\nolimits}
\def\Hess{\mathop{\rm Hess}\nolimits}
\def\Stab{\mathop{\rm Stab}\nolimits}
\def\Pd{\mathop{\rm Pd}\nolimits}
\def\Aut{\mathop{\rm Aut}\nolimits}
\def\Diff{\mathop{\rm Diff}\nolimits}
\def\boFlag{\mathop{\rm Flag}\nolimits}
\def\boFlagSt{\mathop{\rm FlagSt}\nolimits}
\def\dOrb{{\mathop{\bf dOrb}}}
\def\dMan{{\mathop{\bf dMan}}}
\def\mKur{{\mathop{\bf mKur}}}
\def\Kur{{\mathop{\bf Kur}}}
\def\Pic{\mathop{\rm Pic}}
\def\Re{\mathop{\rm Re}}
\def\Im{\mathop{\rm Im}}
\def\re{\mathop{\rm re}}
\def\im{\mathop{\rm im}}
\def\SU{\mathop{\rm SU}}
\def\Sp{\mathop{\rm Sp}}
\def\Spin{\mathop{\rm Spin}}
\def\GL{\mathop{\rm GL}}
\def\ind{\mathop{\rm ind}}
\def\area{\mathop{\rm area}}
\def\U{{\rm U}}
\def\vol{\mathop{\rm vol}\nolimits}
\def\virt{{\rm virt}}
\def\emb{{\rm emb}}
\def\bs{\boldsymbol}
\def\ge{\geqslant}
\def\le{\leqslant\nobreak}
\def\O{{\mathbin{\mathcal O}}}
\def\cA{{\mathbin{\mathcal A}}}
\def\cB{{\mathbin{\mathcal B}}}
\def\cC{{\mathbin{\mathcal C}}}
\def\cD{{\mathbin{\scr D}}}
\def\cDHS{{\mathbin{\scr D}_{\Q HS}}}
\def\cE{{\mathbin{\mathcal E}}}
\def\boE{{\mathbin{\mathbf E}}}
\def\cF{{\mathbin{\mathcal F}}}
\def\boF{{\mathbin{\mathbf F}}}
\def\cG{{\mathbin{\mathcal G}}}
\def\cH{{\mathbin{\mathcal H}}}
\def\cI{{\mathbin{\mathcal I}}}
\def\cJ{{\mathbin{\mathcal J}}}
\def\cK{{\mathbin{\mathcal K}}}
\def\cL{{\mathbin{\mathcal L}}}
\def\cM{{\mathbin{\mathcal M}}}
\def\brM{\,\,\overline{\!\!\cM}}
\def\bcM{{\mathbin{\bs{\mathcal M}}}}
\def\cN{{\mathbin{\mathcal N}}}
\def\cO{{\mathbin{\mathcal O}}}
\def\cP{{\mathbin{\mathcal P}}}
\def\boR{{\mathbin{\mathbf R}}}
\def\cS{{\mathbin{\mathcal S}}}
\def\cT{{\mathbin{\mathcal T}}}
\def\cU{{\mathbin{\mathcal U}}}
\def\cQ{{\mathbin{\mathcal Q}}}
\def\boQ{{\mathbin{\mathbf Q}}}
\def\cW{{\mathbin{\mathcal W}}}
\def\C{{\mathbin{\mathbb C}}}
\def\bQ{{\mathbin{\mathbb Q}}}
\def\bV{{\mathbin{\mathbb V}}}
\def\bE{{\mathbin{\mathbb E}}}
\def\bD{{\mathbin{\mathbb D}}}
\def\boF{{\mathbin{\mathbf F}}}
\def\bF{{\mathbin{\mathbb F}}}
\def\H{{\mathbin{\mathbb H}}}
\def\N{{\mathbin{\mathbb N}}}
\def\Q{{\mathbin{\mathbb Q}}}
\def\R{{\mathbin{\mathbb R}}}
\def\bS{{\mathbin{\mathbb S}}}
\def\Z{{\mathbin{\mathbb Z}}}
\def\sF{{\mathbin{\mathscr F}}}
\def\phi{\varphi}
\def\al{\alpha}
\def\be{\beta}
\def\ga{\gamma}
\def\de{\delta}
\def\io{\iota}
\def\ep{\epsilon}
\def\eps{\epsilon}
\def\la{\lambda}
\def\ka{\kappa}
\def\th{\theta}
\def\ze{\zeta}
\def\up{\upsilon}
\def\vpphi{\varphi}
\def\si{\sigma}
\def\om{\omega}
\def\De{\de}
\def\La{\Lambda}
\def\Si{\Sigma}
\def\Th{\Theta}
\def\Om{\Omega}
\def\Ga{\Gamma}
\def\Up{\Upsilon}
\def\pd{\partial}
\def\ts{\textstyle}
\def\st{\scriptstyle}
\def\sst{\scriptscriptstyle}
\def\w{\wedge}
\def\sm{\setminus}
\def\bu{\bullet}
\def\op{\oplus}
\def\ot{\otimes}
\def\ov{\overline}
\def\ul{\underline}
\def\bigop{\bigoplus}
\def\bigot{\bigotimes}
\def\iy{\infty}
\def\es{\emptyset}
\def\ra{\rightarrow}
\def\Ra{\Rightarrow}
\def\Longra{\Longrightarrow}
\def\ab{\allowbreak}
\def\longra{\longrightarrow}
\def\hookra{\hookrightarrow}
\def\dashra{\dashrightarrow}
\def\t{\times}
\def\ci{\circ}
\def\ti{\tilde}
\def\d{{\rm d}}
\def\dt{{\rm dt}}
\def\D{{\rm D}}
\def\Lie{{\mathcal{L}}}
\def\ha{{\ts\frac{1}{2}}}
\def\md#1{\vert #1 \vert}
\def\bmd#1{\big\vert #1 \big\vert}
\def\ms#1{\vert #1 \vert^2}
\def\nm#1{\Vert #1 \Vert}
%%%%%%%%%%%%%%%%%%%%%%%%%%%%%%%%%%%%%%%%%%%%%%%%%%%%%%%%%%%%%%%%%%%%%%%%
%%%%%%%%%%%%%%%%%%%%%% Text of paper %%%%%%%%%%%%%%%%%%%%%%%%%%%%%%%%%%%
%%%%%%%%%%%%%%%%%%%%%%%%%%%%%%%%%%%%%%%%%%%%%%%%%%%%%%%%%%%%%%%%%%%%%%%%
\title{Conically singular Cayley submanifolds III: Fibrations}
\author{Gilles Englebert}
\date{\today}
\maketitle

\setcounter{tocdepth}{2}
\tableofcontents
\section{Introduction}
This is the third and last in a series of three papers working towards the construction of non-trivial Cayley fibrations using gluing methods. In the previous two papers we described the deformation and desingularisation theory of conically singular Cayley submanifolds. This is all in an effort to understand the perturbative properties of Cayley fibrations of compact $\Spin(7)$-manifolds. Indeed, experience suggests that a fibration of a compact manifold by compact Cayleys must necessarily include singular fibres. However there is hope that the singular fibres may sometimes admit only conical singularities, which are analytically tractable. Results such as the gluing theorem \cite[Thm. 3.15]{englebertConicallySingularCayley2023} then allow us to investigate the deformation theory of fibrations with singular fibres and deduce properties such as stability (in a suitable sense) under small deformations of the $\Spin(7)$-structure. In this paper we will show two stability results, the first of which is the following, about weak fibrations:
\begin{thm}
\label{1_1_fibration_weak}
Let $(M, \Phi)$ be an almost $\Spin(7)$-manifold which is weakly fibred by conically singular Cayleys with semi-stable cones, such that all the Cayleys in the fibration are unobstructed. Let $\Phi_t$ be a smooth deformation of the $\Spin(7)$-structure. Then there is an $\eps > 0$ such that for all $t \in [0, \eps)$ the manifold $(M, \Phi_t)$ can still be weakly fibred. 
\end{thm}
Here semi-stability imposes some mild restrictions on the critical rates of the asymptotic Cayley cone. The notion of weak stability is homological in nature. It is weaker than the usual notion of fibration, since fibres may intersect. The main idea in the proof is that unobstructed Cayleys, both non singular and conically singular, deform smoothly under smooth change of the ambient $\Spin(7)$-structure. Both the deformation theory of Cayley submanifolds and the desingularisation theory play a role here. First, the deformation theory ensures that Cayleys of a given type continue to exist, deform by isotopies, and do not form additional singularities. The gluing theorem ensures that we have a precise quantitative grasp on the geometry as one passes from nearly singular Cayleys in the perturbation to their singular limit. Essentially, this is a result that makes minimal assumptions on the geometry, but is consequently rather weak and mostly topological in nature. We have the following result on strong Cayley fibrations, which needs further assumptions on the geometry.
\begin{thm}
\label{1_1_fibration_strong}
Let $(M, \Phi)$ be an almost $\Spin(7)$-manifold which is strongly fibred by simple conically singular Cayleys, such that all the Cayleys in the fibration are unobstructed. Assume that the fibration is non-degenerate. Let $\Phi_t$ be a smooth deformation of the $\Spin(7)$-structure. Then there is an $\eps > 0$ such that for all $t \in [0, \eps)$ the manifold $(M, \Phi_t)$ can still be strongly fibred. 
\end{thm}
Simpleness and non-degeneracy are technical conditions that allow us to analyse the deformations of Cayleys near the singular limit in detail and avoid certain pathological behaviours. Simpleness ensures that the deformations of the nearly singular compact Cayleys (which in a fibration should be a four-dimensional space) can be described by deformations dominated by the movement of their cone and deformations at a unique other, lower rate. This ensures Fredholmness when we turn to solving the Cayley equation (i.e. the linear p.d.e. governing the first-order deformations of a Cayley) over the manifolds we obtain from the gluing theorem. Non-degeneracy is the requirement that the solutions to the Cayley equations, weighted appropriately, can be bound away from zero. We will show that the fibrations in this situation will remain fibering up to first order, and the weak stability result then allows us to show that this implies global stability.

These two results in turn can then be used to construct fibrations on compact manifolds via gluing of non-compact pieces, as in the programme of Kovalev \cite{KovalevFibration}. We review the examples proposed by Kovalev coming from the twisted connected sum construction and show that they satisfy the conditions of Theorem \ref{1_1_fibration_strong}. Thus in particular we prove:
\begin{thm}
\label{1_1_existence}
There are compact, torsion-free $\Spin(7)$-manifolds of holonomy $G_2$ which admit fibrations by Cayley submanifolds.
\end{thm}
From this we can deduce as a corollary:
\begin{cor}
\label{1_1_existence_g2}
There are compact, torsion-free $G_2$-manifolds of full holonomy $G_2$ which admit fibrations by coassociative submanifolds.
\end{cor}
Note that the examples in the $\Spin(7)$ case are not of maximal holonomy $\Spin(7)$. This is due to our current lack of fibrations on gluing constructions of full holonomy $\Spin(7)$ manifolds.

\subsubsection*{Notation}
We will denote by $C$ an unspecified constant, which may refer to different constants within the same derivation. To indicate the dependence of this constant on quantities $x, y, \dots$, we will write $C(x, y, \dots)$. Similarly, if an inequality holds up to an unspecified constant, we will write $A \lesssim B$ instead of $A \le C B$.
\subsubsection*{Acknowledgements}
This research has been supported by the Simons Collaboration on Special Holonomy in Geometry, Analysis, and Physics. I want to thank my DPhil supervisor Dominic Joyce for his excellent guidance during this project.

\section{Preliminaries}
\subsection{Calibrated Geometry}
We will start by giving a quick outline of calibrated geometry which as an area of study has first been proposed by Harvey and Lawson \cite{HarvLaws}. As an area it is intimately related to study of manifolds of special holonomy, such as the $\Spin(7)$-manifolds that are the objects of interest of this paper. An excellent introduction which goes beyond what we discuss is the book by Joyce \cite{joyceRiemannianHolonomyGroups2007}. 

Let now $(M,g)$ be a Riemannian manifold. Suppose that $\phi \in \Om^k(M)$ is a closed form, such that at each point $p\in M$ and for each $k$-plane $\Pi \in \Gr(T_pM, k)$ the \textbf{calibration inequality}:
\eas
\label{calibration_inequality}
\phi\lvert_\Pi \ \le \dvol_\Pi
\eas 
is satisfied. We then call $\phi$ a \textbf{calibration}. We say that a $k$-dimensional submanifold $N \subset M$ is $\phi$-\textbf{calibrated} if the calibration inequality becomes an equality, i.e.:
\eas
\phi\lvert_N \ = \dvol_N.
\eas
Any calibrated $N$ is volume minimizing in its homology class, which can be seen by an application of Stokes' theorem. Indeed, for $\tilde{N}$ homologous to $N$ we see:
\eas
\vol(N) = \int_N \dvol_N = \int_N \phi = \int_{\tilde{N}}\phi \le \int_{\tilde{N}} \dvol_{\tilde{N}} \le \vol(\tilde{N}).
\eas

We will now review the fundamentals of three calibrated geometries, namely the Calabi--Yau, $G_2$ and $\Spin(7)$ geometries. In all cases this will entail the study of calibrated submanifolds of course, but in fact admitting calibrations of these types also places restrictions on the holonomy of the ambient space, which will lie in $\SU(n), G_2$ and $\Spin(7)$ respectively.
\begin{rem}
In the $G_2$ and $\Spin(7)$ case we do not always require the differential form to be closed. Consequently the form is not a calibration and its calibrated manifolds are not necessarily minimal. 
\end{rem}

\subsection{Calabi--Yau and Fano Geometry}  
\subsubsection*{Calabi--Yau Geometry}
We briefly review some aspects of Calabi--Yau and Fano manifolds which will be relevant to our discussion of Cayley fibrations of $\Spin(7)$ manifolds. For a more in-depth introduction we refer to  \cite[Ch. 6]{joyceCompactManifoldsSpecial2000}.
\begin{dfn}[Calabi--Yau manifold]
Let $(X^{2n}, J, \om, g)$ be a K\"ahler manifold of complex dimension $n$ which admits a nowhere vanishing holomorphic $(n,0)$-form $\Om$. The line bundle of $(n,0)$-forms is called the \textbf{canonical bundle}, so equivalently we may require $(X,J)$ to have holomorphically trivial canonical bundle. If we furthermore have the following normalisation condition which links the complex and symplectic geometry of $X$:
\e
\label{2_1_cy_condition}
\frac{\om^n}{n!} = (-1)^{n(n-1)/2}(i/2)^n \Om \wedge \bar{\Om},
\e
then we call $(X^{2n}, J,\om,g, \Om)$ a \textbf{Calabi--Yau manifold}. The form $\Om$ is called the \textbf{holomorphic volume form}.
\end{dfn}
The holonomy of any Calabi--Yau manifold is contained in $\SU(n)$ and the metric 	$g$ is necessarily Ricci-flat. At any point an $\SU(n)$-structure is isomorphic to the following standard model on $\C^n$ with complex coordinates $z_1 = x_1 +iy_1, \dots, z_n = x_n + iy_n$:
\eas
\om_0 &=  \d x_1 \wedge \d y_1 + \dots \d x_n \wedge \d y_n, \\
g_0 &= \d x_1^2 + \d y_1^2 + \dots +\d x_n^2 + \d y_n^2, \\
\Om_0 &= \d z_1 \wedge \dots \wedge \d z_n.
\eas
We have the following theorem due to Yau (proving a conjecture due to Calabi) which reduces the existence of a Calabi--Yau structure to a question of complex geometry on $(X,J)$. 
\begin{thm}
\label{yau_theorem}
Let $(X,J, \om, g)$ be a compact K\"ahler manifold with trivial canonical bundle. Then there is a unique K\"ahler form $\tilde{\om}$ in the cohomology class of $\om$ (with corresponding metric $\tilde{g}$) and a holomorphic volume form $\Om$ such that $(X,J,\tilde{\om},\tilde{g}, \Om)$ is a Calabi--Yau manifold.
\end{thm}

On Calabi--Yau manifolds there are two calibrations of interest. First we have the real part of the holomorphic volume form $\Re \Om \in \Om^n(M)$, whose calibrated submanifolds are the so called \textbf{special Lagrangians}, which are difficult to construct, and we will not go further into discussing them here. Secondly we have the complex submanifolds in any dimension $1 \le k \le n$, which are calibrated by the form $\frac{\om^k}{k!}$. The calibration inequality is in this case also called the Wirtinger inequality and is a feature of any K\"ahler manifold, not just Calabi--Yau manifolds.

In two complex dimensions, Calabi--Yau manifolds are particularly well understood. Their underlying complex surfaces must either be tori $T^4$ or so called \textbf{K3 surfaces}, which are the only two deformation types of complex surfaces with trivial canonical bundle. We will discuss K3 surfaces in more detail now, see \cite[Section 7.3.3]{joyceCompactManifoldsSpecial2000} for a more in depth discussion. By a result of Kodaira all complex analytic K3 surfaces $K$ belong to a single diffeomorphism type, namely that of a quartic $\{x_0^4 + x_1^4 +x_2^4 + x_3^4 = 0\} \subset \C P^3$. In particular they are simply connected and all have isomorphic cohomology groups, the only non-trivial one being $\Ho^2(K,\Z)$. Since $K$ is a compact closed four-manifold its second cohomology admits a non-degenerate intersection pairing, and this lattice we denote by $\La$. 
Next we recall that analytic K3 surfaces form a $20$-dimensional moduli space. To see this explicitly we define a \textbf{marked K3 surface} to be a K3 surface $K$ together with a choice of lattice isomorphism $h:\Ho^2(K,\Z) \ra \La$. The complex structure of the K3 surface will then be determined locally by its Hodge structure (i.e. how $\Ho^2(K, \C) = \Ho^{0,2}(K, \C)\op \Ho^{1,1}(K, \C)\op \Ho^{2,0}(K, \C)$ splits with respect to the marking $h$). More precisely we define the so called \textbf{period domain}: 
\ea
\label{K3_period_domain}
D_{K3} &= \{ u \in P(\La \ot \C): u^2 = 0,\  u \cdot \bar{u} > 0 \} \\ \nonumber
&\simeq  \{ \Pi \subset \La \ot \R: \langle \cdot, \cdot \rangle\lvert_\Pi > 0\} \subset \Gr_+(2, \La \ot \R). 
\ea
This is the space of all possible complex lines $h(\Ho^{2,0}(K,\C))$ in $\La \ot \C$. The map sending $(K, h)$ to $[h(\Ho^{2,0}(K,\C))] \in D_{K3}$ is called the \textbf{period map}. It is a local but not global diffeomorphism, in particular because the moduli space of marked K3 surfaces is not Hausdorff (while $D_{K3}$ is). The isomorphism to $\Gr_+(2, \La\ot \R)$ follows from identifying $\Ho^{2,0}\op \Ho^{0,2} = \Pi \ot \C$ for a real two-plane $\Pi$. 

Next, for our discussion we will need K3 surfaces with additional structure, so called \textbf{lattice polarised K3 surfaces} (see \cite{beauvilleFanoThreefoldsK32002}). For this we look at the \textbf{Picard group} $\Pic(K,J)$, which is the abelian group of holomorphic line bundles under the tensor product. As K3 surfaces are simply connected, we can think of the Picard group as being embedded in $\Ho^2(K,\Z)$ via the first Chern class $c_1: \Pic(K,J) \ra \Ho^{1,1}(K,\Z)$. Thus, while the intersection form on $\Ho^{2}$ is a topological invariant, we can restrict it to the Picard group to get an invariant of the complex structure, the \textbf{Picard lattice}. This is a lattice of rank $(1,\rho-1)$ where $0 \le \rho \le 20$ is the rank of the Picard lattice. 

Assume now that we are given a sublattice $N \subset \La$ of signature $(1, r-1)$ and an element $A \in N$ with $A\cdot A = 2g-2 > 0 $. We say that a marked K3 surface $(K,J, h)$ is $(N,A)$-\textbf{polarised} if $h(N) \subset \Pic(K,J)$, this embedding is primitive, meaning that $\La/N$ is torsion-free, and $h(A) \in \Pic(K,J)$ is ample. The number $g$ is then called the \textbf{genus} of the polarised K3 $K$. Similar to the period domain of marked K3 surfaces \eqref{K3_period_domain} one can describe a similar period domain for marked polarised K3 surfaces. For this note that as $h(N) \subset \Pic(K,J) \subset \Ho^{1,1}(K)$ the complex line $\Ho^{2,0}(K)$ must be orthogonal to $h(N)$. This motivates the definition of the following domain:
\ea
\label{K3_polarised_period_domain}
D_{N} &= \{ u \in P(N^{\perp} \ot \C): u^2 = 0,\  u \cdot \bar{u} > 0\} \\ \nonumber
&\simeq \{ \Pi \subset N^{\perp} \ot \R: \langle \cdot, \cdot \rangle\lvert_\Pi > 0\} \subset \Gr_+(2, N^\perp \ot \R). 
\ea
The corresponding Torelli theorem states that the period map from before maps the moduli space of marked $(N,A)$-polarised K3 surfaces $\cK^{N,A}$  to $D_{N}$ by a local diffeomorphism. Hence this moduli space has dimension $20-r$. 

The K\"ahler geometry of K3 surfaces is also rather explicit. Suppose that the (non-polarised) K3 surface $(K,J,\om,g)$ has period point $\Pi\in D_{K3}$. We then define the \textbf{root system} corresponding to $\Pi$ as:
\eas
\Delta_{\Pi} = \{ \la \in \La: \la \cdot \la = -2,\  \la \cdot p = 0 \ \ \forall p \in \Pi\}.
\eas
Then the set of \textbf{K\"ahler chambers} of the K3 surface is given by:
\e
\label{K3_Kaehler_chambers}
\{\om \in \La \ot \R: \om \cdot \om >0,\  \om \cdot p =0 \text{ for } p\in \Pi,\  \om \cdot \la \ne 0\  \ \forall \la \in \Delta_{\Pi}\}.
\e
Now the K\"ahler cone is always a connected component of the set of K\"ahler chambers, and thus in particular an open subset of $\Ho^{1,1}(K)$. 
After our discussion of the complex and K\"ahler geometry of a K3 surface, consider now a K3 surface $(S, \om_I, I, g, \Om_I)$ with a chosen Calabi--Yau structure. By Yau's Theorem \ref{yau_theorem} we see that $\om_I, g$ and $\Om_I$ are determined by the complex structure $I$ and the cohomology class $[\om] \in \Ho^2(S)$. We can then write $\Om_I = \om_J + i\om_K$. As suggested by the notation $S$ is also K\"ahler with respect to the forms $\om_J$ and $\om_K$ for new complex structures $J$ and $K$ (meaning that $g(\cdot, \cdot) = \om_J(\cdot, J\cdot) = \om_K(\cdot, K\cdot)$). The three complex structure satisfy the quaternionic relations $I^2 = J^2 = K^2 = IJK = -1$. In fact, for $(a,b,c) \in S^2 \subset \R^3$ any linear combination $aI + bJ + cK$ determines a further complex structure for which $(S,g)$ is K\"ahler for a suitably chosen K\"ahler form. Riemannian manifolds that are K\"ahler in three compatible ways like above are called \textbf{hyperk\"ahler} manifolds. In the K3 case, we can describe the K3 moduli space explicitly:

\begin{prop}
\label{HK_moduli}
The moduli space $\cM^{\text{hk}}$ of hyperk\"ahler K3 surfaces admits a period map:
\e
\label{HK_K3_period}
\cP^{\text{hk}}: \cM^{\text{hk}} \longra D^{{\text{hk}}}_{K3}.
\e
Here $D^{\text{hk}}_{K3}$ is defined as:
\ea
\label{HK_K3_domain}
D^{\text{hk}}_{K3} =& \{(\al_1, \al_2, \al_3): \al_{i} \in  \La \ot \R,\  \al_i \cdot \al_j = a\de_{ij} \text{ with } a >0, \\ \nonumber
 &\text{ for each } \la \in \La \text{ with } \la \cdot \la = -2 \text{ there is } i = 1,2 \text{ or } 3 \text{ such that } \al\cdot \la \ne 0\}.
\ea
\end{prop}

Hyperk\"ahler manifolds admit isometries of a special kind which interchange the complex structures, called \textbf{hyperk\"ahler rotations}. More formally, for us a hyperk\"ahler rotation will be an isometry $\phi: S_1 \ra S_2$ between K3 surfaces $S_1$ and $S_2$ with complex structures $I_1, J_1, K_1$ and $I_2,J_2,K_2$ respectively, so that
\e
\label{hyperkaehler_twist}
  \phi^*I_2 = J_1, \quad \phi^*J_2 = I_1, \text{ and } \phi^*K_2 = -K_1.
\e
Alternatively we can define hyperk\"ahler rotations by their actions on the K\"ahler forms. Indeed the HK rotation $\phi$ from above induces the following action on the K\"ahler forms $(\om_+, \om_-, \om_0)$ corresponding to the distinguished complex structures $(I,J,K)$ of a K3 surface $S$:
\e
\label{HK_rotation_action}
(\om_+, \om_-, \om_0) \longmapsto (\om_-, \om_+, -\om_0).
\e
These special isometries will important in section \ref{4_3_twisted_connected}, when we discuss the construction of $G_2$-manifolds from Calabi--Yau pieces. We will glue asymptotically cylindrical $G_2$-manifolds which have ends modelled on $\R \times S^1 \times S^1 \times S$, where $S$ is a K3 surface. For topological reasons, which will be explained later, we need to identify the two K3 surfaces of either end by a hyperk\"ahler rotation.

\subsubsection*{Fano Geometry}
We will now review some aspects of the geometry of Fano threefolds. More details can be found in the book by Kollar \cite{kollar2013rational} and the survey paper by Beauville \cite{beauvilleFanoThreefoldsK32002}.

\begin{dfn}
\label{Fano_dfn}
A \textbf{Fano manifold} is a compact, complex manifold $X$ with ample anticanonical bundle, meaning that a basis of $\Ho^0(X, (-K_X)^{\ot k})$ gives a well-defined embedding into $\C P^N$ for some $k \ge 1$.
\end{dfn}
Being Fano is quite a restrictive condition. In each dimension $n \ge 1$ there are only finitely many deformation types of Fano $n$-folds. We are mostly concerned with Fano threefolds, of which there are $105$ deformation types. In fact our entire discussion can be adapted to what Corti, Haskins, Nordstr\"om and Pacini \cite{cortiAsymptoticallyCylindricalCalabi2013} call \textbf{semi-Fano} manifolds, however their definition is somewhat involved, so  we will not present it here. Morally speaking, semi-Fanos are desingularisations of mildly singular Fanos. 

We will now recall some properties of (semi)-Fano manifolds that are relevant to our discussion of the twisted connected sum construction of $G_2$-manifolds, mainly following \cite{cortiAsymptoticallyCylindricalCalabi2013}. To begin, assume that $X$ is a Fano three-fold. We can then define the following pairing on $\Ho^2(X, \Z)$:
\eas
\langle \cdot, \cdot \rangle_X: \Ho^2(X, \Z) \times \Ho^2(X, \Z) \longra \Ho^6(X, \Z) \simeq \Z, \quad (a,b)\longmapsto a\cdot b \cdot c_1(-K_X).
\eas
This endows $\Ho^2(X, \Z)$ with a non-degenerate lattice structure. We can write 
\eas
\langle -K_X, -K_X \rangle = 2g-2 > 0,
\eas
where $g$ is the degree of the Fano three-fold. Next, let $S \subset X$ be an anticanonical divisor. It is known that generically this will be a smooth K3 surface \cite{sokurovSmoothnessGeneralAnticanonical1980}. From now on assume that is is. It can then be proven that the restriction map $\Ho^2(X,\Z)\ra \Ho^2(S, \Z)$ is a primitive embedding of lattices, where we consider $\Ho^2(S, \Z)$ with the usual intersection pairing. Thus $S$ is a $(\Ho^2(X,\Z), -K_X)$-polarised K3 surface. From this it is natural to discuss the moduli space of pairs $(X,S)$ where $X$ is a (semi)-Fano threefold and $S \subset X$ is a smooth, anticanonical K3 divisor, where $\Ho^2(X,\Z)$ is isomorphic to a fixed lattice $N$ and $A \in N$ satisfies $A^2 = -K_X^3$. Write this moduli space as $\cF^{N,A}$. This is again a (potentially singular) complex manifold. Of course we have a forgetful morphism:
\eas
s^{N,A}: \cF^{N,A} \longra \cK^{N,A}, \quad (X,S) \longmapsto S.
\eas 
It has the following important property.
\begin{prop}[Thm. 6.8 in \cite{cortiAsymptoticallyCylindricalCalabi2013}]
\label{Fano_K3_surjective}
The image of each connected component of $\cF^{N,A} $ is an open dense subset of $\cK^{N,A}$, and for smooth points $(X, S) \in \cF^{N,A}, S \in \cK^{N,A}$ we have that $s^{N,A}$ is locally a submersion.
\end{prop}
\subsection{\texorpdfstring{$G_2$ and coassociative Geometry}{G2 and coassociative geometry}}
We now discuss the basic definitions of $G_2$-geometry. The interested reader can find accessible introductions to the topic in the books \cite{karigiannisLecturesSurveysG2Manifolds2020} and \cite{joyceRiemannianHolonomyGroups2007}.

Consider $\C^3$ with the standard Calabi--Yau structure $(\C^3, J_0, \om_0, g_0, \Om_0)$. We can define the following three-form, called the \textbf{associative form} on $\R^7 = \R \times \C^3$:
\eas
\phi_0 = \dt \wedge \om_0 + \Re \Om_0.
\eas
Here $t$ denotes the coordinate on $\R$. The stabiliser of this form in $\GL(7)$ is the $14$-dimensional simple Lie group $G_2 \subset \SO(7)$. A $7$-manifold $M$ together with a three-form $\phi \in \Om^3(M)$ such that at each point $(T_pM, \phi_p)$ is isomorphic to the standard model $(\R^7, \phi_0)$ is called a $G_2$-\textbf{manifold}. The \textbf{associative form} $\phi$ induces a metric $g_\phi$ on $M$ via the pullback of the standard metric on $\R^7$. If now $\phi$ is both closed and co-closed, i.e. $\d \phi = 0$ and $\d^\star_\phi \phi =0$, then both $\phi$  and $\star_\phi \phi$ are calibrations. Their calibrated submanifolds are called associatives and coassociatives respectively. This is called the \textbf{torsion-free} case. We then also have that the holonomy of $(M, g_\phi)$ is contained in $G_2$.
\begin{ex}
\label{example_g2_mflds}
Let $(X^6, J,\om,g,\Om)$ be a Calabi--Yau threefold. Consider $M^7 = X \times S^1$ with the coassociative form $\phi =  \d s \wedge \om+ \Re \Om $, where $s$ is the coordinate on $S^1$. This $G_2$-structure is torsion-free, and a special Lagrangian $L \subset X$ gives rise to an associative manifold $L \times \{p\}$ for any $p\in S^1$, whereas a complex surfaces $S^4 \subset X$ gives rise to a coassociative submanifold $S \times \{p\}$. 
\end{ex}

\subsection{\texorpdfstring{$\Spin(7)$ and Cayley Geometry}{Spin(7) and Cayley Geometry}}
In this section we briefly review the results in $\Spin(7)$ and Cayley geometry that are most important for our discussion, based on \cite{englebertConicallySingularCayley2023a} and \cite{englebertConicallySingularCayley2023}. A thorough exposition of $\Spin(7)$ and Cayley geometry is given in the foundational papers \cite{HarvLaws} and \cite{mcleanDeformationsCalibratedSubmanifolds1998}. 

The group $\Spin(7)$ is usually seen as the double cover of $\SO(7)$. In this paper we define it as the stabilizer of the \textbf{Cayley four form}:
\ea
\label{2_1_cayley_form}
\Phi_0 &= \d x_{1234} - \d x_{1256} - \d x_{1278} - \d x_{1357} + \d x_{1368} - \d x_{1458} - \d x_{1467} \nonumber \\
&- \d x_{2358} - \d x_{2367} + \d x_{2457} - \d x_{2468} - \d x_{3456} - \d x_{3478} + \d x_{5678}
\ea 
under the action of $\GL(8, \R)$ on forms in $\La^4 \R^8$. A pair $(M, \Phi)$, where $M$ is an $8$-dimensional oriented manifold and $\Phi \in \Om^4 (M)$, is a $\Spin(7)$-manifold if at each point $p \in M$ there is an oriented isomorphism $T_pM \simeq \R^8$ that takes $\Phi_p$ to $\Phi_0$. Expressed differently, Cayley forms on $M$ are sections of a subbundle $\cA(M) \subset \La^4 T^*M$. Any Cayley form $\Phi$ induces a Riemannian metric $g_\Phi$ on $M$ via pointwise pullback of the standard metric on $\R^8$, as $\Spin(7) \subset \SO(8)$. If $\d \Phi = 0$ we say that the $\Spin(7)$-manifold is \textbf{torsion-free}. In this case the holonomy $\Hol(g_\Phi)\subset \Spin(7)$.
\begin{ex}
\label{examples_spin_7_mflds}
Examples of $\Spin(7)$-manifolds arise from other calibrated geometries, but these examples never have full holonomy. We can think of $\Spin(7)$-geometry as the generalisation of other calibrated geometries of equal or lower dimension.
\begin{itemize}
\item Let $(X^8, J, \om, g, \Om)$ be a Calabi--Yau fourfold. Then the form $\Phi = \Re \Om + \ha \om \wedge \om$ is a closed Cayley form and $g_\Phi = g$. The holonomy of $(X,g)$ will be included in $\SU(4)\subset \Spin(7)$.
\item Let $(M^7, \phi)$ be a $G_2$-manifold, not necessarily torsion-free. We can then look at $X = M \times S^1$ with the form $\Phi = \dt \wedge \phi + \star \phi$. Then $(X, \Phi)$ is a $\Spin(7)$-manifold. If $\phi$ is torsion-free then $\Phi$ will be torsion-free as well, and the holonomy will be a subgroup of $G_2 \subset \Spin(7)$.
\end{itemize}
\end{ex}

Let now $N^4 \subset M$ be any immersed submanifold. The Cayley form $\Phi$ then satisfies the Cayley inequality:
\e
\Phi|_N \le \dvol_{N}.
\e
Here $\dvol_N$ is the volume form induced by the metric $g_{\Phi}$. We say that a manifold $N$ is \textbf{Cayley} if $\Phi|_N = \dvol_N$, and that it is $\al$-\textbf{Cayley} for ($\al \in (0,1)$) if instead  $\Phi|_N \ge  \al\dvol_{N}$ (see also section 3 of \cite{englebertConicallySingularCayley2023a}). If we know that a submanifold is $\al$-Cayley for $\al$ close to $1$, then there is hope that a true Cayley submanifold can be found in its vicinity, as we will see in Theorem \ref{2_2_gluing} for a special case.

\begin{ex}
\label{examples_cayley_mflds}
The examples of 
\begin{itemize}
\item Let $(X^8, J, \om, g, \Om)$ be a Calabi--Yau fourfold, seen as a $\Spin(7)$-manifold with Cayley form $\Phi = \Re \Om + \ha \om \wedge \om$. Then both special Lagrangians (calibrated by $\Re \Om$) and complex surfaces (calibrated by $\ha \om \wedge \om$) are Cayleys. 
\item Let $X = M \times S^1$ with the Cayley form $\Phi = \dt \wedge \phi + \star \phi$, where $(M,\phi)$ is a $G_2$-manifold. For an associative $A^3 \subset M$ (i.e. a manifold calibrated by $\phi$) the product $A \times S^1$ is a Cayley in $X$. Similarly, for $C^4 \subset M$ a coassociative, the product $C\times \{p\}$ is a Cayley for any $p \in S^1$. 
\end{itemize}
\end{ex}
Looking at $(\R^8, \Phi_0)$ again, the Cayley four planes form a $12$ dimensional subset of the $16$ dimensional Grassmannian of oriented four planes. Thus the Cayley condition can be described by four independent equations, given explicitly by the vanishing of a four form $\tau$ derived from the Cayley form. It turns out that for $\al$ sufficiently close to $1$ we define we can define a \textbf{deformation operator} associated to any $\al$-Cayley:

\ea
\label{2_1_deformation_op}
\begin{array}{rl}
F: C^\infty (N,V) & \longrightarrow C^\infty(E_{\cay}), \\
v &\longmapsto \pi_E(\star_N \exp^*_v(\tau|_{N_v})) .
\end{array}
\ea
Here $V \subset \nu(N)$ is an open  neighbourhood of the zero section in the normal bundle of $N$, $E_{\cay}$ is a bundle dependent on $N$ and $\exp_v: N \hookra M$ denotes the exponential of the vector field $v \in \nu(N)$. We denote the image of this map by $N_v$. The precise definition can be found in \cite[Section 3.1]{englebertConicallySingularCayley2023a}, where it is also shown that this operator is elliptic at the zero section and detects Cayley submanifolds, meaning that $F(v) =0$ implies that $N_v$ is indeed Cayley. 

The linearisation of $F$ at the zero section is the \textbf{linearised deformation operator} $\CayD$ or simply the \textbf{Cayley operator}. The operator $F$ admits the following Taylor expansion around the zero section:
\eas
\label{2_1_taylor}
F(v) = F(0) + \CayD[v] + Q(v).
\eas
Here $Q$ is a remainder term, which contains all the quadratic and higher order behaviour of $F$ near $0$. The coefficients of $\CayD$ depend on the data $M,N$ and $\Phi$ in a very precise way.
\begin{prop}
\label{2_1_deformation_op_dependence}
Let $(M, \Phi)$ be a $\Spin(7)$-manifold, where $\Phi$ is part of a smooth family $\{\Phi_s\}_{s\in \cS}$. Then there is an open subset $U \subset \cA(M)$ be an open subset containing the image of $\Phi$, $\al$ sufficiently close to $1$ and smooth bundle maps as follows:
\eas
c_1:& M \times \Cay_\al(4,TM) \times U \longra \Hom(T^*M \ot TM, \La^2 M),  \\
c_0:& M \times \Cay_\al(4,TM) \times (T^*M \ot U)\longra \Hom(TM, \La^2 M).
\eas 
For any immersed $\al$-Cayley $N \subset (M, \Phi_s)$ with an associated linearised deformation operator $\CayD_{N,s}$, we have that:  
\eas
\CayD_{N,s} v (p) = c_1(p, T_pN, \Phi_s(p)) \cdot \nabla v (p) + c_0(p, T_pN, \nabla\Phi_s(p)) \cdot  v (p).
\eas
Here $\nabla$ refers to the Levi-Civita connection for the fixed $\Spin(7)$-structure $\Phi$.
\end{prop}
\begin{proof}
This is a consequence of  \cite[Prop. 3.4]{englebertConicallySingularCayley2023a}, which gives a coordinate expression of $\CayD$ in a carefully chosen frame $\{f_j\}_{1\le j\le 8}$ as:
\eas
\begin{array}{rl}
\CayD[v] =&\pi_E(\beta\sum_{i=1}^4 f_i \times \nabla^\perp_{f_i} v + \sum_{i=1}^4\sum_{j=1}^8 \beta_{ij} f_{j} \times \nabla^\perp_{f_i} v + \nabla_v \tau (f_1, f_2, f_3, f_4)).
\end{array}
\eas
Here $\beta, \beta_{ij}$ depend algebraically on the choice of frame (which depends on $TN$) and $\Phi_s(p)$, and $\nabla^\perp$ is the connection on the normal bundle induced by $\Phi_s$. The product $\times$ that appears also depends pointwise on $\Phi_s$, and the derivative of the form $\tau$ depends pointwise on $\nabla \Phi$ and $TN$. We remark that the Christoffel symbols of $\nabla^\perp$ also depends on $\nabla \Phi$ and that this is included in $c_0$.
\end{proof}
In particular, if two almost Cayley submanifolds are sufficiently close to one another, their deformation operators will differ in a controlled manner.
\begin{cor}
\label{2_1_deformation_op_variation}
Let $N \subset (M, \Phi)$ be an almost Cayley with linearised deformation operator $\CayD_N$. Let $v \in C^\infty(\nu(N))$ be a sufficiently small normal vector field, so that $N_v$ again admits a deformation operator. Identify the normal bundles of $N$ and $N_v$ via orthogonal projection. We can then write: 
\eas
\CayD_{N_v} = \CayD_{N} + \tilde{\CayD}_v,
\eas
where $\tilde{\CayD}_v[w] = a_1(v,N) \cdot \nabla w + a_0(v,N) \cdot w$, and: 
\eas
\md{\nabla^k a_i} \lesssim \md{\nabla^{k+1} v}. 
\eas
\end{cor}
\begin{proof}
This follows from the previous Proposition \ref{2_1_deformation_op_dependence} by realising that the variation in $T_pN$ is governed by the first derivative of $v$, and similarly for higher derivatives. Finally we note that the $c_i$ from the previous proposition only depend on the ambient $\Spin(7)$-structure and not on the submanifold.
\end{proof}

Let $N$ be an almost Cayley submanifold of $(M, \Phi)$ almost $\Spin(7)$. Let $\{\Phi_s\}_{s \in \cS} \subset \mathscr{A}(M)$ be a smooth finite dimensional family of $\Spin(7)$ structures such that $\Phi = \Phi_{s_0}$ for some $s_0 \in \cS$. 
We have a good local theory of the moduli space $\cM(N, \cS)$ of compact Cayley submanifolds isotopic to $N$ for any of the $\Spin(7)$-structures in $\cS$. Combining \cite[Thm 4.9, Thm 4.10]{englebertConicallySingularCayley2023a} we obtain the following extension of the work of McLean \cite{mcleanDeformationsCalibratedSubmanifolds1998} and Moore \cite{mooreCayleyDeformationsCompact2019}:
\begin{prop}
\label{2_1_cpt_cayley_linearisisation}
The map $F$ from \eqref{2_1_deformation_op} is smooth as a map $L^p_{k+1} \ra L^p_k$ and elliptic at the zero section. A neighbourhood of $(N, \Phi)$ in $\cM(N, \cS)$ is homeomorphic to $F^{-1}(0)$. Furthermore, if $\Coker D = \{0\}$ we say that $N$ is \textbf{unobstructed}, and $\cM(N, \cS)$ is then a smooth manifold of dimension $\frac{1}{2}(\si(N) + \chi(N)) - [N] \cdot [N] + \dim \cS$ near $(N, \Phi)$.
\end{prop}

\subsection{Manifolds with ends}
We now briefly recall the definitions of asymptotically conical and conically singular manifolds. More details can be found in \cite[Section 2.3]{englebertConicallySingularCayley2023a}. Recall that an \textbf{asymptotically conical} manifold of rate $\eta < 1$  ($\AC_\eta$) is a Riemannian manifold $(M,g)$ such that away from a compact subset $K \subset M$ we can identify $M \setminus K \simeq (r_0, \infty) \times L$ and the metric satisfies: 
\e
\md{\nabla^i(g-g_{\con})} = O(r^{\eta- 1-i}) \text{ as } r \ra \infty,\label{2_1_ac_conv}
\e
where $g_{\con} = \d r^2 +r^2h$ is a conical metric on $(r_0, \infty) \times L$. For a embedded submanifold $f: A \hookra (\R^n, g)$, where $g$ asymptotically conical of rate $\eta < 1$ and asymptotic to flat $\R^n$, we say that it is an \textbf{asymptotically conical submanifold} of rate $\eta < \la < 1$  if: 
\ea
\md{\nabla^i(f(r,p)-\io(r,p))} \in O(r^{\la-i}), \text { as } r \ra \infty. \label{2_3_extrinsic_ac}
\ea
Here $\io: C \hookra \R^n$ is the embedding of the uniquely determined asymptotic cone of $A$. 

Finally, we say that a continuously embedded topological space $N \subset (M, \Phi)$ is \textbf{conically singular} with rates $\bar{\mu} = (\mu_1, \dots, \mu_l)$, where $1 <\mu_j<2$ ($\CS_{\bar{\mu}}$) asymptotic to cones $C_1, C_2, \cdots, C_l$ if it is a smoothly embedded manifold away from $l$ points $\{z_1, \dots, z_l\}$, and there are parametrisations $\Th_j: (0, R_0) \times L_j \ra M$ of $N$ near $z_j$ such that: 
\e
\md{\nabla^i(\Th_j(r,p)-\io_j(r,p))} \in O(r^{\mu_j-i}), \text { as } r \ra 0.
\e
Here $\io_j$ is the embedding of the cone $C_i$ via a parametrisation $\chi_j$ of $M$ that is compatible with the $\Spin(7)$-structure, i.e. $\D \chi_j^* \Phi(z_j) = \Phi_0$.  In both the $\AC$ and $\CS$ case, if a submanifold is Cayley, then it must be asymptotic to Cayley cones.

We now recall briefly the theory of Sobolev spaces as well as the Fredholm theory on conical manifolds. First of all, let $M$ be an $n$-dimensional manifolds with $l$ conical ends and $E$ a bundle of tensors on $M$. For a collection of weight $\bar{\de} \in \R^l$  and a section $s \in C^\infty_c(E)$ we define the Sobolev norm: 
\e
\nm{s}_{p, k, \bar{\de}} = \left(\sum_{i = 0}^k \int_{M} \md{\nabla^i s\rho^{-w+i}}^p \rho^{-n} \d \mu\right)^\frac{1}{p},
\e
Here $w$ is a weight function that interpolates between the conical ends. On the $j$-th conical end it is given by $\de_i$. We denote the completion of $C^\infty_c(E)$ under this norm by $L^p_{k, \bar{\de}}$. Note that the Banach space structure is independent of the choice of weight function. If there is only one end, we denote the space by $L^p_{k, \de}$. The $C^k_{\bar{\de}}$ spaces are defined as the completion with regards to the norm: 
\e
\nm{s}_{C^k_{\bar{\de}}} = \sum_{i = 0}^k \md{\nabla^i s\rho^{-w+i}}.
\e
There is a Sobolev embedding theorem: 

\begin{thm}[{\cite[Thm 4.8]{lockhartFredholmHodgeLiouville1987}}]
\label{2_3_Sobolev_embedding_acyl_ac}
Suppose that the following hold: 
\begin{itemize}
\item[i)] $k - \tilde{k} \ge n\left(\frac{1}{p}-\frac{1}{\tilde{p}}\right)$ and either:
\item[ii)] $ 1 < p \le \tilde{p} < \infty$ and $\tilde{\de} \ge \de$ \text{(}$\AC$\text{)} or $\tilde{\de} \le \de$ ($\CS$)
\item[ii')]$ 1 <  \tilde{p} < p < \infty$ and $\tilde{\de} > \de$ \text{(}$\AC$\text{)} or $\tilde{\de} < \de$ ($\CS$)
\end{itemize}
Then there is a continuous embedding: 
\e
L^p_{k, \de}(E) \longra L^{\tilde{p}}_{\tilde{k}, \tilde{\de}}(E) 
\e
\end{thm}
We call a linear differential operator $D$ that is asymptotically compatible with the conical structure a \textbf{conical operator}. This means that the coefficients of $D$ tend towards the coefficients of a rescaling invariant operator on the asymptotic cone. For a precise definition we refer to \cite[Section 2.3.3]{englebertConicallySingularCayley2023a}. It is of degree $\nu \in \R$ if $D [ C^{\infty}_{\de}] \subset C^{\infty}_{\de-\nu}$. Such operators define continuous operators between conical Banach spaces, and enjoy a good Fredholm theory.

\begin{thm}
\label{2_1_change_of_index}
Let $D$ be a conical operator of order $r\ge 0$ and degree $\nu \in \R$. Let $1<p<\infty$ and $k\ge 0$. Then it defines a continuous map: 
\eas
D: L^p_{k+r, \de}(E) \longra L^p_{k, \de-\nu}(E).
\eas
If $D$ is elliptic, then this map is Fredholm for $\de$ in the complement of a discrete subset $\cD \subset \R$. This subset is determined by en eigenvalue problem on the asymptotic link.

Denote by $i_\de(P)$ for $\de \in \R \setminus \cD$ the index of the operator $D: L^p_{k+r, \de}(E) \ra L^p_{k, \de-\nu}(E)$. We then have that: $i_{\de_2}(P)-i_{\de_2}(P) = \sum_{\la \in (\de_1, \de_2) \cap \cD} d(\la)$. Here $ d(\la)$ is determined solely by the asymptotic link. 
\end{thm}
\begin{ex}
\label{2_2_ex_weights}
By \cite[Ex 2.25]{englebertConicallySingularCayley2023a} we know that the critical weights of the quadratic complex cone $C_q = \{x^2 + y^2 + z^2 = 0, w = 0\} \subset \C^4$ with link $L \simeq \SU(2)/\Z_2$, the weight between $(-2, 2)$ are: 
\e
d(-1) = 2,\quad d(0) = 8, \quad d(1) = 22,\quad d(-1+\sqrt{5}) = 6.
\e
\end{ex}

Now that we have a well-behaved theory of elliptic operators on manifolds with conical ends, we can write down explicitly the moduli spaces of $\AC_\la$ and $\CS_{\bar{\mu}}$ as zeros of non-linear partial differential operators. First, we consider the case of an $\AC_\la$ almost Cayley submanifold $A \subset (\R^8, \Phi)$. The non-linear deformation operator defined in \ref{2_1_deformation_op}  can be extended to maps between weighted spaces as follows: 
\eas
F_\AC: L^p_{k+1, \la}(V) \times  \cS & \longrightarrow L^p_{k, \la-1}(E_{\cay}).
\eas 
Here $\cS$ is a family of $\Spin(7)$-structures, with $F(v, \tilde{\Phi})$ meaning that $A_v$ is Cayley with regard to the $\Spin(7)$-structure $\tilde{\Phi}$, which may be different from $\Phi$. This map is smooth, and its linearisation at the zero section $\CayD_\AC$ is Fredholm away from a discrete critical set $\cD$ \cite[Thm. 4.16]{englebertConicallySingularCayley2023a}. For instance consider the $\AC_{-1}$ Cayley $A_\ep = \{x^2 + y^2 + z^2 = \eps, w = 0\} \subset \C^4$. From Example \ref{2_2_ex_weights} we see that between rates $-1 < \la < 0$ the operator Fredholm and of constant index, which turns out to be $2$ (cf. \cite[Rm. 4.18]{englebertConicallySingularCayley2023a}). Furthermore it is unobstructed. Hence the nonlinear \textbf{moduli space} of $\AC_{\la}$ Cayley submanifolds isotopic to $A_\eps$ and having the same limiting cone $C_q = \{x^2 + y^2 + z^2 = 0, w = 0\}$ is given by: 
\eas
\cM_\AC^{\la}(A_\eps) = F_\AC^{-1} (0) \simeq \C\setminus \{0\},
\eas
where we identify $A_\eps \ra \eps$. We also define the \textbf{extended moduli space} to be the usual space where we also adjoin the singular cone as an element. So in this case we have: 
\eas
\brM_\AC^{\la}(A_\eps) \simeq \C\setminus \{0\} \cup \{C_q\} \simeq \C.
\eas
On this moduli space we the \textbf{scale} of an element to be its rescaling factor with regards to a fixed cross section of the scaling action. On $\brM_\AC^{\la}(A_\eps)$ such a scale could be: 
\eas
t(A_{\eps}) = \md{\eps}.
\eas
In this case the chosen cross-section of scale one elements is $\{A_\eps: \md{\eps} = 1\}$. All of this is explored in more detail in \cite[Section 4.2]{englebertConicallySingularCayley2023a}. 

Next, there are similar moduli spaces of conically singular Cayley submanifolds. First, let $N \subset (M, \Phi)$ be a $\CS_{\bar{\mu}}$ Cayley for the $\Spin(7)$-structure $\Phi$, which is part of a smooth family $\{\Phi_s\}_{s \in \cS}$. Then we can define the moduli space of $\CS_{\bar{\mu}}$ Cayley manifolds that are $\CS_{\bar{\mu}}$ perturbations of $N$, $\cM_\AC^{\bar{\mu}}(N, \cS)$, to be the zero set of the non-linear operator: 
\eas
F_{\CS, \text{fix}}: L^p_{k+1, \bar{\mu}}(V) \times \cS & \longrightarrow L^p_{k, \bar{\mu}-1}(E_{\cay}).
\eas 
In particular they all have the same singular points and limiting cones, since the  condition $\mu_i > 1$ excludes translations of the singular points, of rate $0$, and deformations of the cone, which are of rate $1$. We can reinsert them by hand by defining two new operators: 
\ea
\label{2_2_cs_deformation_op}
F_{\CS, \text{points}}: L^p_{k+1, \bar{\mu}}(V) \times \cS \times \cF_{\text{points}} & \longrightarrow L^p_{k, \bar{\mu}-1}(E_{\cay}), \\
F_{\CS, \text{cones}}: L^p_{k+1, \bar{\mu}}(V) \times \cS \times \cF_{\text{cones}} & \longrightarrow L^p_{k, \bar{\mu}-1}(E_{\cay}).
\ea 
Here $\cF_{\text{points}}$ and $\cF_{\text{cones}}$ contain the possible perturbations of the singular points and the cones (with fixed points) respectively. Hence, in the case of a unique singular point, if $F_{\CS, \text{cones}}(v, \Phi, \tilde{C}) = 0$, this means that there is a conically singular Cayley, $\CS_\mu$ to the new cone $\tilde{C}$, but with the same singular point as $N$. If instead $F_{\CS, \text{cones}}(v, \Phi, \tilde{p}, \tilde{C}) = 0$, then the Cayley has a possibly different singular point $\tilde{p}$, and a singular cone $\tilde{C} \subset T_{\tilde{p}}M$. We denote the corresponding moduli spaces by $\cM_\AC^{\bar{\mu}}(N, \cS) = \cM_\AC^{\bar{\mu}, \text{points}}(N, \cS)$ and $\cM_\AC^{\bar{\mu}, \text{cones}}(N, \cS)$. Similar to the asymptotically conical case, the operators $F_{\CS, (\star)}$ define smooth maps and have a good Fredholm theory \cite[Thm. 4.32]{englebertConicallySingularCayley2023a}. We conclude our discussion of conical manifolds by proving a useful lemma. Since we now how $\CayD$ varies under perturbation by a vector field by Corollary \ref{2_1_deformation_op_variation}, we can determine precisely what the convergence rate of $\CayD$ to the conical model is depending on the rate of the $\AC$ or $\CS$ manifold.

\begin{lem}
\label{2_2_deformation_op_conical} 
Let $A \subset (\R^8, \Phi)$ be an almost Cayley submanifold which can be seen as a perturbation of a Cayley cone $C \subset \R^8$ by a normal vector field $v \in C^\infty(\nu(C))$ with $\nm{v}_{C^{k+1}_\ga} = 1$ for $\ga \in \R$. We identify the tensor bundles on $A$ and $C$ so that the Cayley operator $\CayD$ of $A$ and the Cayley operator of the cone $\CayD_{\con}$ are both defined on $C$. For any rate $\ze \in \R$ we then have the pointwise estimate:
\eas
\md{(\CayD-\CayD_{\con})s}_{C^{k}_{\ze-1} } \lesssim r^{\ga-1}\md{s}_{C^{k+1}_{\ze}}.
\eas 
\end{lem}

\subsection{Desingularisation of Cayley submanifolds}

We now have all the necessary tools at our disposal to discuss a desingularisation theorem that will allow us to study fibrations in a neighbourhood of their singular fibres. This is a slight modification of \cite[Thm. 3.15]{englebertConicallySingularCayley2023}, in that we separate the $\AC$ and $\CS$ deformations at a different rate. Here we exclusively consider cones which have no critical rates between the translations at $0$ and the rotations at $1$, which we call \textbf{semi-stable} cones. 

\begin{thm}[Gluing Theorem]
\label{2_2_gluing}
Let $(M, \Phi)$ be an almost $\Spin(7)$ manifold and $N$ a $\CS_{\bar{\mu}}$-Cayley in $(M, \Phi)$ with singular points $\{z_i\}_{i = 0, \dots, l}$ and rates $1 < \mu_i < 2$, modelled on the semi-stable cones $C_i = \R_+ \times L_i \subset \R^8$ . Assume that $N$ is unobstructed in $\cM_{\CS}^{\bar{\mu}} (N, \{\Phi\})$. For a fixed $k\le l$, assume for each $i \le k$ that the $L_i$ are unobstructed as associatives (i.e. that the $C_i$ are unobstructed cones), and that $\cD_{L_i} \cap (1, \mu_i] = \emptyset$. For $1 \le i \le k$, suppose that $A_i$ is an unobstructed $\AC_\la$-Cayley with $\la < 0$, such that $\cD_{L_i} \cap [\la, 0) = \emptyset$. Let $\{\Phi_s\}_{s \in \cS}$ be a smooth family of deformations of $\Phi = \Phi_{s_0}$. Then there are open neighbourhoods $U_i$ of $C_i \in  \overline{\cM}_{\AC}^{\la}(A_i)$, an open neighbourhood $s_0 \in U \subset \cS$  and a continuous map:
\e
\Ga: U \times \cM_{\CS}^{\bar{\mu}} (N, \{\Phi\}) \times \prod_{i=1}^k U_i \longra \bigcup_{I \subset \{1, \dots, k\}} \cM_{\CS}^{\bar{\mu}_I}(N_I , \cS).
\e
Here we denote by ${\bar{\mu}_I}$ the subsequence, where we removed the $i$-th element for $i\in I$ from ${\bar{\mu}}$. Moreover, $N_I$ denotes the isotopy class of the manifold obtained after desingularising the points $z_i$ for $i \in I$ by a connected sum with $A_i$. 

This map is a local diffeomorphism of stratified manifolds. Thus away from the cones in $\brM_{\AC}^{\la}(A_i)$ it is a local diffeomorphism onto the non-singular Cayley submanifolds. It maps the point $(s, \tilde{N}, \tilde{A}_1, \dots, \tilde{A}_k)$  into $\cM_{\CS}^{\bar{\mu}_I}(N_I , \cS)$, where $I$ is the collection of indices for which  $\tilde{A}_i = C_i$. This corresponds to partial desingularisation.
\end{thm}

The submanifold $\Ga(\Phi_s, N, \bar{A})$ is constructed by first producing an almost Cayley submanifold $N^{s,\bar{A}}$ which is sufficiently close to being Cayley so that an iteration scheme allows to perturb it to a nearby exact Cayley. The submanifolds $N^{s,\bar{A}}$ is explicitly given as a gluing of two kinds of pieces via a partition of unity, the first being a deformation of $N$ to a $\CS_{\bar{\mu}}$-Cayley with respect to the $\Spin(7)$-structure $\Phi_s$, and the second being the asymptotically conical desingularisations of the conical points. We recall the precise error bounds on $\tau|_{N^{s, \bar{A}}}$.

\begin{prop}[Prop. 3.6 \cite{englebertConicallySingularCayley2023a}]
\label{2_2_glue_bounds}
For a sufficiently small scale $t$ of the $\bar{A}$ and for $s \in \cS$ sufficiently close to the initial $\Spin(7)$-structure, $p > 4$ and $1 < \ga <  \ga_{\max}$, $k \in \N$, we have:
\e
\label{2_2_tau_estimate}
\nm{F_{N^{s,\bar{A}}}(0,s)}_{L^p_{k, \ga-1, \bar{A}}} < C_F (t^{\nu(\ga_{\max}-\ga)}).
\e
Here $C_F > 0$ is a constant that only depends on the geometry of $(M, \Phi)$ and not on $\bar{A}$, and $0 <\nu < 1$ is a constant appearing in the gluing construction.
\end{prop}

Afterwards we perform the following iteration scheme for each of the resulting pre-gluings $N^{s,\bar{A}}$ simultaneously, one for each choice of $N$ and $\bar{A}$. We first choose what we call \textbf{pseudo-kernels} for all $N^{s,\bar{A}}$, which are approximations to the kernel of $\CayD$ that we can control a priori. Indeed for this we fix compactly supported pseudo-kernels $\ka_{\AC}, \ka_{\CS} \subset L^p_{k+1, \ga}(\nu)$ for $N$ and $\bar{A}$ respectively. For sufficiently small $t > 0$ we can then think of the supports of these kernels as being embedded in $N^{s,t\bar{A}}$, and thus we can pull back the normal vector fields onto $N^{s,t\bar{A}}$. As the Cayley operator restricted to either end of $N^{s,t\bar{A}}$ is approaching the operator on the respective conical piece, this new space of sections $\ka_t = \ka_{\AC} \op \ka_{\CS}$ will still be a good approximation of the true kernel in the following sense. If $v \perp \ka_t$ in some adapted $L^2$-norm, then if additionally $v \in \ker \CayD_{N^{s,t\bar{A}}}$ we automatically get $v = 0$. We can now describe the iteration scheme in detail. For this consider a sequence $\{v_i\}_{i \in \N} \subset L^p_{k+1,\ga, t\bar{A}}(\nu)$ with $v_0 = 0$ and:
\ea
\label{2_3_iteration_scheme}
\CayD_{t} v_{i+1} = -F_{t}(0) - Q_{t}(v_i),   \ v_i \perp^{L^2_{\de}} \ka_t, \ i \ge 0.
\ea 
We need to invert the operator $\CayD_{t}$ relative to some vector space of dimension $\ker \CayD_t$, such as the pseudo-kernel (this requires both $N$ and $\bar{A}$ to be unobstructed, in which case $\CayD_t$ is automatically surjective). To estimate the norm of $v_{i+1}$ from the norm of $v_i$ we need the following to results, which also hold without the unobstructedness assumption.
\begin{lem}[Prop 3.12 in \cite{englebertConicallySingularCayley2023}]
\label{2_3_estimates_D}
Let $\de \in \R$ be a non-critical rate for the cone of $N$. Then there is a constant $C_D > 0$ such that we have for all sufficiently small $t$ and $u \in L^p_{k+1,\de,t\bar{A}}$ with $u \perp^{L^2_{\de}} \ka_t$:
\ea
\nm{u}_{L^p_{k+1, \de,t\bar{A}}} &\le C_D \nm{\CayD_{t}u}_{L^p_{k, \de-1,t\bar{A}}}.
\ea 
\end{lem}
\begin{lem}[Prop 3.13 in \cite{englebertConicallySingularCayley2023}]
Let $\de \in \R$ be given. Then there is a constant $C_Q > 0$ such that we have for all sufficiently small $t$ and $u,v \in L^p_{k+1,\de,t\bar{A}}(\nu)$ with $u,v$ of sufficiently small norm:
\ea
\label{2_3_estimates_Q}
\nm{Q_t(u)-Q_t(v)}_{L^p_{k+1, \de,t\bar{A}}} &\le C_Q\nm{u-v}_{L^p_{k, \de-1,t\bar{A}}}(\nm{u}_{L^p_{k, \de-1,t\bar{A}}}+\nm{v}_{L^p_{k, \de-1,t\bar{A}}}).
\ea 
\end{lem}
If $C_DC_Q\nm{F_t(0)}_{L^p_{k,\ga-1,\bar{A}}}$ is sufficiently small, then the iteration scheme converges. Noteworthy is that the estimate \ref{2_3_estimates_D} and \ref{2_3_estimates_Q} are independent of the scale $t$ of the asymptotically conical manifolds. From the iteration scheme we obtain vector fields $v_{s, \bar{A}} \in C^\infty(\nu(N^{s,\bar{A}}))$ such that $\exp(v_{s, \bar{A}})$ is $\Phi_s$-Cayley for any $N$ and $\bar{A}$. Moreover we have the bound : 
\ea
\label{2_3_iteration_error_formula}
\nm{v_{s, \bar{A}}}_{L^p_{k+1, \ga, \bar{A}}} \le C_It^{\nu(\ga_{\max}-\ga)}. 
\ea
Here $C_I$ is a constant independent of the scale of $\bar{A}$ that varies continuously with $N$ and $s$. We can perform the same iteration scheme using unweighted Sobolev spaces for any almost Cayley $N \subset (M,\Phi)$ which admits a pseudo-kernel $\ka$ satisfies the conditions :
\ea
&\bullet \nm{u}_{L^p_{k+1}} \le C_D \nm{\CayD_{t}u}_{L^p_{k}} \text{ whenever }  u \perp^{L^2}\ka \nonumber \\
&\bullet \nm{Q_t(u)-Q_t(v)}_{L^p_{k+1}} \le C_Q\nm{u-v}_{L^p_{k}}(\nm{u}_{L^p_{k}}+\nm{v}_{L^p_{k}}).\label{2_3_iteration_scheme_convergence}\\
&\bullet C_D C_Q \nm{F(0)}_{L^p_{k}} < \eps \text{ for some sufficiently small } \eps(C_D,C_Q) > 0 \nonumber
\ea
In this more general case the iteration scheme takes the form:
\e
\label{2_3_iteration_scheme_general}
\CayD v_{i+1} = -F(0) - Q(v_i),   \ v_i \perp^{L^2} \ka, \ i \ge 0.
\e
Finally, the bound on $v_\infty = \lim_{i\ra \infty}v_i$ will be of the form:
\e
\label{2_3_iteration_bound_general}
\nm{v_\infty}_{L^p_{k+1}} \le C_I \nm{F(0)}_{L^p_k},
\e
where $C_I > 0$ is a fixed constant.

\subsection{Fibrations}

Let $(M, \Phi)$ be a fixed $\Spin(7)$-manifold, and assume that $N$ is a compact, unobstructed Cayley submanifold such that every element of the moduli space $ \cM(N, \Phi)$ is unobstructed. Then $ \cM(N, \Phi)$ is a smooth manifold, which in general will be non-compact. Various kind of behaviours could in principle arise, but it is expected that generically only conically singular degenerations occur. Under such genericity assumptions, we expect $\cM(N, \Phi)$ to decompose as:
\eas
 \cM(N,\Phi)= K \sqcup \bigsqcup_{i=1}^n \Ga(\{\bar{A} \in \cM^{\AC}(\bar{A_i}):  0 < t(\bar{A}) \le 1 \}, \cM^{\CS}(N_i), \Phi).
\eas
Here $K$ is a compact set of non-singular Cayley submanifolds, and the rest is given as desingularisations of conically singular Cayleys. Again we expect that generically both the conically singular and the asymptotically conical manifolds are unobstructed. We can include the conically singular Cayleys to form the \textbf{completed moduli space} :
\eas
\brM(N,\Phi)= \cM(N,\Phi) \sqcup \bigsqcup_{i=1}^n  \cM^{\CS}(N_i, \Phi).
\eas
The topology is induced from the topology on the completed moduli space of asymptotically conical manifolds. In other words if $N_k = \Ga(A_k, \hat{N}_k)$ with $A_k$ limiting to the cone $C_k$ as $k\ra \infty$ and $\hat{N}_k \ra \hat{N}$, then $N_k \ra \hat{N}$ as well in the completed moduli space. This gives $\brM(N,\Phi) $ a well-defined topology as the gluing map $\Ga$ is continuous.
In fact this space is a stratified manifold where the full-dimensional open stratum is exactly $\cM(N,\Phi) $. The lower-dimensional strata are $\cM^{\CS}(N_i, \Phi)$, and by unobstructedness they are of codimension $\dim \cM^{\AC}(\bar{A_i}) $. From this discussion it is natural to define the following concept of a Cayley fibration.
\begin{dfn}
\label{2_3_caley_fib}
A \textbf{Cayley fibration} of a compact $\Spin(7)$ manifold $(M,\Phi)$ is a homeomorphism  $\ev: \Univ (\brM(N,\Phi) ) \simeq M$, for some smooth Cayley submanifold $N$. Here $\Univ (\cM) $ is the \textbf{universal family} of a moduli space of submanifolds $\cM$. As a topological space it is the union of all $N \in \cM$ with the topology induced from the embeddings of the $N$ into the ambient manifold. Furthermore, $\ev$ is the evaluation map that sends a point in a Cayley to itself, but seen as a point of $M$.
\end{dfn}

Ideally we would like that Cayley fibrations do not contain any singular fibres at all. However this assumption is unrealistic and we likely need to admit some singularities.
It turns out however that these singular points make it difficult to prove differentiability of the fibration under smooth perturbation. This is why we now introduce a weakened version of the fibration property. Here, stability under change of the $\Spin(7)$-structure relies only on showing continuity of the fibration under perturbation. It uses the notion of $\textbf{pseudo-cycles}$ from \cite[Section 7.1]{mcduffJholomorphicCurvesQuantum1994}. They allow us to define the degree of the evaluation map $\ev: \Univ (\bar{\cM}(N,\Phi) ) \ra M$, even though the domain is not necessarily a smooth manifold. The key is that pseudo-cycles are essentially images of smooth manifolds which can be singular in codimension at most $2$. This is sufficiently small so that a the push-forward of the fundamental class can still be defined. More precisely, a pseudo-cycle from a smooth manifold $X$ of dimension $n$ to a compact smooth manifold $M$ is a smooth map $f: X \ra M$ such that the boundary of $f(X)$ is of dimension at most $n-2$. More precisely, the we define this boundary as the set of all limit points (in $M$) of sequences $f(x_k)$, such that $x_k$ does not have a convergent subsequence in $X$. In our situation we will take $X = \Univ(\cM(N,\Phi))$ and $f = \ev$, so that the boundary of $f(X)$ consists of all the points in $M$ which lie in a conically singular Cayley. We say that two pseudo-cycles $f: X \ra M$ and $g: Y \ra M$ of dimension $n$ are \textbf{bordant}, if there is a further pseudo-cycle with boundary $h: W \ra M$ of dimension $n+1$ such that the boundary of $W$ is exactly $X \sqcup Y$, and $h$ restricts to $f$ and $g$ on $X$ and $Y$ respectively. Pseudo-cycles of a given dimension $n$, taken up to bordism, form a group, which we denote $B^n(M)$. It is related to the homology of $M$ by a group morphism $ [\cdot]: B^n(M)\ra \Ho^n(M)$. In other words, each pseudo-cycle defines a homology class. More specifically when $n = \dim M$, we can define the \textbf{degree} of a pseudo-cycle $f: X \ra M$ as $\deg f = k$ where $[f] = k \cdot [M]$, $[M]$ being the fundamental class of the compact smooth manifold $M$. This corresponds to the usual definition of the degree when $X$ is a smooth manifold. We are now able to define weak Cayley fibrations.

\begin{dfn}
\label{2_3_caley_weak_fib}
A \textbf{weak Cayley fibration} of a compact $\Spin(7)$ manifold $(M,\Phi)$ is a well defined pseudo-cycle $\ev: \Univ (\brM(N,\Phi) ) \ra M$, for some smooth Cayley submanifold $N$, where $\ev$ is required to have degree $1$. Here $\ev$ is the evaluation map that sends a point in a Cayley submanifold to the corresponding point in the ambient $M$.
\end{dfn}

Note that requiring the evaluation map $\ev$ to be a pseudo-cycle puts some restrictions on the possible local models near the singular fibres. Indeed, the singular Cayleys need to be of codimension at most $2$ in $\brM(N,\Phi)$. Thus for the unobstructed case, this means that $\dim \cM^{\AC}(A,\Phi_0) \ge 2$. This is for instance satisfied for the asymptotically conical model $A_\eps = \{x^2 + y^2 + z^2 = \eps, w = 0\}$ in $\C^4$, which has $\cM_{\AC}^\la (A_\eps) \simeq \C \setminus \{0\}$.

\section{Stability of fibrations}
\subsection{Weak fibrations}
Let $N \subset (M, \Phi)$ be an almost Cayley with a pseudo-kernel $\ka \subset C^\infty(\nu(N))$ that satisfies the convergence conditions \eqref{2_3_iteration_scheme_convergence} of the iteration scheme \eqref{2_3_iteration_scheme_general}, meaning in particular that 
\eas
 C_D C_Q \nm{F_N(0)}_{L^p_{k+1}} < \epsilon
\eas 
for a fixed $\eps > 0$. This will still be true for an $L^p_{k}$-neighbourhood of submanifolds around $N$, where the pseudo-kernel is $\ka$ parallelly transported and suitably projected. Consider now a smooth family of nearby almost Cayleys $\{N_t\}_{t\in\cT}$ with pseudo-kernels that satisfy the convergence conditions. We would like to investigate the dependence of the resulting Cayleys on the starting almost Cayley. For this, note that we can recast the deformation problem on the nearby submanifold $N_t$ as a deformation problem on $N$, but where the smooth differential operator $F_N$ is perturbed smoothly $F_{t}$. In a similar fashion we have that perturbing a $\Spin(7)$-structure $\Phi$ in a family $\{\Phi_s\}_{s\in \cS}$ gives rise to a further perturbation of the differential operator to $F_{s,t}$, where we set that $F_N = F_{s_0,t_0}$. 
\begin{lem}
\label{3_1_smooth_fixed_point}
Let $N$ be a compact, non singular almost Cayley submanifold of $(M, \Phi)$, with deformation operator $F: U\subset L^p_k (\nu(N)) \longra L^p_{k-1}(E)$, where $U$ is an open neighbourhood of $0\in L^p_k(\nu(N))$. Assume that $\ka$ is a pseudo-kernel such that $(N,\ka)$ satisfies the convergence criteria for the iteration scheme \eqref{2_3_iteration_scheme_general}. Let $F_{s,t}: U \longra L^p_{k-1}(E)$ be a family of smooth perturbations for $t\in \cT, s \in \cS$ as described above. Then there is a unique element $v_{s,t} \in U$ such that $v_{s,t}\perp \ka_{s,t}$ and $F_{s,t}(v_{s,t}) = 0$, which depends smoothly on $s,t$.
\end{lem}
\begin{proof}
First of all, note that $F_{s,t}: U\rightarrow L^p_{k-1}(\nu(N))$ is a smooth family of Banach maps. Hence the convergence criteria will also be satisfied for $(F_{s,t}, \ka_{s,t})$ with slightly larger constants $C_D, C_Q$, provided that $(s,t)$ remain sufficiently close to $(s_0, t_0)$. Thus the iteration converges to a unique solution $v_{s,t} \in \ka_{s,t}^\perp$ to $F_{s,t}(v_{s,t}) = 0$. As the constants are only slightly increased in this neighbourhood, we also see by the bound \ref{2_3_iteration_bound_general} that $\nm{v_{s,t}}_{L^p_{k}} \le 2 C_I\nm{F_N(0)}_{L^p_k}$, independent of $s,t$. We now use an implicit function argument to show that $v_{s,t}$ varies smoothly in $s$ and $t$ when it exists. We first note that we can assume $\ka_{s,t}$ to be identical, by precomposing $F_{s,t}$ with a suitably chosen automorphism of $L^p_k(\nu(N))$ that varies smoothly in $s,t$. We still call the resulting maps $F_{s,t}$ and the constants $C_D$ and $C_Q$ remain unchanged. We then look at the smooth map:
\eas
A: (\ka^\perp\cap U)\times \cS \times \cT &\longra \ka^\perp\cap U\\
(v,s,t) &\longmapsto (D_{s,t}|_{\ka^\perp})^{-1}(-F_{s,t}(0)-Q_{s,t}(v)) - v.
\eas
We clearly have $A(w,s,t)=0$ exactly when $w = v_{s,t}$. To prove smoothness of $v_{s,t}$ it is thus sufficient to show that $\partial_v A(v_{s,t},s,t): \ka^\perp \longra \ka^\perp$ is an isomorphism. One can show explicitly that:
\eas
\partial_v A(v_{s,t},s,t) = (D_{s,t}|_{\ka^\perp})^{-1} \partial_v Q_{s,t}(v_{s,t}) - \id.
\eas
From the quadratic bound on $Q_{s,t}$ we see that 
\eas
\nm{\partial_v Q_{s,t}(v_{s,t})}_{op} \le 2C_Q\nm{v_{s,t}}_X \le 4C_Q\nm{F(0)}_{L^p_k}.
\eas
From the bound on $D$, we see that $\nm{(D_s|_{\ka^\perp})^{-1} \partial_v Q_s(v(s))}_{op} \le 4C_D C_Q \nm{F(0)}_{L^p_k}$. In particular, if we further reduce $C_F$ so that $ 4C_D C_Q \nm{F(0)}_{L^p_k}<1$, we can assure that $\partial_v A(v(s),s)$ is an isomorphism.
\end{proof}

The previous result shows that a collection of compact and non singular  unobstructed Cayley submanifolds vary smoothly under change of the ambient $\Spin(7)$ structure. We now need to analyse how nearly singular Cayleys are perturbed. For this, consider an unobstructed $\CS_\mu$ ($1<\mu<2$) Cayley $N \subset (M,\Phi)$ with one singular point. Assume that we have a matching $\AC_\la$ ($\la < 1$) Cayley in $\R^8$, so that the nearby non-singular Cayleys in $M$ are given as $\Ga(\Phi, N,A)$. If $\{\Phi_s\}_{s\in \cS}$ is a smooth perturbation of the $\Spin(7)$-structure, we can then look at $\Ga(\Phi_s, N,A)$. Let $v_{s} \in C^\infty(\nu(\Ga(\Phi,N,A)))$ be the normal vector field that describes the perturbation of $\Ga(\Phi,N,A)$ to $\Ga(\Phi_s,N,A)$. We claim that it can be decomposed into two contributions as follows:
\e
\label{3_1_decomp}
v_{s} = v_{\CS,s} + \tilde{v}_s. 
\e

Here $v_{\CS,s}$ gives the deformation from $N^{s_0, A}$ to $N^{s,A}$, the pre-glued almost Cayley obtained from gluing $N_s$ and $A$. This can be thought of a gluing of the perturbation vector field that takes $N$ to $N_s$ with the perturbation vector field that takes $A$ to translated and rotated $A$, which is determined by how the conical point moves between $N$ and $N_s$. The error term $\tilde{v}_s$ is then the sum of the perturbations from $\Ga(\Phi,N,A)$ to $N^{s_0, A}$ and from $N^{s,A}$ to $\Ga(\Phi_s, N,A)$. Now by our gluing theorem \ref{2_2_gluing} we know that $\nm{\tilde{v}_s}{L^p_{k+1,\de}} < C t^\al$ for some constants $\de > 1, \al > 0, C >0$. In particular, since $t^\al \ra 0$ as $t\ra 0$, we know that the dominant term is $v_{\CS,s}$. We are now ready to prove the stability result for weak Cayley fibrations.

\begin{thm}[Stability of weak fibrations]
\label{3_1_stability_weak}
Let $(M, \Phi)$ be a $\Spin(7)$-manifold that is weakly fibred by $\Univ (\brM(N,\Phi))$, and suppose that $\{\Phi_s\}_{s \in \cS}$ is a smooth family of $\Spin(7)$-structures with $\Phi = \Phi_{s_0}$. Assume that all the Cayleys in $\brM(N,\Phi)$ are unobstructed, and that the cones in the conically singular degenerations of $N$ are semi-stable and unobstructed. Then there is an open set $s_0 \in U \subset \cS$ such that $M$ is weakly fibred by $\Univ(\brM(N, \Phi_s))$ for any $s \in U$. 
\end{thm}

\begin{proof}
Note first that all the Cayleys in $\brM(N,\Phi)$ persist under a sufficiently small perturbation of the $\Spin(7)$-structure. To see this, we apply the iteration scheme from the proof of Theorem \ref{2_2_gluing} simultaneously to all the Cayleys in $\brM(N,\Phi)$ in the following way. First, we fix pseudo-kernels $\ka_{\CS}(\hat{N})$ for all the conically singular Cayleys $\hat{N} \in \brM(N,\Phi) \sm \cM(N,\Phi)$. As we assumed all the $\CS$ manifolds to be unobstructed, their moduli spaces are of strictly lower dimension than $\cM(N,\Phi)$ (as $\dim \cM^{\AC}(A) \ge 1$ by rescaling). These moduli spaces can be non-compact as well, but only in that further conical singularities can appear. Thus only finitely many conical singularities can appear, and both $\brM(N,\Phi)$ and $ \brM(N,\Phi) \sm \cM(N,\Phi) $ must be compact. In particular we can bound the values of the constants $C_D, C_Q $ uniformly for all conically singular Cayleys that appear. The same is true for the non-singular Cayleys that are a fixed distance away from the singular points, as they form a compact set as well. Finally by the estimates \ref{2_3_estimates_D} and \ref{2_3_estimates_Q} we see that the remaining non singular Cayleys, which are desingularisations of the conically singular ones have bounded $C_D$ and $C_Q$ as well. Here we note that this is exactly because we adapt our Banach spaces to the scale of the glued manifold. In conclusion the values of the constants $C_D$ and $C_Q$ are uniformly bounded for all Cayleys in the weak fibration. In particular, for sufficiently small perturbations of the $\Spin(7)$-structure that fix the singular points, we can ensure that $2C_D$ and $2C_Q$ are still valid constants, and that the initial error $C_F$ is arbitrarily small as well. Hence for small perturbations of the ambient manifold all Cayleys persist simultaneously, and furthermore we get a family of vector fields $v_s \in \Map(\Univ (\brM(N,\Phi)), TM)$ for $s \in \cS$. These vector fields need not be continuous a priori, as they are defined separately on each Cayley as the limit vector field $v_\infty$ obtained in the iteration scheme. However by Lemma \ref{3_1_smooth_fixed_point} above we immediately see that they fit together to form a smooth vector field on the open subset of $\Univ (\brM(N,\Phi))$ given by the union of all non-singular Cayleys. Similarly we see that on a singular stratum of $\Univ (\brM(N,\Phi))$ with fixed kinds of conical singularities the vector fields also fit together to form a single smooth vector field along that stratum. What is not a priori known is the regularity of the global vector field along the normal direction of a singular stratum, i.e. what happens as a non-singular Cayley degenerates towards a singular Cayley. We can now use the bounds on the desingularisations in Equation \eqref{2_3_iteration_error_formula} to show exactly that. Consider for this a conically singular Cayley $N \subset M$, so that the nearby fibres of the fibration at time $s \ge 0$ are given by its desingularisations. Assume that under the change of $\Spin(7)$-structure we map $N$ to $N_s$. We will then choose the identification $\phi: \Univ (\brM(N,\Phi)) \simeq \Univ (\brM(N,\Phi_s))$  as topological spaces so that $\Ga(\Phi, N, \bar{A})$ is sent to $\Ga(\Phi_s, N_s, \bar{A})$. We can now analyse the behaviour of the vector fields $v_s \in \Map(\Univ (\brM(N,\Phi)), TM)$ near the singular fibres. As we have seen from \eqref{3_1_decomp}, the vector field $v_{\Ga(\Phi, N, \bar{A})}$ that describes the perturbation of $\Ga(\Phi, N, \bar{A})$ decomposes as follows:
\e
v_{\Ga(\Phi, N, \bar{A})} = v_{\CS,s} + \tilde{v}_s.
\e
Here $v_{\CS,s}$ is a glued vector field, obtained by combining the vector fields $v_{N,s}$ that take $N$ to $N_s$ and $A$ to a rotated and translated $A$. In particular, this component approaches $v_{N,s}$ as $t\ra 0$. The other component,$\tilde{v}_s$, satisfies $\nm{\tilde{v}_s}_{L^p_{k+1,\de}} < C t^\al$ from our gluing theory, and hence also $\md{\tilde{v}_s}_{C^0} < t^\al\rho^\delta$. Thus $ev_{\Phi_s}$ is a continuous map, even as one approaches the singular Cayleys, and the vector fields $v_s \in \Map(\Univ (\brM(N,\Phi)), TM)$ are in fact continuous, and vary continuously with $s$. 

We showed that $\phi$ is a smooth map on the non-singular stratum, and maps the singular strata homeomorphically to singular strata with the same singularities. Since we also showed that $\ev_{\Phi_s}$ are continuous maps, we see that $\ev_{\Phi_s}|_{\Univ (\cM(N,\Phi))}$ remain pseudo-cycles, since the boundaries of $\ev_{\Phi_s}$ remain of codimension at least $2$. If we now consider the manifold $W = \Univ(\cM(N, \{\Phi_r\}_{r\in [0,s]}))$ and the evaluation map $\ev: W \ra M$, we see that it forms a bordism pseudo-cycle between $\ev_{\Phi_0}$ and $\ev_{\Phi_s}$. So in particular, if the degree of $\ev_{\Phi_0}$ was $1$, it will remain $1$ for $\ev_{\Phi_s}$.
\end{proof}

\subsection{Strong fibrations}

We showed in the previous section that weak fibrations are stable under perturbation of the $\Spin(7)$ structure. This relied on the fact that the perturbation vector fields (which describe how a given Cayley perturbs under change of the $\Spin(7)$ structure to a nearby Cayley for the new structure) remain continuous under the collapse of nearly singular Cayleys to their conically singular limits. In other words, the nearly singular Cayleys deform in a similar fashion than the singular Cayleys. This means that by perturbing the $\Spin(7)$ structure, the entire completed moduli space (including the conically singular Cayleys) varies continuously, even at the singular fibre. Proving the stability of strong fibrations means improving this result by showing that these vector fields, which are continuously differentiable, have bounded $C^1$ norm as the neck size shrinks to zero and one approaches a singular limit (as we will see later, the regions away from the singularities as well as the conically singular Cayleys themselves are easy to handle, essentially because their moduli space are compact). As a toy example, we consider a fibration of $\R^2$ by lines, which we see as the projection map:
\eas
f: &\ \ \R^2 \longra \R \\
& (r,t) \longmapsto t
\eas
\begin{figure}
\label{figure_folding_over}
\includegraphics[scale=0.45]{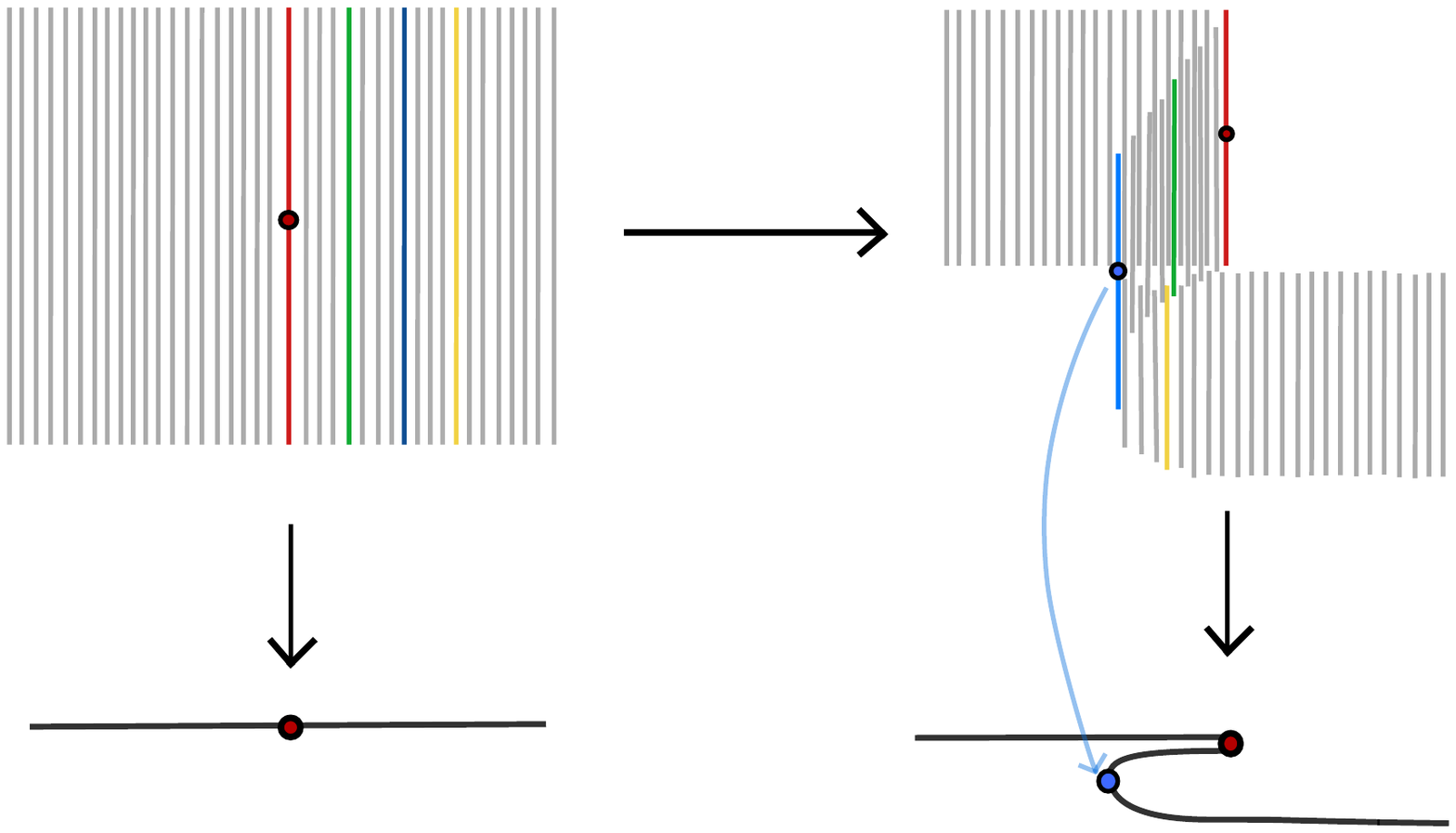}
\caption{Folding over for arbitrarily small time. The fold has width $O(s^{1/(1-\al)})$.}
\end{figure}
Here we think of $t$ as the distance from the singular fibre $f^{-1}(0)$, and of $r$ as the radial distance from the singular point $(0,0) \in \R^2$. By the weak stability result above, we know that any perturbation of this fibration gives a homotopy of the initial identity map, which we call $h: (-\eps, \eps) \times \R^2 \longra \R$. Note that here the new fibres are not the pre-images $h_s^{-1}(t)$ for $t \in \R$, but the images $h(s,t, \cdot)$ for some fixed $s,t \in \R$. Let's say that the singular fibre and point remain fixed, so that $h(s, r, 0) = 0$ for all $s \in (-\eps, \eps)$. We proved above that for $t\in \R$ sufficiently close to $0$, the value of  $h(s,r,t)$ remains close to $h(0,r,t)$ in that $\md{h(s,r,t)-h(0,r,t)} \le s\md{t}^\al \md{r}^\ga$ ($0 < \al < 1 <\ga$). We realise quickly then that this does not imply that $h(s, \cdot, \cdot)$ is $C^1$ on all of $\R^2$ for $s\ne 0$, even if we assume that it is differentiable everywhere at time $s = 0$ and differentiable away from the singular fibre for all time. Indeed, we consider:
\eas
h(s,r,t) = t - s\md{t}^\al \md{r}^\ga.
\eas 

Then clearly $\partial_t h \ra \infty$ as $t \ra 0$ for some fixed $r \ne 0$. Thus the fibres in this fibration start to move very quickly relative to one another, even though they do not perturb very much for any finite time. And indeed the fibration property is not preserved, as for $s \ne 0$ the weak fibration $h_s$ admits fibres that intersect. Just consider the initial fibres $t = \eta$ and $t = \eps$ for $0 < \eps <\eta< s^{1/(1-\al)} $. One can see that there is $r_{\eta, \eps}$ such that $h_s(\eta, r_{\eta, \eps}) = h_s(\eps, r_{\eta, \eps})$. 
Thus we need to analyse the equivalent of $\partial_t h$ for the Cayley fibration problem, which are the \textbf{infinitesimal deformation vector fields}. For a given Cayley fibration $f: M^8 \ra B^4$, they are the normal vector fields to a Cayley $N = f^{-1}(b)$ which are lifts of tangent vector $v \in T_bB$ in the base. Seen differently, they are the first order variation of a family of Cayleys parametrised by curves in $B$. Finally, they can also be seen as solutions to the linearised Cayley problem $\CayD w = 0$ for $w \in C^\infty(\nu(n))$, and this is the perspective we will use. To show boundedness of the infinitesimal deformation vector fields, we solve the linearised Cayley equation $\CayD w = 0$ on the desingularisations $N\glue{\Phi}\bar{A}$ via gluing. In the following we assume:
\begin{itemize}
\label{3_2_conditions}
\item  $N \subset (M, \Phi)$ is an unobstructed $\CS_\mu$ Cayley ($1 < \mu < 2$) with a unique singular point, with semi-stable cone $C \subset \R^8$ .
\item  $A \subset \R^8$ be an unobstructed $\AC_\la$ Cayley ($\la <0$) with asymptotic cone $C$. 
\item  There is a critical rate $\ze = \max\{\cD_C \cap (-\infty, 0)\}$ such that $\CayD_{\AC}$ is an isomorphism at rates just below $\ze$.
\end{itemize}
Under these assumptions, the deformation vector fields of $N\glue{\Phi}\bar{A}$ split into two classes:
\begin{itemize}
\item The deformations of rate $\ze$, which come from varying $A$, i.e. $w_{\AC} = \partial_t N\glue{\Phi}\bar{A_t}$ for a family $\{A_t\}_{t\in(-\eps, \eps)}$. These correspond to moving \textit{orthogonal} to the singular locus in $B$.
\item The deformations of rate $0$, which come from varying $N$, i.e. $w_{\CS} = \partial_t N_t\glue{\Phi}\bar{A}$ for a family $\{N_t\}_{t\in(-\eps, \eps)}$. These correspond to moving \textit{parallel} to the singular locus in $B$.
\end{itemize}

We will now show in turn that these infinitesimal deformation fields remain bounded in suitable weighted Sobolev spaces. 

\subsubsection*{Deformations in the normal directions}

First we look at the deformations of nearly singular Cayley submanifolds that are coming from variations in $A$, i.e. deformations of rate $r^\ze$ with $\ze  < 0$. Note that as the neck size $t \ra 0$, we can find vector fields $w_{\AC,t}$ as above with $\min \md{w_{\AC,t}} \ra 1$ but $\max \md{w_{\AC,t}} = O(t^\ze)$. In this sense they are very different from vector fields coming from parallel movement, which are of constant magnitude as we approach the singular limit.

Suppose that $\{\Phi_s\}_{s\in\cS}$ is a smooth family of $\Spin(7)$ structure on $\R^8$ which are all $\AC_\eta$ ($\eta < 1$) to the flat $\Phi_0$. For a fixed $s_0 \in \cS$, let $A \subset (\R^8, \Phi_{s_0})$ be an unobstructed $\AC_\la$ $\al$-Cayley submanifold ($\eta< \la < 1$) asymptotic to the cone $C = \R^+ \times L$. Suppose that $\al$ is sufficiently close to $1$, so that $A$ admits a linearised deformation operator $\CayD_\AC$. For a given weight $\ze \in \R$, denote by $\cI_{\AC}^{\ze}(A)$ the solutions $w\in C^\infty_\ze(\nu(A))$ to $\CayD_{\AC} w = 0$, i.e. which have decay rate at most $r^\ze$. More precisely, fix an identification of the end of $A$ with $(R_0, \infty)\times L$. We can then define, when it exists:
\e
\partial_\ze w = \lim_{r \ra \infty} r^{\ze+1}M_r^*(w|_{\{r\}\times L}).
\e
Here $M_r$ denotes the map from the unit sphere in $\R^8$ to the sphere of radius $r$, also in $\R^8$. Hence $\partial_\ze w \in C^\infty(\nu(L\subset S^7))$ is extracting the component $r^\ze \si$ of exactly rate $r^\ze$. The vector field $w$ is called $\ze$-\textbf{non-zero} if $\partial_\ze w \ne 0$. Similarly we call it $\ze$-\textbf{non-vanishing} if $\partial_\ze w$ is nowhere vanishing. We define  $\cI_{\AC}^{\ze}(A)$ to be the subspace $\ker \CayD_{\AC} \subset C^\infty_\ze(\nu(A))$. Note that on the end we can write $\CayD_{\AC} = \frac{\d}{\d r} - r^{-1}B(r)$, with $B(r) \ra B_\infty$, the limiting operator on the link. If $w \in \cI_{\AC}^{\ze}(A)$, then $\partial_\ze w $ is well-defined and an eigensection of the limiting operator $B_\infty$ with eigenvalue $\ze$. By Lemma \ref{2_2_deformation_op_conical}, the asymptotic behaviour of $B(r)$ can be given more precisely as:
\e
\label{3_2_convergence_limit_op}
\nm{B(r)-B_{\infty}}_{\text{op}} = O(r^{\la-1}).
\e
Here the operator norm is taken with regards to an arbitrary Sobolev norm on the cross section. From this we deduce the asymptotic expansion of infinitesimal deformation vector fields.
\begin{prop}
\label{3_2_asymptotic_expansion}
Let $w \in \cI^{\ze}_{\AC}(A)$ with $\CayD_{\AC}$ unobstructed at sufficiently large rates $\ze-\eps < \ze$ for $\eps > 0$. Then there is $\eps, R > 0$ such that for $r > R$ and $p\in L$ we have:
\eas
w(r,p) = (\partial_\ze w)(p) r^\ze + \de w,
\eas
where $\de w \in C^\infty_{\ze-\eps}(\nu)$.
\end{prop}
\begin{proof}
Recall that the $\AC_\la$ condition gives us an identification of the end $A\setminus K$ with $L \times (r_0, \infty)$, where $K \subset A$ is compact. Then define  $\de w = w - \al(r) \si r^\ze$ for $\si = \partial_\ze w$ and a cut-off function $\al: [r_0, \infty) \ra \R$ such that $\al = 1$ for large $r$ and $\al(r_0) = 0$. We then compute, using the fact that $\CayD_{\AC} =\frac{\d}{\d r} - r^{-1}(B_\infty+\de B(r))$ with $\nm{\de B(r)}_{\text{op}} = O(r^{\la-1})$: 
\eas
0 &= \CayD_{\AC}w  = \CayD_{\AC}( \al(r) \si r^\ze + \de w) \\
&= r^{\ze-1}(r\partial_r \al \si- \de B(r)[\al\si])+\CayD_{\AC}\de w.
\eas
In particular this implies that for $r$ large we have: 
\eas
\CayD_{\AC}\de w = \de B(r)\si r^{\ze-1} \in C^\infty_{\ze-1+(\la- 1)}(\nu),
\eas
 so that from unobstructedness at the rate $\ze + \la -1 < \tilde{\ze} < \ze$ we see that there is $\tilde{w} \in L^p_{k,\tilde{\ze}}(\nu)$ with $\CayD_{\AC} \tilde{w} = \CayD_{\AC} \de w$ (not necessarily unique).
In particular this means that $\tilde{w} = \de w + u$, where $\CayD_{\AC}u=0$. However there are no non-zero infinitesimal deformation vector fields with rate in $(\tilde{\ze}, \ze]$ which satisfy $\partial_\ze u = 0$, hence $\de w$ itself must have decay in $O(r^{\ze -\eps})$ for sufficiently small $\eps > 0$. From elliptic regularity for the operator $\CayD_{\AC}$ at rate $\ze -\eps$ we can now deduce that $\de w\in C^\infty_{\ze -\eps}(\nu)$.
\end{proof}

We can now prove that both the existence and $\ze$-nowhere-vanishing of infinitesimal deformation vector fields are stable under $\AC_\la$ perturbations with $\la < 0$.

\begin{prop}
\label{3_2_perturbation_ac}
Let $\{A_t\}_{t\in \cT}$ be a smooth family of $\AC_\la$ perturbations of $A = A_{t_0}$. Assume that all the elements of $\cI_{\AC}^{\ze}(A, \Phi)$ are $\ze$-nowhere-vanishing and that the operator $\CayD$ is an isomorphism at rate $\ze-\eps$ for some small $\eps > 0$. Then there is an open neighbourhood of $(s_0,t_0) \in \cS \times \cT$ and $R>0$ such that all the elements of $\cI_{\AC}^{=\ze}(A_t, \Phi_s)$ are still $\ze$-nowhere-vanishing.
\end{prop}
\begin{proof}
Take a solution $w \in \cI_{\AC}^{=\ze}(A, \Phi)$, which we can therefore write as 
\eas
w =  (\partial_\ze w)r^\ze \al + \tilde{w}.
\eas
Here $\al: A \ra \R$ is some cut-off function that is zero for $r \le R_1$ and one for $r \le R_2$. We have that $\tilde{w} \in L^p_{k, \ze-\eps}$ for some small $\eps > 0$ by Proposition \ref{3_2_convergence_limit_op}. We can now write the linearised Cayley equations for $(A_t, \Phi_s)$ as:
\e
\label{3_2_perturbed_AC_inf}
\left(\frac{\d}{\d r} - r^{-1}(B_{\infty}+ \tilde{B}_{t,s}(r))\right) w = 0. 
\e
Here $B_{\text{cone}}$ is the limiting operator on the link, and $\tilde{B}_{t,s}(r)$ is the error term introduced by $(A_t, \Phi_s)$. By Proposition \ref{2_1_deformation_op_dependence} it satisfies $\nm{\tilde{B}_{t,s}(r)}_{op} \in O(r^{\la-1})$, where $\la < 1$ is the asymptotic rate of $A$. In particular, if we make the ansatz $w_{t,s} =  (\partial_\ze w)r^\ze \al + \tilde{w}_{t,s}$ as above, this reduces Equation \eqref{3_2_perturbed_AC_inf} to:
\eas
\CayD_{t,s} \tilde{w} = r\partial_r \al (\partial_\ze w) + \tilde{B}_{t,s}(r)[\al\partial_\ze w]. 
\eas
Now, since $\al$ will eventually be equal to $1$ for sufficiently large radii, this means that the right hand side in the above expression will decay like $\tilde{B}_{t,s}(r)[\partial_\ze w]$, i.e. in $O(r^{\ze+(1-\la)})$ by what we just said. In particular, finding $\tilde{w}$ amounts to solving the linearised equation in $L^p_{k+1, \ze-\eps}(\nu)$, which can be done uniquely by assumption.
\end{proof}

Let $N \subset (M, \Phi)$ be an unobstructed $\CS_\mu$ Cayley ($1 < \mu < 2$) and $A \subset (\R^8, \Phi_0)$ an unobstructed $\AC_\la$ Cayley ($\la < 0$) which satisfy the assumptions of the gluing theorem \ref{2_2_gluing}. In particular they admit the same asymptotic cone $C \subset \R^8$ and the same critical rates $\cD_C \subset \R$. We thus get a family $N^{tA}$ of compact almost Cayley submanifolds of $M$ for $t>0$ small, obtained by gluing a copy of $A$ rescaled by a factor $t$ onto the conical singularity on $N$. By Proposition \ref{2_2_glue_bounds} this approximate Cayley satisfies $\nm{\tau|_{N^{tA}}}_{L^p_{k, \ga,t}} \le t^{\ga_{\max}-\ga}$, where $\ga_{\max} > 1$ is determined by the additional gluing parameter $0<\nu<1$ and the rates $\mu$ and $\la$. Furthermore, there are perturbations $N \glue{} tA$ of $N^{tA}$ which are truly Cayley and that satisfy $N \glue{} tA = \exp({v_t})$ for normal vector fields $v_t \in C^\infty(\nu(N^{tA}))$ with $\nm{v_t}_{L^p_{k, \ga,t}} \le 2t^{\ga_{\max}-\ga}$. Since we know that $N^{tA}$ is $L^p_{k, \ga,t 
}$-close to the cone in the intermediate region, we can deduce from Proposition \ref{2_1_deformation_op_variation} that the linearised deformation operator satisfies:
\e
\label{3_2_N_t_op}
\CayD = \frac{\d}{\d r} - r^{-1}(B_\infty + \de B_t(r)).
\e
Here $\nm{\de B_t(r)}_{op} = O(r^{\ga_{\max}-\ga})$ by the gluing estimates \ref{2_2_glue_bounds}.

We will now perform an additional gluing construction for the infinitesimal deformation vector fields defined on the glued manifolds $N \glue{} tA$. Let $w_{\CS} \in \cI^{\ze}_{\CS}(N)$ be an infinitesimal deformation of $N$ of rate $\ze \in \cD_C$, where $\ze = \max\{(-\infty, \la)\cap \cD\}$. This is a solution to $\CayD_{\CS}[w_{\CS}] = 0$ so that $\partial_\ze w_{\CS} = \si$ is non-zero, except if $w_{\CS}=0$. Similarly let $w_{\AC}\in \cI^{\ze}_{\AC}(A)$ be an infinitesimal deformation of $A$ of rate $\ze$ as well that shares the same limit eigensection $\si$. We then claim that a suitably glued vector field $w_t$ on $N^{tA}$ will be a good approximation of a true solution to the equation $\CayD_{N^{tA}} w = 0$. For this, we recall that $N^{tA}$ can be divided into three pieces as follows: 
\eas
N^{tA}= N^{tA}_u \sqcup N^{tA}_m \sqcup N^{tA}_l. 
\eas
Here $N^{tA}_u$, the upper region, is just $\{p \in N: \rho(p) > tr_0\}$, a truncated version of $N$. Then we have the lower region, $N^{tA}_l$, which is $\{p \in tA: \md{p} < tR_0\}$ embedded into $M$ via a $\Spin(7)$-parametrisation $\chi: \R^8 \ra M$. Finally the middle region $N^{tA}_m$ is interpolating the two pieces between the radii $t r_0$ and $R_0$. The gluing parameter $0<\nu<1$ determines where the interpolation happens, in between the radii $\ha t^\nu$ and $t^\nu$. We can however also think about this decomposition in a different way, namely: 
\eas
N^{tA}= N^{tA}_{\CS} \cup N^{tA}_{\AC}.
\eas
Here $N^{tA}_{\CS} = \{p \in N: \rho(p) > t^\nu\}$ extends the upper region from before down to radius $t^\nu$ and still agrees with $N$, but $N^{tA}_{\AC}$ now includes everything up until radius $t^{\nu''}$, where $0 <  \nu''< \nu' < \nu < 1$ are two further parameters. The gluing of the infinitesimal deformation vector fields will be performed between $t^{\nu'}$ and $t^{\nu''}$. If we now consider $t^{-1}\chi^{-1}(N^{tA}_{\AC}) \subset (\R^8, \Phi_0)$, we see that it is a non-compact almost Cayley, that agrees with $A$ for radii below $\ha t^{\nu-1}$, but extends to radius $t^{\nu'-1}$. In fact, the $\CS_{\mu}$ condition implies that this non-compact Cayley can be extended to a an $\AC_\la$ Cayley extending all the way to infinity, such that the resulting family $A_t \subset \R^8$ converges to $A$ in $C^\infty_\la$. We can do the same for the family of $\Spin(7)$-structures $t^{-1}\cdot \chi^{-1}(\Phi|_{B_r(p)})$, and they will form an $\AC_\eta$ family for $\eta < \la$ . See also Lemma 3.7 in \cite{englebertConicallySingularCayley2023a} for a more precise description of how to do this. First, since $A$ is unobstructed, the Proposition \ref{3_2_perturbation_ac} shows that we can find a smooth family of perturbations $\tilde{w}_{t,\AC} \in C^\infty_\ze(\nu(A_t))$ of $w_{\AC}$ with $\partial_\ze \tilde{w}_{t, \AC} = \partial_\ze w_{\AC}$, so that $\tilde{w}_{t,\AC} \ra w_{\AC}$ in $C^k_{\ze}$ as $t\ra 0$. Hence we also get infinitesimal deformation vector fields over $N^{tA}_{\AC}$ which we will also denote by $w_{\AC,t}$. Now choose a smooth cut-off function $f_{\mathop{cut}}: \R \ra [0,1]$, such that:
\eas
f_{\rm{cut}}|_{(-\infty, \nu'']} = 0, \quad f_{\rm{cut}}|_{[\nu', +\infty) } = 1.
\eas
Using $f_{\rm{cut}}$ we define a partition of unity on $N^{tA}$ as follows:
\e
\label{3_2_alpha_t}
\phi(p) = \left\lbrace \begin{array}{cc}
f_{\rm{cut}}\left(\frac{\log(\rho(p))}{\log(t))}\right),&\text{ if } p \in \Psi^{tA}(L \times (r_0t_i, R_0)) \\ 
0,& \text{ if } p \in  N^{tA}_u , \\
1,& \text{ if } p \in N^{tA}_l, \\
\end{array} \right.
\e
Here $\Psi^{tA}$ is a parametrisation that identifies $L\times (r_0t, R_0)$ with the gluing region of $N^{tA}$. We can now define $w_t$ as the interpolation: 
\eas
w_t = \phi w_{\CS} + (1-\phi)t^{\ze} w_{\AC,t}.
\eas
This normal vector field is a good approximation to an infinitesimal deformation on $N^{tA}$ in that it is almost a solution to the linearised Cayley equation:
\begin{prop}
\label{3_2_error_estimate}
We have for $A$, $N$, $w_t$ and $\ze$ as above that there is some $\eps, \al >0$ such that:
\eas
\nm{\CayD w_t}_{L^p_{k,\ze+\eps-1, t}} \lesssim t^\al\nm{w_t}_{L^p_{k+1,\ze, t}}.
\eas
\end{prop}
Here we note that $L^p_{k+1,\ze, t}$ is the highest weight space such that the family $w_t$ has sub-polynomial volume growth. Indeed, we have:
\eas
\nm{w_t}_{L^p_{k+1,\tilde{\ze}, t}} = \left\lbrace \begin{array}{cc}
O(1),& \tilde{\ze} <  \ze,\\
\Th(\log t),& \tilde{\ze} = \ze, \\
O(t^{\ze-\tilde{\ze}}),& \tilde{\ze} > \ze
\end{array}\right..
\eas
To see this, we note that the intermediate conical regional, which has mass proportional to $\log t$,is the dominant term. Proposition \ref{3_2_error_estimate} shows that the decay of $\CayD w_t$ is faster than the expected $O(r^{\ze-1})$ decay. In fact, it shows that the decay rate is $\ze+\eps$, for some small $\epsilon > 0$. For the proof we use the following auxiliary result comparing the Cayley operator on $N^tA$ to the conical operator on the gluing region of $w_t$. We denote this region by $G_t = (t^{\nu'}, t^{\nu''}) \times L \subset N^{tA}_m$.
\begin{lem}
\label{3_2_difference_est}
We have for $A$, $N$, $w_t$ and $\ze, \ga$ as above that there is some $\eps,\al > 0$ such that if $s \in C^\infty(\nu(N^{tA}))$ is a normal vector field then:
\eas
\md{(\CayD -\CayD_{\con})s|_{G_t}}_{L^p_{k,\ze+\eps-1, t}} \le t^\al\md{s|_{G_t}}_{L^{p}_{k+1,\ze, t}}
\eas
\end{lem}
\begin{proof}
We note that for $(r,p) \in (t^{\nu'}, t^{\nu"})\times L$ we have that the perturbation vector field $v$ taking the cone to $N^{tA}$ satisfies:
\eas
\md{\nabla^k v} \le r^{\ga_{\max}-k} .
\eas
Here $\ga_{\max}$ is the approximation rate of the gluing (as in Proposition \ref{2_2_glue_bounds}). We have $\ga_{\max}-1 >0$, and so by Corollary \ref{2_2_deformation_op_conical} we see that for $\eps = \ha(\ga_{\max}-1)$ and $\al =  \ha(\ga_{\max}-1)\nu''$ we get:
\eas
r^{-\ze-\eps+1}\md{(\CayD -\CayD_{\con})s|_{G_t}} &= r^{-\ze-\ha(\ga_{\max}-1)+1}\md{(\CayD -\CayD_{\con})s|_{G_t}} \\
&\le r^{-\ze-\ha(\ga_{\max}-1)+1} r^{\ga_{\max}-1}\md{s|_{G_t}}_{C^1_1} \\
&\le  r^{-\ze+\ha(\ga_{\max}-1)+1} \md{s|_{G_t}}_{C^1_1} \le t^\al \md{s|_{G_t}}_{C^1_{\ze}}.
\eas
This is the case $k = 0$, and the higher order cases are entirely analogous.
\end{proof}
Next we need to take a second look at the asymptotic expansion of Lemma \ref{3_2_asymptotic_expansion} for $w_{\AC}$ and adapt it to a suitable estimate on $N^{tA}$.
\begin{lem}
\label{3_2_de_w}
For $w_{\AC, t}$ and $w_{\CS}$ as above, we can write them on the gluing region $G_t$ as follows:
\eas
\tilde{w}_{\AC,t} = r^\ze \partial_\ze \tilde{w}_{\AC}+ \de \tilde{w}_{\AC,t}, \quad w_{\CS} = r^\ze \partial_\ze w_{\CS}+ \de w_{\CS}.
\eas
Denote by $\de w_{\AC,t}$ the vector fields induced by $\de \tilde{w}_{\AC,t}$. Then there is $\eps > 0$ such that $\de\tilde{w}_{\AC,t} \in C^\infty_{\ze-\eps}$ (with $C^k_{\ze-\eps}$-norms bounded uniformly in $t$) and $\de w_{\CS} \in C^\infty_{\ze+\eps}$. On $N^{tA}$ we have furthermore for some $\eps' < \eps$:
\eas
\nm{ \de w_{\AC,t}|_{G_t}}_{L^p_{k,\ze+\eps', t}} &\le t^{\eps\left(1-\nu'\right) -\nu'\eps'}\nm{w_{t}}_{L^p_{k, \ze, t}}, \\
\nm{ t^\ze\de w_{\CS}|_{G_t}}_{L^p_{k,\ze+\eps', t}} &\le t^{\nu'(\eps-\eps')}\nm{w_{t}}_{L^p_{k, \ze, t}}.  
\eas
\end{lem}
\begin{proof}
From Lemma \ref{3_2_asymptotic_expansion} it is clear that we can find $\eps > 0$ and $\de \tilde{w}_{\AC, t}$ such that:
\eas
\tilde{w}_{\AC,t} = r^\ze \partial_\ze \tilde{w}_{\AC}+ \de \tilde{w}_{\AC,t}.
\eas
The uniform boundedness (in $t$) of the $C^k_{\ze-\eps}$ norms of $\tilde{w}_{\AC,t}$ follow from the fact that the constructions in the proof of Lemma \ref{2_2_deformation_op_conical} is continuous with respect to the operator $\CayD$. Since we are working with a family $\{A_t\}_{t\in [0,\eps)}$ where $A_0 = A$, at least for sufficiently small $t$ we will have boundedness. The result in the conically singular case is completely analogous to the $\AC$ case, so that $\de w_{\CS} \in C^\infty_{\ze+\eps}$ can be constructed. Note that in the $\CS$ case, stronger decay means higher rate, whereas in the $\AC$-case stronger decay means lower rate. 
Now we move onto the bounds on $N^{tA}$. In the following we set $\vol(G_t) = \int_{G_t} r^{-4}\dvol$.
\eas
\nm{t^\ze\de w_{\AC, t}}_{L^p_{k,\ze +\eps',t}} &\lesssim \vol(G_t)t^\ze \max_{G_t} \md{\de w_{\AC,t}}_{C^k_{\ze +\eps',t}} \\
&\lesssim \vol(G_t)t^\ze \max_{G_t} \md{\de \tilde{w}_{\AC,t}(t^{-1}\cdot)}_{C^k_{\ze +\eps',t}}\\
&\le \vol(G_t)t^\ze \left(\frac{r}{t}\right)^{\ze - \eps} r^{-\ze-\eps'}\\
&\le \vol(G_t)t^{\eps\left(1-\nu'\right) -\nu'\eps'} \lesssim t^{\eps\left(1-\nu'\right) -\nu'\eps'}\nm{w_t}_{L^p_{k,\ze,t}}.
\eas
Here we used the fact that $\vol(G_t) = O(\nm{w_t}_{L^p_{k,\ze,t}})$. Note that the exponent of $t$ is positive for sufficiently small $\eps'$. Finally, the calculation for the conically singular case is more direct:
\eas
\nm{\de w_{\CS}}_{L^p_{k,\ze +\eps',t}} &\lesssim \vol(G_t) \max_{G_t} \md{\de w_{\CS}}_{C^k_{\ze +\eps',t}} \\
&\lesssim  \vol(G_t)  r^{\ze+\eps-\eps'} \lesssim t^{\nu'(\eps-\eps')} \nm{w_t}_{L^p_{k,\ze,t}}.
\eas
\end{proof}
\begin{proof}[Proof of Prop. \ref{3_2_error_estimate}]
Note first that $\CayD w_t$ is only non-zero in the gluing annulus $G_t$, since $w_t$ is interpolating between two exact solutions in this regions. From the expression:
\eas
w_t =  \phi w_{\CS} + (1-\phi)t^{\ze} w_{t,\AC}
\eas
and using the fact that $\CayD = \CayD_{\con} + \de \CayD$ with $\de \CayD$ a small perturbation (see Lemma \ref{2_2_deformation_op_conical}) we can compute the following:
\eas
\CayD w_t &= (\CayD_{\con} + \de \CayD) w_t \\
 &= \CayD_{\con}(\phi w_{\CS} + (1-\phi)t^{\ze} w_{t,\AC}) + \de \CayD w_t
\eas
Already, we see from Proposition \ref{3_2_difference_est} that for $\eps > 0$ sufficiently small there is a $\al > 0$ such that:
\eas
\nm{\de \CayD w_t}_{L^p_{k,\ze+\eps, t}} \le t^\al \nm{w_t}_{L^p_{k+1,\ze+1, t}}.
\eas
Furthermore, if we use the asymptotic expansions from Lemma \ref{3_2_de_w} and the fact that $\partial_\ze w_{\CS} = \partial_\ze \tilde{w}_{\AC, t}$ we see that:
\eas
\CayD_{\con}(\phi w_{\CS} + (1-\phi)t^{\ze} w_{t,\AC}) &= \CayD_{\con}(r^\ze(\partial_\ze w_{\CS}) \phi+ \de w_{\CS} + (1-\phi)t^{\ze} \de  w_{t,\AC}) \\
&= \CayD_{\con}(\phi \de w_{\CS} + (1-\phi)t^{\ze} \de  w_{t,\AC}),
\eas 
since $\partial_\ze w_{\CS}$ is by definition a $\ze$-eigensection, and thus $r^\ze(\partial_\ze w_{\CS})$ an infinitesimal Cayley deformation of the cone. Now, since both $ \de w_{\CS}$ and $t^{\ze} \de  w_{t,\AC}$ have $L^p_{k,\ze+\eps',t}$ norms bounded by $t^{\al'}\nm{w_t}_{L^p_{k,\ze,t}}$ for some $\al', \eps' > 0$, we see that we get the desired expression:
\eas
\nm{\CayD w_t}_{L^p_{k,\ze+\min\{\eps,\eps'\}-1, t}} \lesssim t^{\min\{\al, \al'\}}\nm{w_t}_{L^p_{k+1,\ze, t}}.
\eas
\end{proof}

So now we found a solution up to order $r^\ze$ to the linearised Cayley equation. We now solve the equation in $L^p_{k,\ze+\eps,t}$, for which we recall that the inverse of the Cayley operators $\CayD_{N^{tA}}$ (modulo kernel) have operator norms uniformly bounded in $t$ as in Lemma \ref{2_3_estimates_D}. More precisely there are subspaces $\ka_t \subset C^\infty_c(\nu(N^{tA}))$ such that for any $u \in L^p_{k,\ze+\eps,t}(\nu(N^{tA}))$ with $u \perp \ka_t$ (for a suitably chosen inner product) we have:
\e
\label{3_2_inverse_op_bound_eq}
\nm{u}_{L^p_{k,\ze+\eps,t}} \lesssim \nm{\CayD u}_{L^p_{k-1,\ze+\eps-1,t}}
\e
This relies on the fact that $\ze+\eps$ is not a critical rate, and that both $\CayD_{\AC}$ and $\CayD_{\CS}$ are unobstructed at rate $\ze+\eps$. We proved this in \cite{englebertConicallySingularCayley2023} as Proposition 3.12. In the same paper we also show that when both operators on the pieces are surjective, then so is $\CayD$ on $N^{tA}$ (Prop. 3.14 in \cite{englebertConicallySingularCayley2023}). In particular, this means that there is a unique $u_t \perp \ka_t$ such that:
\eas
\CayD u_t = \CayD w_t,
\eas
and furthermore:
\eas 
\nm{u_t }_{L^p_{k,\ze+\eps,t}} &\lesssim \nm{\CayD u_t}_{L^p_{k-1,\ze+\eps-1,t}}\\
&=  \nm{\CayD w_t}_{L^p_{k-1,\ze+\eps-1,t}} \\
&\lesssim t^\al \nm{w_t}_{L^p_{k,\ze,t}}.
\eas

Thus in particular we get a normal vector field $w_t-u_t$, which is an infinitesimal deformation vector field on $N^{tA}$ and which we understand up to second order (orders $\ze$ and $\ze+\eps$). We now need to make the leap to a deformation vector field on the true Cayley $N\glu tA$.
\begin{prop}
\label{3_2_gluing_inf_def}
Let $N \subset (M, \Phi)$ be an unobstructed $\CS_\mu$ Cayley ($1 < \mu < 2$) with a unique singular point, and assume that its cone $C \subset \R^8$ is semi-stable. Let $A \subset \R^8$ be an unobstructed $\AC_\la$ Cayley ($\la <0$) with matching cone and sufficiently small scale. Assume that the operator $\CayD_{\AC}$ is an isomorphism just below the critical rate $\ze = \max\{\cD_C \cap (-\infty, \la)\}$ and that all deformations of $A$ are of rate exactly $\ze$. We then have that for any two matching infinitesimal deformation vector field $w_{\CS} \in \cI^{\ze}_{\CS} (N)$ and $w_{\AC} \in \cI^{\ze}_{\AC} (A)$ there are gluings $w_{\CS}\glue{t}w_{\AC} \in \cI(N\glu tA)$ such that (after identifying $\nu(N^{tA})\simeq \nu(N\glu tA)$) we have:
\ea
w_{\CS}\glue{t}w_{\AC} = w_t + \de w_t.
\ea
Here $\nm{w_t}_{L^p_{k,\ze,t}} = O(\md{\log t})$ and $\nm{\de w_t}_{L^p_{k,\ze+\eps,t}} \le t^\al \nm{w_t}_{L^p_{k,\ze,t}}$, with $\al > 0$. In particular this implies that $\md{\delta w_t} \ll \md{w_t}$ as $t \ra 0$.
\end{prop}
\begin{proof}
We first note that by the gluing estimates \ref{2_2_glue_bounds} and Lemma \ref{2_2_deformation_op_conical} we have that: 
\eas
\nm{\CayD_{N^{tA}}-\CayD_{N\glu tA}}_{\ope} \lesssim t^{\ga-1},
\eas
for some $\ga > 1$. In particular for sufficiently small $t$ we also get:
\eas
\nm{(\CayD_{N^{tA}}|_{\ka^\perp})^{-1}-(\CayD_{N\glu tA}|_{\ka^\perp})^{-1}}_{\ope} \lesssim t^{\ga-1}.
\eas
Thus the same procedure as above will allow us to prove that we can perturb our infinitesimal deformation vector field $w_t-u_t$ on $N^{tA}$ to an infinitesimal deformation vector field $w_{\CS}\glue{t}w_{\AC}$ on $N\glu tA$, with a just a further $C^\infty_{\ze+\eps}$-perturbation, whose norm we can bound in exactly the same way. This concludes the proof.
\end{proof}

\subsubsection*{Deformations in the parallel directions}

We will now discuss the deformations of nearly singular Cayley submanifolds that can be interpreted as running parallel to the singular locus in the base of the fibration. Whereas in the previous section we looked at deformations of rate $\ze < 0$ coming from the $\AC$ piece, we now look at the deformations of the next higher rate $0$, which can be understood as coming from translations of the conically singular points in the $\CS$ piece. The key difference however compared to the previous section is that the Cayley operator $\CayD_{\CS}$ on $N$ is never unobstructed at rates above $0$. This is because we assume unobstructedness of $\CayD_{\CS}$ at rates slightly below $0$ where $\ind \CayD_{\CS} \le 4$ by the fibration property. However the multiplicity of the critical rate $0 \in \cD$ is at least $8$ whenever the Cayley cone is not a plane (in which case it is $4$). Thus for all the conical models we are interested in we see by Theorem \ref{2_1_change_of_index} that $\ind_\de \CayD_{\CS} < 0$ for any $\de > 0$. In this situation, the following problem appears: if $\xi \in L^p_{k, \de-1}(E)$, then there is always a unique $w \in L^p_{k+1, \de}(\nu(N^{tA}))$ such that $\CayD w = \xi$ with $w \perp^{L_2} \ker \CayD$, independent of $\de$. Now for rates $\ze < \de < 0$ we have the important estimate:
\e
\label{3_2_inverse_bound}
\nm{w}_{L^p_{k+1, \de}} \le C \nm{\xi}_{L^p_{k, \de-1}}
\e
with $C > 0$ independent of $t$. However this is not true any more  if $\de$ falls outside this range. Indeed if $\de > 0$ then $\CayD_{\CS}$ admits non-trivial obstructions, i.e. there are elements $\xi_{\ob} \in L^p_{k, \de-1}(E)$ which are not in the image of $\CayD_{\CS}$. Since these obstructions disappear once $\de$ crosses $0$, this means that $\CayD_{\CS} w_{\ob} = \xi_\ob$ for some $w_{\ob} \in C^{k+1}_0(\nu(N))$. Thus from the perspective of $w_{\ob}$, the rate of decay of $\CayD_{\CS} w_{\ob}$ is higher than expected ($-1+\de$ instead of $-1$). This is the reason why the estimate \ref{3_2_inverse_bound} cannot hold as is. Indeed, we see for $\de$ just slightly positive that the kernel of $\CayD_{\AC}$ is $(d(\ze)+d(0))$-dimensional, while $\CayD_{\CS}$ has trivial kernel. In particular from Theorem \ref{2_3_estimates_D} we see that the bound \ref{3_2_inverse_bound} does hold, but only if we have $w \perp^{L^2_{\de\pm\eps}} \ka_t$ where $\ka_t$ is $(d(\ze)+d(0))$-dimensional. Since $d(\ze)+d(0) \ge 5$ that means that asymptotically there are some elements in $L^p_{k, \de-1,t}(E)$ on the glued manifold which simply do not admit a small pre-image under $\CayD_{N^{tA}}$. Hence if we want to proceed as in the previous section we need to avoid $\CayD w_t$ having too large a component in this \say{bad sector}.
\begin{ex} 
\label{3_2_example_local}
Consider the model fibration:
\eas
f_0: \C^4 &\longra \C^2 \\
(x,y,z,w) &\mapsto (x^2 + y^2 + z^2, w).
\eas
It is modelled on the quadratic cone $C_q = \{x^2 + y^2 + z^2 = 0, w = 0\}$, for which we know \cite[Ex. 2.26]{englebertConicallySingularCayley2023a} that $d(-1) = 2, d(0) = 8,  d(1)=22$ and there are no other critical rates in the range $[-1, 1]$. If this local model were part of a fibration of a compact $\Spin(7)$-manifold, then we would have for $\eps > 0$ small:
\eas
{\ind}_{-\eps} \CayD_{\CS} = {\ind}_{-\eps}\CayD_{\AC} = 2.
\eas
Thus in particular $\ind_{+\eps}\CayD_{\CS} = 2 + d(0) = 2+8=10$. This means that $\CayD w_t$ needs to lie in a codimension $10-\dim \ker \CayD_{N^{tA}} = 6$ subspace of $L^p_{k, \eps-1}(E)$ in order to perturb $w_t$ to a true solution with a small perturbation (i.e. using the bound \ref{3_2_inverse_bound}).
\end{ex}

We now go back to the deformation theory of $N$ as an unobstructed $\CS_\mu$ Cayley with moving points and cones. Using the notation of \ref{2_2_cs_deformation_op}, by solving the deformation problem we get a smooth submanifold $P \subset \cF_{\text{points}}$ of possible vertex locations and cone deformations of neighbouring $\CS_\mu$ Cayleys. We remark that the higher order deformation of the cone in a given $\CS$ Cayley is already determined by the translation applied to the point. Hence if $N$ has singular point $p$ and cone $C \subset \R^8$, then $T_{(p,C)}P$ can be identified with the possible translation directions of the conically singular point. This will be a subspace $\Th_N \subset T_pM$ of dimension $\dim \cM^\mu_{\CS}(N)$. Now we can decompose the kernel of $\CayD_{\AC}$ at rate $\eps > 0$ as follows: 
\e
\label{3_2_kernel_ac}
\ker \CayD_{\AC} = \cI^{\ze}_{\AC}(A) \op (\Th_N \op \Th_N^\perp).
\e 
Here $\cI^{\ze}_{\AC}(A)$ are the deformations of rate $\ze < 0$ that we discussed in the previous section, $\Th_N$ are the unobstructed directions of the conically singular problem, and $\Th_N^\perp$ is its orthogonal complement within the space of translations of the cone. We can choose compactly supported approximations $u_0$ of elements $u \in \Th_N^\perp$ such that $u \perp \Th_N$ pointwise. This then gives rise to a splitting of the pseudo kernel $\ka_t$ as follows: 
\eas
\ka_t = \ka_{\ze,t} \op \ka_{0,t} \op \ka_{\ob,t}.
\eas
Note that all sections in the family of pseudo-kernels are entirely supported on $N^{tA}_{\AC}$. We can now consider, as in the previous section, two infinitesimal deformation vector fields $w_{\AC} \in \cI^0_{\AC}(A), w_{\CS} \in \cI^0_{\CS}(N)$ with matching boundary conditions $ \si = \partial_0 w_{\AC} = \partial_0 w_{\CS}$. We then pre-glue them together as before to obtain $w_t \in C^\infty(\nu)$ with:
\eas
\nm{\CayD w_t}_{L^p_{k,\eps-1, t}} \lesssim t^\al\nm{w_t}_{L^p_{k+1,\eps, t}},
\eas
where $\eps > 0$ is a small constant. We now perturb slightly, so that $w_t \perp \ka_{\ob,t}$. For this we note that $\ka_{\ob,t}$ consists of compactly supported normal vector fields which are pointwise orthogonal to the unobstructed perturbation directions in the $\CS$ deformation problem (such as $w_{\CS}$). Hence $\nm{\pi_{L^2_{\pm\eps}}[w_t]}_{L^p_{k+1,\eps, t}} \ra 0$ as $t \ra 0$. Thus we end up with a perturbation $\tilde{w_t}$ such that we still have $\nm{\CayD \tilde{w}_t}_{L^p_{k,\eps-1, t}} \lesssim t^\al\nm{\tilde{w}_t}_{L^p_{k+1,\eps, t}}$ and additionally $\tilde{w}_t \perp \ka_{\ob,t}$. We are now in a position to run the argument from the previous section again using the bound \ref{3_2_inverse_bound} to obtain:

\begin{prop}
\label{3_2_parallel_def}
Let $N \subset (M, \Phi)$ be an unobstructed $\CS_\mu$ Cayley ($1 < \mu < 2$) with a unique singular point, and assume that its cone $C \subset \R^8$ is semi-stable. Let $A \subset \R^8$ be an unobstructed $\AC_\la$ Cayley ($\la <0$) with matching cone and sufficiently small scale. We then have that for any two matching infinitesimal deformation vector field $w_{\CS} \in \cI^{0}_{\CS} (N)$ and $w_{\AC} \in \cI^{0}_{\AC} (A)$ there are gluings $w_{\CS}\glue{t}w_{\AC} \in \cI(N\glu tA)$ such that (after identifying $\nu(N^{tA})\simeq \nu(N\glu tA)$) we have:
\ea
w_{\CS}\glue{t}w_{\AC} = w_t + \de w_t.
\ea
Here $\nm{w_t}_{L^p_{k,\ze,t}} = O(\md{\log t})$ and $\nm{\de w_t}_{L^p_{k,\eps,t}} \le t^\al \nm{w_t}_{L^p_{k,0,t}}$. 
\end{prop}

\subsubsection*{Stability of strong fibrations}

Now we have shown that if the nearly singular fibres of a Cayley fibration admit deformations of exactly two different rates, namely the normal deformations at rate $\ze < 0$ and the parallel deformations at rate $0$, then under the change of $\Spin(7)$ structure the infinitesimal deformation vector fields are perturbed by adding additional terms which are in $L^p_{k+1, \ze+\eps, t}(\nu)$ and $L^p_{k+1, \eps, t}(\nu)$ respectively, and are always bounded uniformly in $t$. We can now use this to show that the strong fibration property is stable under perturbation, given some additional assumptions. 

Let us assume that we have a strong Cayley fibration $f: (M,\Phi) \ra B$ as in definition \ref{2_3_caley_fib} with discriminant $\Delta \subset B$ of dimension $l = 1,2$. Already this means that all compact and conically singular Cayley fibres of $f$ are unobstructed in their respective moduli spaces. The singular cones of the conically singular Cayleys share their set of weights $\cD\subset \R$, and we let $\zeta = \max\{\cD \cap(-\infty, 0)\}$. We then require additionally that for each conically singular Cayley $N \subset M$ the Cayley operator $\CayD_N$ has index $4$ just below the critical weight $\ze$ and is unobstructed. We call such $N$ \textbf{simple}, and it can be seen that $N$ is simple exactly when $A$ is simple, so our discussion of the normal perturbations apply.

Consider now an atlas $\{(U_\al, \phi_\al)\}_{\al\in A}$ of of the base $B$, where $\phi_\al: U_\al \ra B_1(0) \subset \R^4$ is a diffeomorphism. If $U_\al \cap \Delta \ne 0$, then we further assume that this chart is compatible with the gluing map $\Ga$ in the sense that we identify $B_1(0) \simeq U_{\CS,\al}\times U_{\AC, \al}$, where $U_{\CS,\al} \subset \cM^\mu_{\CS}(N)$ and $U_{\AC,\al} \subset \brM^\la_{\AC}(A)$ with the condition that $f^{-1}(\phi^{-1}_\al(\tilde{N}, \tilde{A})) = \Ga(\tilde{N}, \tilde{A})$. On top of this we consider the framings $\{e_{i,\al}\}_{i=1,2,3,4}$ of $TB|_{U_{\al}}$ such $e_{i, \al}= \partial_i$. For each $b\in B \setminus \De$ we thus get four infinitesimal deformation vector fields $w_{1,\al}, \dots, w_{4,\al}$, which are just the lifts of $e_{i,\al}(b)$ via $f$. For $b\in \Delta$, we can again find a local frame, such that $T_b(U_{\CS,\al}) = \spn\{e_{1,\al}(b),\dots, e_{l,\al}(b)\}$ and $T_b(U_{\AC,\al})_{b} = \spn\{e_{l+1,\al}(b),\dots, e_{4,\al}(b)\}$. Note that in this case, $e_{1,\al},\dots, e_{l,\al}\in C^\infty_0(\nu)$ and $e_{l+1,\al},\dots, e_{4,\al}\in C^\infty_\ze(\nu)$. Note that at or near a singular point $w_{1,\al},\dots,w_{l,\al}$ are what we above called the parallel infinitesimal deformation vector fields and $w_{l+1,\al},\dots, w_{4,\al}$ are the orthogonal deformation vector fields. 

For a given $b \in U_\al$ we will now consider the following function in $C^\infty(f^{-1}(b))$: 
\e
\label{3_2_det}
{\det}_{\al,b} = \det(w_{1},\dots, w_{l},\rho^{-\ze} w_{l+1},\dots,\rho^{-\ze} w_{4}).
\e
We now make a final assumption on the initial fibration, namely that $\det_{\al,b}$ is bounded from below and above by constants $C$ and $C^{-1}$ respectively, for all $\al\in I$ and $b\in B$. This means that the infinitesimal vector fields never vanish for any Cayley in the fibration, and that the deformations of the cone of rate $\ze$ and $0$ have no zeros as well. We call such a fibration \textbf{non-degenerate}. Consider now a variation of the $\Spin(7)$-structure $\{\Phi_s\}_{s\in(-\eps, \eps)}$ with $\Phi_0 = \Phi$. The key insight here is the following:
\begin{lem}
If $\det_{\al,b}(s) > 0$ for all $\al, b$ and $0 \le s < s_{\max}$ then the universal families $ \Univ(\brM(N,\Phi_s))$ form strong fibrations of $M$.
\end{lem}
\begin{proof}
As we have $\dim \Delta \le 2$, we may apply the  weak fibration theorem \ref{1_1_fibration_weak} to conclude that small perturbations of $(M,\Phi)$ are still weakly fibreing. Thus if $N\subset M$ is any Cayley then
\eas
\ev_s: \Univ(\brM(N,\Phi_s)) \longra M
\eas
are homotopic maps for $s \in [0, s_{\max})$ of degree one in the sense of pseudo-cycles. Each of these maps is stratified smooth, thus in particular a local diffeomorphism on the open stratum as $\det_{\al, b}(s) > 0$ for every time $s \le s_{\max}$. Furthermore the maps $\ev_s$ remain orientation preserving for the same reason. But an orientation preserving local diffeomorphism of degree one must necessarily be a global diffeomorphism, as the algebraic count (which in this case is simply the naive count) of pre-images of a point equals the degree.
\end{proof}

 We then have that far away from the singular fibres, 
the $w_i$ perturb smoothly in $L^p_{k+1}(\nu)$. In particular, by compactness of $M$, we can insure that $\det_{\al,b} > \ha C > 0$ up until some $s_{\max, 0} > 0$, for all fibres (that are a given distance away from the singularities) simultaneously. Near the singular fibres, we see from the gluing results \ref{3_2_gluing_inf_def} and \ref{3_2_parallel_def} that $w_1,\dots, w_l$ perturb by continuously varying additional terms $\de w_{1,s},\dots, \de w_{l,s} \in L^p_{k, \eps, t}(\nu)$ and $w_{l+1},\dots,w_4$ additional terms $\de w_{l+1,s},\dots, \de w_{4,s} \in L^p_{k, \ze+\eps, t}(\nu)$. Now since $L^p_{k, \de+\eps,t} \hookra C^0_{\de,t}$ are continuous embeddings with bounded embedding constants, we see that $\det(w_{1},\dots, w_{l},\rho^{-\ze} w_{l+1},\dots,\rho^{-\ze} w_{4})$ varies continuously in $C^0$, uniformly in $t$. In other words, 
\eas
\nm{{\det}_{\al,b}(s)-{\det}_{\al,b}(0)}_{C^0} \le C_{\det}\md{s},
\eas 
where $C_{\det} > 0$ is independent of the neck size $t$. In particular, for a given chart we can find $s_{\max, \al} > 0$ such that $\det_{\al, b}(s) > \ha C > 0$ for any $s < s_{\max, \al}$. Hence, since we can cover $B$ with finitely many charts, this means that for small times $0 \le s < s_{\max}$ the we maintain the fibration property of the non-singular fibres. From this we can also deduce that the singular fibres do not intersect the non-singular fibres. Indeed, assume that for some $0 \le s < s_{\max}$ some singular fibre $\hat{N}$ intersected a non-singular fibre $N$. Then by what we just proved, the fibres near $N$ are locally still fibering, thus in particular for $t > 0$ sufficiently small, $\hat{N}\glu tA$ will intersect another non-singular fibre, which is of course impossible, as the non-singular fibres are still fibering for time $s$. Finally, as the singular fibres form a compact, smooth manifold, we can show that their infinitesimal deformation vector fields deform continuously in $C^0$ as well, and hence singular fibres will remain intersection-free as well. We thus proved:

\begin{thm}[Stability of strong fibrations]
\label{3_2_stability_strong}
Let $(M, \Phi)$ be an almost $\Spin(7)$-manifold which is strongly fibred by conically singular Cayleys which are simple, such that all the Cayleys in the fibration are unobstructed. Assume that the fibration is non-degenerate. Let $\Phi_t$ be a smooth deformation of the $\Spin(7)$-structure. Then there is an $\eps > 0$ such that for all $t \in [0, \eps)$ the manifold $(M, \Phi_t)$ can still be strongly fibred. 
\end{thm}

\section{Gluing construction of a Kovalev-Lefschetz fibration}

In the following we will give examples of strong fibrations by Cayleys of a family of non-torsion free $\Spin(7)$ manifolds given as products of coassociative fibrations as presented in Kovalev\cite{KovalevFibration} with $S^1$. These are fibrations $f: M^7 \times S^1 \ra S^3\times S^1$ obtained by gluing together two complex fibrations on a so called \textbf{twisted connected sum} $(M,\phi_T)$, where $T > T_{\min}$ and $\phi_T$ is a $G_2$-structure that admits a cylindrical neck of size $O(2T)$. The manifold $X^8 = M \times S^1$ does not admit a global holomorphic structure, but away from an interpolating region in the middle of the neck we can think of $f$ as being a complex fibration with the neighbourhoods of the singular points locally modelled on the fibration $f_0: (x,y,z,w)\mapsto (x^2+y^2+z^2, w)$ from the previous section, see also Example \ref{3_2_example_local}.

To apply our Theorem \ref{3_2_stability_strong} we need that all fibres in the fibration $f$ are unobstructed in their respective moduli spaces. Furthermore we need to check simpleness of the fibres and that the fibration is non-degenerate. We will turn our attention first to the unobstructedness. 

\subsection{The complex quadric}

For a moment let us focus on the local model $f_0$ near a singular point, given by the following holomorphic fibration:
\eas
f_0: \C^4 &\longra \C^2 \\
(x,y,z,w) &\longmapsto (x^2 + y^2 + z^2, w).
\eas
So we have $f^{-1}_0(0, \eta) \simeq C_q$ is a quadric cone and the nearby nonsingular fibres are the asymptotically conical Cayleys $A_\eps = f^{-1}(\eps, 0)$. We note that the holomorphic normal bundle $\nu^{1,0}(A_\eps)$ is trivial. To see this, consider the following two nowhere vanishing sections of $\nu^{1,0}(A_\eps)$:
\ea
\label{4_1_trivialisation}
s_{1,\AC}(x,y,z,w) &= \partial_w, \nonumber \\
s_{2,\AC}(x,y,z,w) &= \frac{\delbar_x + \delbar_y + \delbar_z}{\md{(x,y,z,0)}^2}.
\ea
Note that $s_{1,\AC}$ is an infinitesimal deformation corresponding to a translation and thus of rate $0$, whereas $s_{2,\AC}$ is of rate $-1$ and corresponds to a variation in the parameter $\eps$. In other words, the normal bundle of $A_\eps$ is trivial exactly because of the existence of the fibration $f_0$. Next, we want to prove a Liouville theorem for $A_\eps$. 
\begin{prop}[Liouville theorem]
\label{4_1_liouville_AC}
Any bounded holomorphic function on $A_\eps$ is constant.
\end{prop}
\begin{proof}
We can embed:
\eas
A_\eps \subset \bar{A}_\eps \subset \C P^4,
\eas
where $\bar{A}_\eps$ is the completion $\{x^2+y^2+z^2 = \eps u^2, w=0\} \subset \C P^4$ (here $\C P^4$ has homogeneous coordinates $[x:y:z:w:u]$). Since $\bar{A}_\eps$ is compact and nonsingular, any bounded holomorphic function on $\bar{A}_\eps$ is automatically constant. Now, using a removable singularities theorem in higher dimensions (such as \cite[Thm. 1.23]{voisinHodgeTheoryComplex2002}) and noting that $ \bar{A}_\eps \setminus A_\eps \subset \C P^4$ is a non-singular subvariety, we can extend any bounded holomorphic function on $A_\eps$ to a holomorphic function on $\bar{A}_\eps$, which concludes the proof.
\end{proof}
We can now use this to prove unobstructedness of $A_\eps$ as a Cayley.
\begin{prop}
\label{4_1_unobstructed_ac}
The $\AC_\la$ Cayley $A_\eps \subset (\R^8, \Phi_0)$ for $\eps \in \C\setminus 0$ is unobstructed and has no infinitesimal deformations at rate $-2 < \la < -1$. 
\end{prop}
\begin{proof}
Following \cite[Prop. 3.5]{mooreCayleyDeformationsCompact2019} we can write the Cayley operator on a complex $\AC_\la$ surface in a $(\R^8, \Phi_0)$ as:
\eas
\CayD_{\AC}=\delbar + \delbar^*: C^\infty_{\la}(\nu^{1,0}(A_\eps) \op \La^{0,2}A \ot \nu^{1,0}(A_\eps)) \longra C^\infty_{\la-1}(\La^{0,1}A_\eps \ot \nu^{1,0}(A_\eps)). 
\eas
Thus if $(u,v) \in C^\infty_{\la}(\nu^{1,0}(A_\eps) \op \La^{0,2}A_\eps \ot \nu^{1,0}(A_\eps))$ satisfies $\delbar u +\delbar^*v =0$ then the pair $(u,v)$ corresponds to an infinitesimal Cayley deformation vector field. If in addition we have $\delbar u = 0$ and $\delbar^* v = 0$, then $(u,v)$ is in fact an infinitesimal complex deformation \cite[Cor. 4.7]{mooreCayleyDeformationsCompact2019}. To start, we will prove that for $\la < -1$ any infinitesimal Cayley deformation is necessarily an infinitesimal complex deformation. For this note that if $\delbar u +\delbar^*v =0$ , then we automatically have $\delbar^* \delbar u = -\delbar^* \delbar^* v = 0$. Thus, since $\delbar u \in L^2(\La^{0,1}A_\eps \ot \nu^{1,0}(A_\eps))$ by our choice of rate $\la$, we find that:
\eas
0 = \int_{A_\eps} \langle \delbar^* \delbar u, u\rangle \dvol = \int_{A_\eps} \langle  \delbar u, \delbar u\rangle \dvol = \nm{\delbar u}_{L^2}.
\eas
In particular $\delbar u = 0$, which entails $\delbar^* v = 0$, and thus any infinitesimal Cayley deformation is in fact also infinitesimal complex. Now, since there are no bounded and non-constant holomorphic functions on $A_\eps$ by Proposition \ref{4_1_liouville_AC}, there are no other infinitesimal complex deformations of rate less than or equal to $-1$ besides constant multiples of $s_2$. Thus the kernel of the Cayley operator must be trivial at this rate:
\eas
{\ker}_{\la}\CayD_{\AC} = 0.
\eas
Finally, we prove the surjectivity of the Cayley operator at rate $\la$. This is equivalent to injectivity of its formal adjoint:
\eas
\CayD^*_{\AC}=(\delbar^*,\delbar): C^\infty_{-4-\la}(\La^{0,1}A_\eps \ot \nu^{1,0}(A_\eps)) \longra C^\infty_{-5-\la}(\nu^{1,0}(A_\eps) \op \La^{0,2}A_\eps \ot \nu^{1,0}(A_\eps)) . 
\eas
But we have $C^\infty_{-4-\la}(\La^{0,1}A_\eps \ot \nu^{1,0}(A_\eps)) = C^\infty_{-4-\la}(\La^{0,1}A_\eps \ot \C^2)$ since the normal bundle is trivial. So if $\delbar v = 0$ and $\delbar^* v = 0$, then in fact $v$ is a harmonic $1$-form with values in $\C^2$, as $A_\eps$ is K\"ahler. Now $v$ is square-integrable (by our assumption on the rate), and thus we can invoke \cite[Thm. 0.14]{lockhartFredholmHodgeLiouville1987}, which says that in this situation square-integrable harmonic one-forms are in one-to-one correspondence with elements of $\Ho^1(A_\eps) = 0$. Thus we get $v = 0$, and the Cayley operator is surjective.
\end{proof}

\subsection{Complex fibrations of Calabi--Yau fourfolds}

\begin{prop}
\label{4_2_unobstructed_cpt}
Let $f: M^4 \ra B^2$ be a complex fibration, where $M$ is a smooth Calabi--Yau fourfold and $B$ is a smooth, complex two-dimensional base. If a fibre $F$ is diffeomorphic to a non-singular K3 surface then it is unobstructed as a Cayley submanifold and has a four-dimensional moduli space.
\end{prop}
\begin{proof}
First, we have from Proposition \ref{2_1_cpt_cayley_linearisisation} that the index of a fibre $F$ as above is given by:
\eas
\ind \CayD_F = \ha(\si(F)+\chi(F)) - [F]\cdot[F] = \ha(-16 + 24) -0 = 4.
\eas
Here the self-intersection number $[F]\cdot[F]$ vanishes by the fibration property. The fibre $F$ admits at least $4$ Cayley deformations by perturbing to nearby fibres, which is equal to the index of the elliptic problem. Hence, showing unobstructedness is equivalent to showing that there are exactly $4$ infinitesimal Cayley deformations. Now by \cite[Lemma 4.7]{mooreDeformationTheoryCayley2017}, we have that infinitesimal Cayley deformations are necessarily infinitesimal complex deformations. However, because $F$ is part of a complex fibration locally, the holomorphic normal bundle $\nu(F) = \cO(F) \op \cO(F)$ is trivial and has $\Ho^0(\nu(F)) \simeq \C^2$ by compactness of $F$. This concludes the proof, as holomorphic normal vector fields are exactly the infinitesimal Cayley deformations.
\end{proof}
\begin{prop}
\label{4_2_unobstructed_cs}
Let $f: M^4 \ra B^2$ be a complex fibration, where $M$ is a smooth Calabi--Yau fourfold and $B$ is a smooth, complex two-dimensional base. Suppose that the fibration is modelled near a singular point on the complex quadric fibration 
\eas
f_0: \C^4 &\longra \C^2 \\
(x,y,z,w) &\longmapsto (x^2 + y^2 + z^2, w).
\eas
Assume furthermore that each singular fibre contains at most two singular points and that the non-singular fibres are diffeomorphic to non-singular K3-surfaces. Finally the singular locus $\Delta \subset B$ should take the form of a transverse intersection of smooth submanifolds. In that case each Cayley in the fibration is unobstructed in its moduli space. Here we allow for moving points and cones in the $\CS$ moduli space.
\end{prop}
\begin{proof}
We denote a non-singular fibre of the fibration by $F$, a singular fibre with a unique singular point by $F_{s}$ and a singular fibre with two singularities by $F_{ss}$. The expectation is that non-singular fibres are generic, fibres with one singularity appear in codimension $2$ and fibres with two singularities appear in codimension $4$. We will now show more precisely that the indices of the deformation problems are given by:
\eas
{\ind} \CayD_{F} = 4, \quad {\ind}_{1+\eps} \CayD_{\CS, F_s} = 2, \text{ and }\quad {\ind}_{1+\eps} \CayD_{\CS,F_{ss}} = 0.
\eas
Here $\eps > 0$ is small and the operators $ \CayD_{\CS, F_s} $ and $ \CayD_{\CS, F_s} $ are including the deformations of the points and cones. We first note that the equality ${\ind} \CayD_{F} = 4$ is the contents of Proposition \ref{4_2_unobstructed_cpt}. Next, the critical rates of the quadratic cone $C_q$ in the range $(-2,2)$ are known from Example \ref{2_2_ex_weights}, and have multiplicities:
\eas
d(-1) = 2,\quad d(0) = 8, \quad d(1) = 22,\quad d(-1+\sqrt{5}) = 6.
\eas
Thus using Theorem \ref{2_1_change_of_index} we see that the index of the problem with varying cones and points at rate $1+\eps$ is in fact equal to the index of the operator with fixed points and cones, but at rate $-1+\eps$. Hence, by gluing one or two matching $\AC$-manifold $A_\eps$ onto the conically singular point with we obtain non-singular $F \simeq F_s \glu A_\eps \simeq F_{ss} \glu A_\eps \glu A_\eps$. Thus we have (using ${\ind}_{-1+\eps} \CayD_{\AC} = 2$ from Proposition \ref{4_1_unobstructed_ac} ):
\eas
{\ind}_{1+\eps} \CayD_{\CS, F_s}  = {\ind}_{-1+\eps} \CayD^{\fix}_{\CS, F_s}  = \ind \CayD_{F}- {\ind}_{-1+\eps} \CayD_{\AC} = 4-2 = 2
\eas
and 
\eas
{\ind}_{1+\eps} \CayD_{\CS, F_{ss}}  &= {\ind}_{-1+\eps} \CayD^{\fix}_{\CS, F_{ss}}  \\
&= \ind \CayD_{F}- {\ind}_{-1+\eps} \CayD_{\AC} - {\ind}_{-1+\eps} \CayD_{\AC}\\
&= 4-2 -2= 0.
\eas
From this it is also clear we we should not expect unobstructed fibres with three or more singularities, as they would have strictly negative virtual dimension. We have now proven the index claims. In order to prove unobstructedness in the singular case (the compact case has been taken care of in Proposition \ref{4_2_unobstructed_cpt}) it is thus sufficient to prove that the space of infinitesimal Cayley deformations is exactly real two-dimensional or zero-dimensional respectively. First, consider the a fibre with a single conical singularity $F_s = f^{-1}(b) \setminus \{p\}$, with the conically singular point $p$ removed. Let $\partial_1, \partial_2 \in T_bB$ be two tangent vectors, where we assume $\partial_1 \in T_b \Delta$ and  $\partial_2 \not\in T_b \Delta$. As the differential $\D f$ only vanishes at the conically singular points, we see that the holomorphic sections of $\nu(F_s)$ given by:
\eas
s_1 = \D f^{*}[\partial_1], \\
s_2 = \D f^{*}[\partial_2], 
\eas
are nowhere vanishing, and thus span $\nu(F_s) = \cO(F_s) \op \cO(F_s)$. We note that $s_1$ has rate $O(1)$ when approaching the vertex $p$, as it comes from deforming the conically singular manifold to a nearby conically singular one (i.e. moving within $\Delta \subset B$). However as $\im \D f(p) = T_b \Delta$, we see that $\md{s_2}$ must diverge as we approach the cone. Indeed from the local model $f_0$ we see that $s_2$ must be asymptotic to $s_{2,\AC}$ from equation \eqref{4_1_trivialisation}, and thus of rate $O(r^{-1})$. Now we are in a position to repeat the proof of Proposition \ref{4_1_unobstructed_ac}. We first note that $F_s$ also has a Liouville theorem. Suppose that $h: F_s \ra \C$ is a bounded holomorphic function. Then blowing up $F_s$ at the conically singular point, we obtain a non-singular $\pi: \tilde{F}_s \ra N$ which is a biholomorphism away from a single exceptional and non-singular curve $Q = \pi^{-1}(p)$. We can then apply a removable singularities theorem in higher dimension \cite[Thm 1.23]{voisinHodgeTheoryComplex2002} to conclude that $h$ extends to a holomorphic function on $\tilde{F}_s$. Thus $h$ must be constant in the first place. Hence the only complex deformations of rate $0$ or above are the deformations coming from moving $F_s$ within the fibration. Now can use the same integration by parts argument that we used for the $\AC$ case to show that there are no further deformations which are Cayley but not complex. 

For the singular fibres with two singularities we again see that $\nu(F_{ss}) = \cO(F_{ss}) \op \cO(F_{ss})$. However now $F_{ss} = f^{-1}(p)$, where $p \in \Delta$ is a transverse intersection point. Thus deforming $p \in \Delta$ within $\Delta$ results in one singularity persisting, with the other one being resolved. Thus our discussion from above shows that all normal sections of $F_{ss}$ necessarily blow up with rate $O(r^{-1})$ near one of the singular points. In particular the conically singular fibres with two singularities are rigid and therefore unobstructed. 
\end{proof}

\subsection{Coassociative fibrations on twisted connected sums}
\label{4_3_twisted_connected}

In this section we briefly describe the twisted connected sum construction, first described by Kovalev \cite{kovalevTwistedConnectedSums2003} following an idea by Donaldson, and later extended by Corti-Haskins-Nordstr\"om-Pacini \cite{cortiManifoldsAssociativeSubmanifolds2015}. It gives rise to torsion-free $G_2$ manifolds via perturbation of an explicit small torsion glued $G_2$ manifold. From their construction, these pre-glued manifolds $M$ admit natural fibrations by coassociatives. Our stability theorem \ref{3_2_stability_strong} allows us to perturb the induced Cayley fibration on $M \times S^1$, which ultimately allows us to prove the existence of coassociative fibrations of $G_2$ manifolds as well.

\subsubsection*{Cylindrical Calabi--Yau 3-folds}

Let $(S,I,\om_\infty,g_\infty,\Om_\infty)$ be a K3 surface with a chosen hyperk\"ahler structure. Assume that $(X^6, J,\om,g, \Om)$ is a non-compact Calabi--Yau threefold. We say that $X$ is \textbf{asymptotically cylindrical} of rate $\la < 0$ or ($\ACyl_\la$), limiting to the hyperk\"ahler surface $S$ if there is a compact subset $K \subset X$ and a diffeomorphism $f: X \setminus K \ra \R_{>0} \times S^1 \times S$ with the following properties for all $k \ge 0$.
\eas
\md{\nabla^k(g - (g_\infty + \dt ^2 + \d s^2))}   &= O(e^{\la t}),\\
\nabla^k( \om - (\om_\infty + \dt\wedge \d s) )&= \d \si, \text{ where } \md{\nabla^k \si} =  O(e^{\la t}), \\
\nabla^k( \Om - (\dt+i\d s)\wedge \Om_\infty) &= \d \Si, \text{ where } \md{\nabla^k \Si} = O(e^{\la t}).
\eas
Here $\R_{>0}\times S^1$ has coordinates $(t,s)$ and $\md{\cdot}, \nabla$ are defined with respect to the product metric on $\R_{>0} \times S^1 \times S$. Asymptotically cylindrical Calabi--Yau threefolds can be constructed from compact Fano three folds using the following theorem.

\begin{thm}[Thm. 2.6 in\cite{cortiAsymptoticallyCylindricalCalabi2013}]
\label{acyl_from_fano}
Let $Z$ be a compact K\"ahler threefold with a morphism $f : Z \ra \C P^1$, with a smooth connected reduced fibre $S$ that is an anticanonical divisor, and let $V = Z \setminus S$. If $(S,J,\om_S,g_S,\Om_S)$ is a hyperk\"ahler structure on the complex surface $(S,J)$ such that $[\om_S] \in \Ho^{1,1}(S)$ is the restriction of the K\"ahler class on $Z$, then there is a CY3 structure $(V,J,\om,g,\Om)$ on $V$ which is asymptotically cylindrical to the CY3 cylinder $\R \times S^1 \times S$ induced by the hyperk\"ahler structure on $(S,J,\om_S,g_S,\Om_S)$.
\end{thm}

Now we will discuss briefly how to obtain such $f,Z$ and $S$ as in the theorem above. In fact Corti, Haskins, Nordstr\"om and Pacini impose extra conditions on maps $f: Z \ra \C P^1$ which make them more suitable for the twisted connected sum construction. 

\begin{dfn}[Building block]
A non-singular complex algebraic threefold $Z$ together with a projective morphism $f: Z \ra \C P^1$ is called a \textbf{building block} if the following conditions are satisfied:
\begin{enumerate}
\item The anti-canonical class $-K_Z \in \Ho^2(Z, \Z)$ is primitive, i.e. not an integer multiple of another class in $\Ho^2(Z, \Z)$.
\item The pre-image $S = f^{-1}(\infty)$ is a non-singular K3 surface and $S \sim -K_Z$ as divisors. 
\item If $k: \Ho^2(Z, \Z) \ra \Ho^2(S, \Z)$ is the map induced by the embedding $S \hookra Z$, then $\im k \hookra \Ho^2(S,\Z)$ is primitive, i.e. $\Ho^2(S,\Z) / \im k$  is torsion-free as an abelian group.
\item The groups $\Ho^3(Z, \Z)$ and $\Ho^4(Z, \Z)$ are torsion-free.
\end{enumerate}
\end{dfn}
There are multiple ways to construct building blocks. The first was introduced by Kovalev in \cite{kovalevTwistedConnectedSums2003} and starts with a Fano threefold as in Definition \ref{Fano_dfn}. This was later extended in \cite{cortiManifoldsAssociativeSubmanifolds2015} by Corti, Haskins, Pacini, Nordstr\"om to what they call semi-Fano threefolds, which can be thought of as desingularisations of certain mildly singular Fano varieties. They outnumber Fano threefolds by several orders of magnitude. Finally there is a different type of building block coming from K3 surfaces with non-symplectic involutions \cite{kovalevK3SurfacesNonsymplectic2011} which yields different examples still. 

In all these examples we obtain a building block $f: Z \ra \C P^1$ where the generic fibre of $f$ is a smooth K3 surface. Singular fibres appear in complex codimension $1$, but in general we cannot say much about the kinds of singularities that appear. Hence we will go through the first construction of building blocks (starting from (semi-)Fano threefolds) and give an example where we can determine the singularities explicitly.

Let now $X$ be a Fano threefold, such as for instance a smooth quartic in $\C P^4$. Then a generic anticanonical divisor in $X$ (which is effective by the Fano property) is a smooth K3 surface by a classical result of Šokurov \cite{sokurovSmoothnessGeneralAnticanonical1980}. We then make the assumption that the linear system $\md{-K_X}$ contains two non-singular members $S_0, S_\infty$ such that $C = S_0 \cap S_\infty$ is a transverse intersection, and thus a non-singular curve. In this case the pencil described by $S_0$ and $S_\infty$ exhausts $X$ and has base locus exactly $C$. If we now blow up $X$ at $C$ to obtain a new manifold $Z$, the pencil generated by the proper transforms $\tilde{S}_0$ and $\tilde{S}_\infty$ of $S_0$ and $S_\infty$ respectively will be base point free. Thus we obtain a holomorphic map $f: Z \ra \C P^1$ with generically smooth K3 fibres such that $f^{-1}(0) = \tilde{S}_0$ and $f^{-1}(\infty) = \tilde{S}_\infty$.   

\begin{prop}[Proposition 3.17 in \cite{cortiManifoldsAssociativeSubmanifolds2015}]
The map $f: Z \ra \C P^1$ determines a building block.
\end{prop}

\begin{ex}
\label{quartic_building_block}
Consider the following quartic polynomial on $\C P^4$ with homogeneous coordinates $[x_0:x_1:x_2:x_3:x_4]$:
\eas
P = x_0^4 + x_1^4 + x_2^4 + x_3^4 + x_4^4 + x_3^3(x_0 + 10 x_1 + 100 x_2).
\eas
Consider the smooth complex submanifold $X = \{P = 0\} \subset \C P^4$. Then, using the adjunction formula, we can see that the canonical bundle $\om_Q$ of $Q$ is given by:
\eas
\om_X = (\om_{\C P^4} \otimes \cO_{\C P^4}(Q))\lvert_X = (\cO_{\C P^4}(-5) \otimes \cO_{\C P^4}(4))\lvert_X = \cO_{\C P^4}(-1)\lvert_X.
\eas
In particular the anticanonical bundle $\om_X^\star = \cO(1)\lvert_X$ is ample, and the anticanonical divisors are exactly the hyperplane sections of $X$. So we can take for instance: 
\eas
S_0 = &\{x_3 = 0\}\cap X \simeq \{x_0^4 + x_1^4 + x_2^4 + x_4^4  = 0\} \subset \C P^3\\
S_\infty = &\{x_4 = 0\}\cap X \simeq \{x_0^4 + x_1^4 + x_2^4 + x_3^4 + x_3^3(x_0 + 10 x_1 + 100 x_2)\} \subset \C P^3.
\eas
Both are smooth K3 surfaces. They intersect transversely in a curve $C \simeq \{x_0^4 + x_1^4+x_2^4 = 0\} \subset \C P^2$. A general element of the pencil generated by $S_0$ and $S_\infty$ is the intersection of $X$ with the plane $\{ax_3 + bx_4\} = 0 \subset \C P^4$. The base point free pencil induced in $Z$ can be described outside the exceptional divisor as the map:
\eas
f:& \ Z\setminus E \longra \C P^1 \\
  & \ [x_0:x_1:x_2:x_3:x_4] \longmapsto [x_3:x_4].
\eas
There will be no singularities of members of this new pencil on the exceptional divisor as there are no singularities for the hyperplane sections down on $X$. Thus we can restrict our search for singularities to the complement of the exceptional divisor, i.e. we can work in the original quartic $X$. A point $x = [x_0:x_1:x_2:x_3:x_4]$ on $X \setminus C$ will be singular for a hyperplane section exactly when $\D P(x) \in \spn\{\d x_3, \d x_4\}$. Thus the singular points can be described as the subvariety $S\subset \C P^4$ defined by the set of equations:
\eas
\left\lbrace
\begin{array}{l}
P = 0,\\
\partial_0 P = 0, \\
\partial_1 P = 0, \\
\partial_2 P = 0 
\end{array}
\right.
\Leftrightarrow\left\lbrace
\begin{array}{l}
x_0^4 + \dots + x_4^4 + x_3^3(x_0 + 10 x_1 + 100 x_2) = 0,\\
4x_0^3 + x_3^3 = 0, \\
4x_1^3 + 10x_3^3 = 0, \\
4x_2^3 + 100x_3^3 = 0.
\end{array}
\right.
\eas
By Bézout's theorem, the algebraic count of solutions to this system of equations (meaning that we count points with their scheme theoretic multiplicity) is the product of the degrees of the equations, hence $3^3 \cdot 4$. In this specific case, the full number of solutions is attained, hence all of them have multiplicity one. To see this, note first that any non-zero solution must have $x_3 \ne 0$. So we are free to set $x_3 = 1$, and solve the final three equations $x_i^3 = c_ix_3^3$ ($c_i \in \C\setminus 0$ ) first. We are left with $3^3$ distinct possibilities for the tuple $(x_0, x_1, x_2)$. Now for each such choice we can solve the first equation for $x_4$ in exactly four different ways, as it reduces to an equation of the form $x_4^4 = c(x_0,x_1,x_2)$ where $c \ne 0$.

Notice also that no two solutions lie in the same hyperplane section, as they all have different values of $[x_3:x_4]$. Indeed once $x_3$ is chosen, this determines $x_0, x_1, x_2$ up to a choice of a third root of unity. Now $x_0 + 10x_1 + 100x_2$ can never take identical values for $x_0, x_1, x_2$ differing only by a multiple of a root of unity. This explains the slightly odd choice of $x_0 + 10x_1 + 100x_2$ instead of something more symmetric like $x_0 +x_1 + x_2$ for instance. In the latter case permuting $x_0, x_1, x_2$ while keeping $x_3, x_4$ the same maps the singular set onto itself, and thus multiple singularities appear on one fibre.

Now as mentioned above, all points of $S$ have multiplicity $1$. That means that if $[x_0:x_1:x_2:x_3:x_4]\in S$ we have that $\dim \cO_{S,p} = 1$, where $\cO_{S,p}$ denotes the local ring of $S$ at $p$. Now fix a singular point $p \in S$, which is obtained by intersecting $X$ with a hyperplane $\Pi$ in $\C P^4$. By choosing a coordinate chart for an affine plane in $\Pi$ and after a further translation we can assume that our singular fibre looks like $f^{-1}(0)$ for a polynomial map $f: \C^3 \longra \C$ which additionally satisfies $f(0) = 0$ ($0_{\C^3}$ corresponds to $p\in S$) and $\D f(0) = 0$ (thus every term in $f$ is at least of second order). In this picture we see that:
\eas
\cO_{S,p} = \frac{\C[x_0,x_1,x_2]}{(\partial_0 f, \partial_1 f, \partial_2 f)}. 
\eas
We now claim that if the dimension of this local ring is $1$, then we can choose coordinates such that $f = x_0^2 + x_1^2 + x_2^2 + O(x^3)$. In particular it suffices to show that if the quadratic terms of $f$ do not form a non-degenerate quadratic form, then $\dim  \cO_{S,p} > 1$. Suppose that this is the case, so that after a linear change of coordinates we can assume that $x_0^2$ does not appear as a term in $f$. Then we clearly have $\partial_0 f = c x_0^2 + O(x_0^3, x_1, x_2)$, and similarly $\partial_1 f$ and $\partial_2 f$ do no contain a linear term proportional to $x_0$. Thus $1$ and $x_0$ are non-zero and linearly independent elements of $\cO_{S,p}$, and thus $\dim  \cO_{S,p} > 1$.

So in particular we have proven that all the singularities that appear in this example of a building block are isomorphic to the quadratic cone singularity $x^2 + y^2 + z^2 = 0$ in $\C^3$, as all the singularities have multiplicity one. This can alternatively also be checked explicitly by looking at the defining equations of $S$ in more detail. 
 
The property of only having $A_2$-type singularities, all in separate fibres, is Zariski open, i.e. it is true for $X$ in an open subset of its deformation type and for $S_0, S_\infty$ in an open subset of the corresponding linear anticanonical system. Since this moduli space is irreducible as a complex variety, it is generically true for Fano threefolds arising from quartics in $\C P^4$. Thus we showed the same is true for a dense Zariski open subset of $\cF^{N,A}$, where $A = -K_X \in N = \Ho^2(S_0)$.
\end{ex}

\subsubsection*{Twisted connected sum construction of $G_2$-manifolds}
Now we have established the basic properties of building blocks, which by Theorem \ref{acyl_from_fano} can be used to construct $\ACyl$ Calabi--Yau threefolds. Starting from a building block $f_Z: Z \ra \C P^1$ with chosen K3 fibre $K = f_Z^{-1}(\infty)$ (seen as a complex manifold), we can choose a hyperk\"ahler structure $(\om_\infty, g_\infty, \Om_\infty)$ compatible with $(K,I)$, under the condition that $[\om_\infty] \in \Ho^{1,1}(K)$ is the restriction of a K\"ahler class on the ambient $Z$. This is an open condition, but may be non-trivial.

By taking the product with $S^1$ we get an asymptotically cylindrical $G_2$ manifold $M = X \times S^1$ with associative form $\phi$, defined by $\phi = \dt \wedge \om + \Re \Om$ as in Example \ref{example_g2_mflds}. As $Z\setminus S$ and $X$ are biholomorphic, we see that $X$ is also fibred by generically smooth K3 surfaces via the same map $f_X = f_Z: X \ra \C P^1 \setminus \{\infty\}$. This induces a corresponding fibration by coassociative submanifolds on $f: M \ra S^1 \times \C P^1 \setminus \infty$.

The fibration on the cylindrical end of $M$ is diffeomorphic to the projection map $\pi: \R_{>0}\times S^1 \times S^1 \times K \mapsto \R_{>0}\times S^1 \times S^1$. By the $\ACyl_{\la}$-condition (with $\la < 0$) the metric on the link converges exponentially to $g_{S^1} \times g_{S^1} \times g_\infty$. The key idea of the twisted connected sum construction is to take two cylindrical $G_2$ manifolds $M_+, M_-$ with isometric asymptotic hyperk\"ahler K3's $K_+, K_-$ and glue them together by a diffeomorphism for $T > 0$:
\ea
\label{g_2_gluing_map}
G: (T,T+1)\times S^1 \times S^1 \times K_+ &\longra (T, T+1)\times S^1 \times S^1 \times K_- \\
 (t, \th_a, \th_b, p) &\longmapsto (2T+1-t, \th_b, \th_a, \mathfrak{r}(p)),
\ea
where $\mathfrak{r}: K_+ \ra K_-$ is a suitably chosen isometry. We exchange the two circles with the gluing diffeomorphism so that the fundamental group of the glued manifold becomes finite, and thus the holonomy will be exactly $G_2$ by a result of Joyce \cite[Prop. 10.2.2]{joyceCompactManifoldsSpecial2000}. In terms of the hyperk\"ahler structure on $(K_{\pm},\om_{\pm}^1, \om_{\pm}^2, \om_{\pm}^3, I_\pm, g_\pm)$ the asymptotic associative form can be written as:
\eas
\phi_{\infty,\pm} &= \d\th_a \wedge \om_{\C \times K_{\pm}} + \Re \Om_{\C \times K_{\pm}} \\
&= \d \th_a \wedge ( \dt \wedge \d \th_b+\om_\pm^1) + \Re (\d \th_b - i\dt)\wedge(\om_\pm^2 +i\om_\pm^3) \\
&= \d \th_a \wedge \dt \wedge \d \th_b + \d\th_a \wedge \om_\pm^1  + \d \th_b \wedge \om_\pm^2 + \dt \wedge \om_\pm^3.
\eas
In particular to ensure that $\phi_{\infty,\pm}$ match up on the overlap, we need:
\eas 
\mathfrak{r}^{*} \om_-^1 = \om_+^2, \quad \mathfrak{r}^{*} \om_-^2 = \om_+^2, \quad \mathfrak{r}^{*} \om_-^3 =- \om_+^3,
\eas
which is equivalent to asking that $\mathfrak{r}$ is a hyperk\"ahler rotation between $K_\pm$ as in equation \eqref{hyperkaehler_twist}. Let the parametrisations of the ends of $M_{\pm}$ as cylinders be denoted by $\phi_\pm: \R_{>0} \times S^1 \times S^1 \times K_{\pm} \hookra M$. For $T > 1$ we consider the following truncated manifolds:
\eas
M_{T,\pm} = M_\pm \setminus \phi_\pm([T+1, \infty)\times S^1 \times S^1 \times K_{\pm}).
\eas
Over the cylindrical end  $\R_{>0}\times L\simeq  \R_{>0} \times S^1 \times S^1 \times K_{\pm}$ we have both the $G_2$ structure $\phi_\pm$ induced from $M_\pm$ as well as the product $G_2$-structure $\phi_{\infty, \pm}$. Define one an for all a smooth cut-off function $f_{cut}: \R \times [0,1]$ such that $f\lvert_{(-\infty,0]} = 0$ and  $f\lvert_{[1,\infty)} = 1$. We can now define a (non torsion free) $G_2$-structure interpolating between the two on $M_{T,\pm}$ by declaring it equal to $\phi_\pm$ away from the cylindrical end, and on the end by the formula:
\eas
\phi_{T,\pm}(t,p) = f_{cut}(t-T) \phi_\pm + (1-f_{cut}(t-T))\phi_{\infty,\pm}. 
\eas
For $T \gg 1$ we will have $\md{\phi_\pm-\phi_{\infty,\pm}}_{(T,T+1)\times L} = O(e^{\la T})$ small. Hence $\phi_{T,\pm}(t,p)$ will be a small perturbation of $\phi_{\infty,\pm}$, and thus again a $G_2$-form. Notice that $\phi_{T,\pm}$ is exactly equal to $\phi_{\infty,\pm}$ on $(T+1,T+2)\times L$ and torsion-free everywhere except over the interpolation region $(T,T+1)\times L$, where $\md{\partial \phi} \in O(e^{\la T})$. Thus, after choosing a hyperk\"ahler rotation matching up $M_\pm$ we can glue $M_{T,\pm}$ over the regions $(T+1,T+2)\times L \subset M_{T,\pm}$ using the gluing map $G$ from \eqref{g_2_gluing_map} to obtain a $G_2$-manifold $(M_{T}, \phi_{T,\mathfrak{r}})$. This can be perturbed to a torsion-free $G_2$ manifold.
\begin{thm}[3.12 in \cite{cortiManifoldsAssociativeSubmanifolds2015}]
\label{G2_deformation}
Let $(X_\pm ,J_\pm,\om_\pm,g_\pm,\Om_\pm)$ be two asymptotically cylindrical Calabi--Yau 3-folds whose asymptotic ends are of the form $\R_{>0}\times S^1\times K_\pm$ for a pair of hyperk\"ahler K3 surfaces $K_\pm$, and suppose there exists a hyperk\"ahler rotation $\mathfrak{r}: K_+ \ra K_-$. Define closed $G_2$–structures
$\phi_{T,\mathfrak{r}}$ on the twisted connected sum $M_{\mathfrak{r}}$ as above. 
For sufficiently large $T$ there is a torsion-free perturbation of $\phi_{T,\mathfrak{r}}$ within its cohomology class.
\end{thm}

It can be shown that this perturbation will become arbitrarily small as $T$ increases. The most difficult aspect of the gluing construction is certainly finding pieces with compatible Calabi--Yau cylindrical ends. This we call the \textbf{matching problem}. The asymptotically cylindrical Calabi--Yau threefolds we consider come from building blocks, which in turn come from (semi-)Fano threefolds with a choice of anticanonical K3 divisors. We now give an outline of the matching procedure from \cite{cortiManifoldsAssociativeSubmanifolds2015}. Consider the deformation types of two Fano manifolds $Y_\pm$ which are polarised by the lattices $N_+ \subset \Lambda$ and $N_- \subset \La$ respectively. Assume that $N_\pm$ has signature $(1, r_\pm)$. Recall the forgetful morphisms $s^{N_\pm}: \cF^{N_\pm} \ra \cK^{N_\pm}$ which takes pairs $(\tilde{Y}_\pm,S)$ of Fano threefolds in the deformation type of $Y_\pm$ and anticanonical K3 divisors $S \subset \tilde{Y}_\pm$ to the polarised K3 moduli space $S \in \cK^{N_\pm}$. We know from Proposition \ref{Fano_K3_surjective} that this morphism is dominant on each irreducible component of $\cF^{N_\pm}$. This gives us our first restriction on K3 surfaces which can be used in the matching, as they must lie in open dense subsets $U_\pm \subset D_{N_\pm}$ which are determined by $s^{N_\pm}$ and a reference marking. The next step is to consider the hyperk\"ahler structure. To make the discussion simpler, we assume that the lattices $N_\pm$ have trivial intersection and are orthogonal to one another. This way we can avoid introducing the construction of an \textbf{orthogonal pushout} of two lattices and also have more concise notation. Define $T = (N_{+} \op N_{-})^{\perp}$. Consider the following subset of the hyperk\"ahler K3 domain:
\eas
D = D^{\text{hk}}_{K3} \cap ((N_+ \ot \R)^+ \times (N_- \ot \R)^+ \times (T \ot \R)).
\eas
Here $L^+$ denotes the positive cone of a lattice $L$. The submanifold $D$ is real $20$-dimensional. Now as the period domain of $N_{\pm}$-polarised K3 surfaces can be identified with positive two-planes in $N_\pm^\perp \ot \R$ as in Equation \eqref{K3_polarised_period_domain}, there are two natural projection maps $\pi_\pm: D \ra D_{N_\pm}$ given by:
\eas
\pi_\pm(\om_+, \om_-, \om_0) = \spn \langle \om_{\mp}, \pm \om_0\rangle \in D_{N_\pm}.
\eas 

Recall from Equation \eqref{hyperkaehler_twist} that the hyperk\"ahler rotation of $(\om_+, \om_-, \om_0)\in D^{\text{hk}}_{K3}$ is exactly $(\om_-, \om_+, -\om_0)$. Hence $\pi_-$ is just the mapping to the complex structure of the hyperk\"ahler rotated K3. So in particular candidates for asymptotic hyperk\"ahler K3 surfaces must be contained in the subset $\pi_+^{-1}(U_\pm) \cap \pi_-^{-1}(U_-)$. It turns out that the image of $\pi_\pm$ in $D_{N_\pm}$ is a real $(20-r_\pm)$-dimensional submanifold in a real $2(20-r_\pm)$-dimensional space. So it is not clear a priori that $\im \pi_\pm$ and $U_\pm$ even intersect. However one can show that $\im \pi_\pm$ is an embedded totally real submanifold. Since it is also of maximal dimension it must intersect any open Zariski dense subset, such as $U_\pm$. In face, one can show that the complement of $\pi_\pm^{-1}(U_\pm)$ is a finite collection of real analytic subsets of positive codimension.
Finally, we also need to take into account not only the complex geometry of the two K3 surfaces, but their K\"ahler geometry as well. Indeed if the hyperk\"ahler structure is given by the triple $(\om_+, \om_-, \om_0)$, then the complex geometry is determined by $(\om_-, \om_0)$ (via $\pi_+$ as above), while the K\"ahler class will be $\om_+$, and similarly for the hyperk\"ahler rotation. So if $\Amp_\pm$ are the ample cones of the polarised K3 surfaces we need that the set:
\eas
\cA = \{ (\om_+, \om_-, \om_0): \om_\pm \in \Amp_\mp\}
\eas
is non-empty. In good cases this can be shown to be a (Euclidean) open subset of $D$. Thus, since $\pi_+^{-1}(U_\pm) \cap \pi_-^{-1}(U_-)$ is open dense, there must be a hyperk\"ahler structure satisfying all the conditions, and thus the matching is possible. We note at this point that imposing a finite number of open dense conditions on either complex K3 surface does not impact the matching procedure. 

\begin{ex}
\label{example_matching_building_block}
Consider the twisted connected sums of two building blocks in the deformation type of example \ref{quartic_building_block}. In this case the matching is possible by the work of Corti, Haskins, Nordstr\"om and Pacini \cite[Prop. 6.18]{cortiManifoldsAssociativeSubmanifolds2015} (the matching is what they call \textbf{perpendicular}, and thus the condition on the K\"ahler classes are automatically satisfied) and we can see from Table 5 in \cite{cortiManifoldsAssociativeSubmanifolds2015} that the resulting $G_2$-manifold will have $b^2 = 0$ and $b^3 = 155$. The example is investigated in more detail as Example Number one in section 7 of \cite{cortiManifoldsAssociativeSubmanifolds2015}.
\end{ex}

\subsubsection*{Coassociative fibrations on connected sums}
We now turn our attention to the fibrations that come automatically with the twisted connected sum construction of $G_2$ manifolds.

\begin{prop}[Prop. 2.18 in \cite{KovalevFibration}]
The fibrations $f_\pm: M_\pm \ra S^1\times \C P^1 \setminus \infty$ join together to form a fibration $f_T: M_{\mathfrak{r}} \ra S^3$.
\end{prop} 
\begin{proof}
Gluing $M_{T,\pm}$ identifies the K3 fibres of the two fibrations by construction. On the level of the base space, this reduces to a gluing of two solid tori $S^1 \times D^2$ (where $D^2$ is the two dimensional disk with boundary) along their boundaries via the map
\eas
S^1 \times S^1 \longra S^1\times S^1, \quad  (a,b) \longmapsto (b,a).
\eas 
This gluing is diffeomorphic to $S^3$, and the decomposition into Tori is in fact a Heegard splitting of $S^3$. Consider now a fixed K3 surface on the overlap $f_\pm^{-1}(t,\th_a, \th_b)$. It is coassociative with respect to $\phi_{\pm}$ by construction and we easily see that it also is coassociative with respect to $\phi_{\infty,\pm}$. Thus is remains coassociative for any linear combination $c \phi_{\pm} +(1-c)\phi_{\infty,\pm}$ with $0 \le c \le 1$.
\end{proof}

Let us now apply the twisted connected sum construction to building blocks with the additional property that the fibres of the map $f: Z \ra \C P^1$ are either non-singular K3 surfaces or have conical singularities modelled on the cone $C_q  = \{x^2 +y^2 +z^2 = 0\}$. A possible building block arises from a quartic in $\C P^4$ as explained in Example \ref{quartic_building_block}. Ensuring the matching up of two building blocks $(Z_\pm, f_\pm)$ is an involved procedure, as we already mentioned above. But since the additional condition we impose is Zariski open dense on the moduli space $\cK^{N,A}$, the matching goes through as without change. Thus we can find two matching building blocks so that the glued fibrations $f_T$ also have the same kind of complex conical singularities. On the tubular intermediate region all the fibres will be smooth K3s.

\begin{figure}
\label{figure_pieces}
\begin{center}
\includegraphics[scale=0.2]{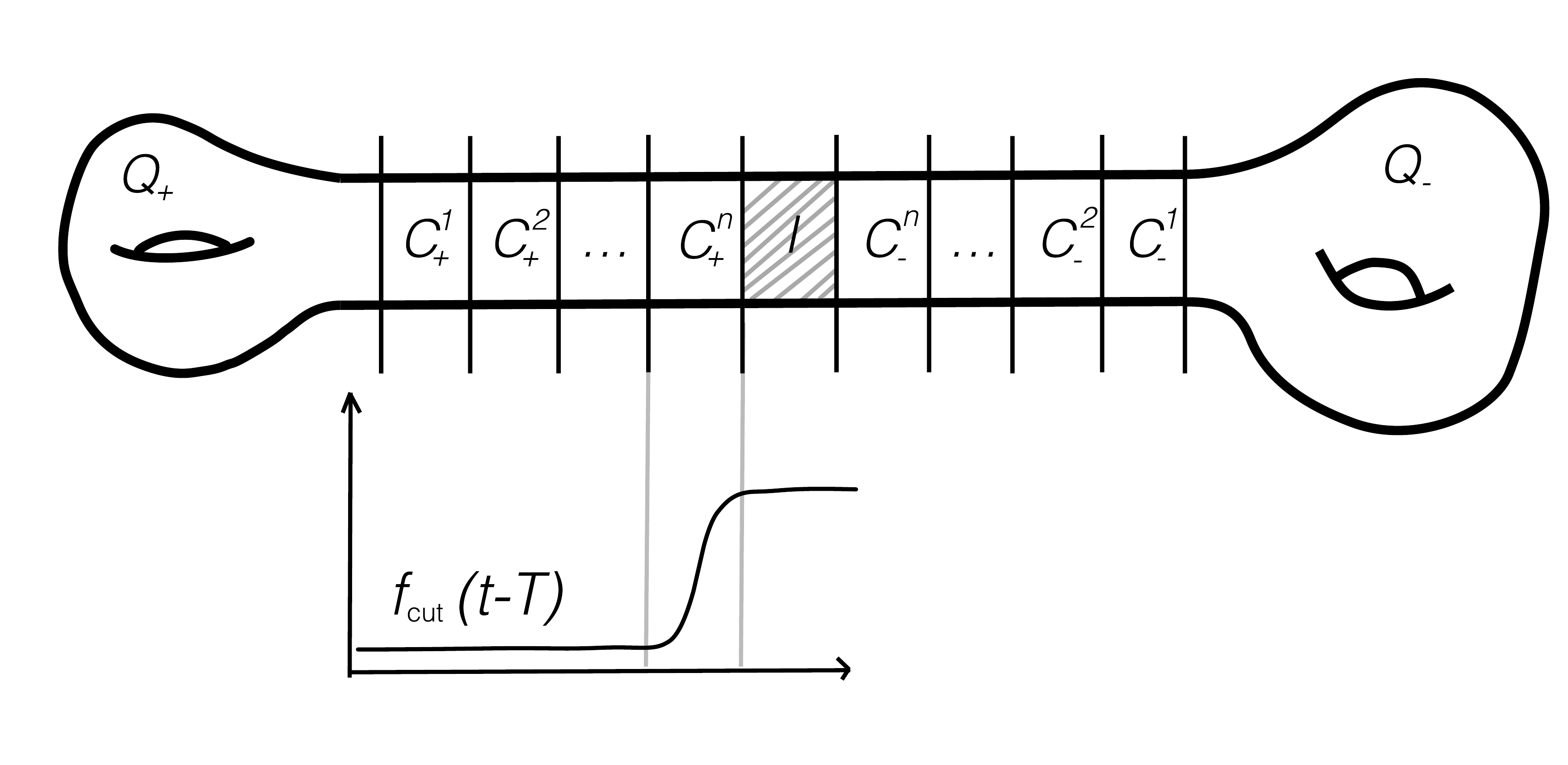}
\end{center}
\caption{Decomposition of the twisted connected sum into $2n+3$ pieces.}
\end{figure}
The upshot is that we are given a smooth twisted connected sum manifold $M^7$ and for any $T \gg 1$ sufficiently large a $G_2$-structure $\phi_T$ such that $\nm{\nabla \phi_T}_{L^p_{k}} \le e^{\la T}$. These come with coassociative fibrations by (possibly singular) K3 surfaces, which over either end are products of complex fibrations of Calabi--Yau threefolds with a circle $S^1$. In the gluing region there is no complex structure for which this is true, however this region of finite volume, and converges to a cylindrical $G_2$ manifold. If we take the product with $S^1$ once more we obtain a $\Spin(7)$-manifold $(X = M \times S^1, \Phi_T = \dt \wedge \phi+\star \phi)$, which admits a torsion-free deformation $\tilde{\Phi}_T$, also of product type, from Theorem \ref{G2_deformation}. The deformation required to achieve torsion freeness becomes smaller as $T$ increases, in the sense that $\nm{\Phi_T - \tilde{\Phi}_T}_{L^p_k} = O(e^{\la T})$. Now geometrically either end of $X$ can by construction be considered a Calabi--Yau fourfold with a fibration by complex surfaces, and the gluing region is approaching a cylindrical $\Spin(7)$-manifold of finite volume. In particular we can use Propositions \ref{4_2_unobstructed_cpt} and \ref{4_2_unobstructed_cs} to see that away from the gluing region the fibres, both compact and singular, are unobstructed as Cayleys. Regarding the gluing region we note that unobstructedness is an open condition in the choice of $\Spin(7)$-structure. Thus while we cannot apply Proposition \ref{4_2_unobstructed_cpt} directly, we see that it applies to the limiting cylindrical $\Spin(7)$ structure. Hence the fibres in the gluing region will be unobstructed for all $T \gg 1$. Now we are in a position to prove the stability of the fibration as we pass from the $\Spin(7)$ structure $\Phi_T$ to the torsion-free $\Spin(7)$-structure $\tilde{\Phi}_T$. For this, we imagine cutting up $(M_n,\Phi_n)$ (for $n \in \N$), a manifold of diameter approximately $2n$ into $2n+3$ pieces. These pieces are first of all the two compact pieces $Q_\pm \subset X_\pm$. Next we have for either side the $n$ pieces $C_{k,\pm}= \phi_\pm^{-1}((k-\ha,k+1\ha)\times L)\times S^1$ for $0 \le k \le n$. Finally we have the glued piece $I = \phi_\pm^{-1}((n+\ha,n+2\ha)\times L)\times S^1$. Notice that both $Q_\pm$ as well as the $C_{k,\pm}$ for $0\le k\le n-1$ when seen as $G_2$-manifolds remain constant as $n$ increases. The two pieces $C_{n,\pm}$ are where the interpolation between the $\ACyl_\la$ structure $\phi_\pm$ and the exactly cylindrical $G_2$ structure $\phi_{\infty,\pm}$ happens. Finally $I$ is exactly cylindrical, independent of $n$. Now, since we checked unobstructedness of all the fibres in the fibration we may apply Theorem \ref{1_1_fibration_strong} to each piece separately, as long as we can ensure non-degeneracy of the fibration. This is clearly satisfied for any piece without a singular fibre. The finitely many pieces with singular fibres can be considered as Calabi--Yau fourfolds with fibrations by complex surfaces and Morse-type singularities. From this we can conclude non-degeneracy, since we know the local model near the singular point.

Thus for each piece there is a maximal $s_{\max}\in (0,1]$ such that for each $0 \le s< s_{\max}$ the fibration property of the fibres in just that piece is preserved for $\Phi_{n,s} = \Phi_n + s(\tilde{\Phi}_n - \Phi_n)$. Now, looking at one of the pieces after they reached their steady state, we see that as $n$ increases, $s_{\max}$ for that piece increases and reaches $1$ eventually. This is because  $\nm{\tilde{\Phi}_n-\Phi_n} \le e^{\la n}$ with $\la < 0$. Eventually $\tilde{\Phi}_n$ will line in the open ball about $\Phi_n$ where $s_{\max} = 1$. Like this we see that for any choice of finitely many pieces, we can ensure the stability of the fibration on the union of these pieces for any sufficiently large $T$.

 On the other hand we have $(C_{k,\pm}, \phi\lvert_{ C_{k,\pm}}) \ra (I,\phi\lvert_I)$ in $C^\infty$. In fact if we consider the path of $\Spin(7)$ structures $\ga(T) =(-T+\ha)^*\Phi_\pm\lvert_{[T-\ha,T+\frac{3}{2}]\times L}$ (where $(-T+\ha)^*$ is the pullback by translation) on $[0,2]\times L$, then the $\ACyl$ condition on $M_\pm$ gives us that $\nm{\ga(T)-\Phi_{\infty}}_{L^p_k} \le e^{\la T}$, for $\la < 0$. Stability of the fibration is true in a quantitative sense, meaning that there is a ball $B(\Phi_{\infty},\eps)$ around $\Phi_{\infty}$ where the moduli space of Cayleys for the given $\Spin(7)$-structure is still fibreing. Thus there will be a ball of radius $\eps - e^{\ha T}$ around $\ga(T)$ so that the same is true. But now, since the distance $\nm{\tilde{\Phi}_T-\Phi_T}_{L^p_k} \le e^{\la T} < 1-e^{\la T} $ for $T$ sufficiently large, the torsion-free $\Spin(7)$-structure will stay within ball of fibering $\Spin(7)$ structures around $(C_{k,\pm}, \phi_{k,\pm})$. In this way we can thus prove stability of the fibration for all pieces with index above a minimal $n_{\min}$. The previous argument then takes care of the finitely many remaining pieces.
Thus we have proven the following:
\begin{thm}[Existence of strong Kovalev-Lefschetz fibrations on $\Spin(7)$ manifolds]
There are compact, torsion-free $\Spin(7)$-manifolds of holonomy $G_2$ which admit fibrations by Cayley manifolds.
\end{thm} 

From this we can deduce that the stability result also holds for the initial $G_2$-manifold, using the following auxiliary result.

\begin{lem}
Let $(M,g,\tau)$ be a manifold together with a calibration $\tau$. Assume that $\tau = \psi + \rho$, where $\psi$ is another calibration and $\rho$ is a closed form. Let $N \subset M$ be $\tau$-calibrated submanifold such that $\int_N \rho = 0$. Then $N$ is $\psi$-calibrated.

\end{lem}
\begin{proof}
We have that by assumption $\dvol_N = \tau|_N$. Now note that:
\eas
\int_N \dvol_N = \int_N \tau|_N = \int_N \psi|_N + \rho|_N=\int_N \psi|_N\le  \int_N \dvol_N.
\eas
If there is a point with $\psi|_N(p) <  \dvol_N(p)$, then we must have $\int_N \psi|_N<  \int_N \dvol_N$, a contradiction. Thus $\psi|_N= \dvol_N$ and $N$ is $\psi$-calibrated.
\end{proof}
Now we set $\tau = \Phi$, $\psi = \star \phi$ and $\rho = \d s \wedge \phi$ in the previous proposition, where $s \in S^1$ is a coordinate on the circle in $X = M \times S^1$. As a K3 fibre $N$ of the initial fibration are contained in $M \times \{s\}$ for a single point $s \in S^1$, we clearly have $\int_{N} \d s \wedge \phi = 0$. As the perturbed Cayleys are continuous deformations of the initial (possibly conically singular) Cayleys and the new $G_2$-structure $\tilde{\phi}$ is cohomologous to $\phi$, we still have $\int_{\tilde{N}} \d s \wedge \tilde{\phi} = 0$ by Stokes' Theorem. Hence the Cayleys for $\tilde{\Phi}$ are also calibrated by $\star\phi$, meaning that their tangent planes are contained in $M \times \{s\}$ and they are in fact coassociative. Thus we have shown:

\begin{cor}[Existence of coassociative fibrations]
There are compact, torsion-free $G_2$-manifolds of full holonomy which admit fibrations by coassociative submanifolds.
\end{cor}

\begin{ex}
Consider the $G_2$-manifold obtained by gluing two copies of the quartic building block from Example \ref{quartic_building_block}, as in Example \ref{example_matching_building_block}. This $G_2$-manifold has Betti numbers $b_2 = 0$ and $b^3 = 155$. Furthermore, as the the conical singularities are stable and no fibre has more than one singular point, the resulting coassociative fibration will have $2\cdot 3^3 \cdot 4 = 216$ connected components of singular coassociatives. 
\end{ex}

\begin{rem}
We expect that similarly well-behaved fibrations can be constructed on many more semi-Fano threefolds, thanks to a discussion with Mark Haskins. This is because the anticanonical system is usually quite large (e.g. it has complex dimension 4 in the case of the quartic in $\C P^4$), and thus it should be possible to avoid bad singularities by choosing a suitably generic pencil and invoking a Bertini-type theorem. Morally speaking this should be true, as having Morse type singularities in separate fibres is a codimension $2$ condition. 
\end{rem}

\addcontentsline{toc}{section}{References}
\bibliographystyle{alpha}
\bibliography{library}

\begin{thebibliography}{CHNP15}

\bibitem[Bea02]{beauvilleFanoThreefoldsK32002}
Arnaud Beauville.
\newblock Fano threefolds and {{K3}} surfaces.
\newblock \url{https://doi.org/10.48550/arXiv.math/0211313}, November 2002.

\bibitem[CHNP13]{cortiAsymptoticallyCylindricalCalabi2013}
Alessio Corti, Mark Haskins, Johannes Nordstr{\"o}m, and Tommaso Pacini.
\newblock Asymptotically cylindrical {{Calabi}}--{{Yau}} 3--folds from weak
  {{Fano}} 3--folds.
\newblock {\em Geometry \& Topology}, 17(4):1955--2059, July 2013.

\bibitem[CHNP15]{cortiManifoldsAssociativeSubmanifolds2015}
Alessio Corti, Mark Haskins, Johannes Nordstr{\"o}m, and Tommaso Pacini.
\newblock G 2 -manifolds and associative submanifolds via semi-{{Fano}} 3
  -folds.
\newblock {\em Duke Mathematical Journal}, 164(10):1971--2092, July 2015.

\bibitem[Eng23a]{englebertConicallySingularCayley2023a}
Gilles Englebert.
\newblock Conically singular {{Cayley}} submanifolds {{I}}: {{Deformations}}.
\newblock \url{https://doi.org/10.48550/arXiv.2309.07830}, September 2023.

\bibitem[Eng23b]{englebertConicallySingularCayley2023}
Gilles Englebert.
\newblock Conically singular {{Cayley}} submanifolds {{II}}:
  {{Desingularisations}}.
\newblock \url{https://doi.org/10.48550/arXiv.2312.04477}, December 2023.

\bibitem[HL82]{HarvLaws}
Reese Harvey and H.~Blaine Lawson.
\newblock Calibrated geometries.
\newblock {\em Acta Mathematica}, 148:47--157, 1982.

\bibitem[Joy00]{joyceCompactManifoldsSpecial2000}
Dominic Joyce.
\newblock {\em Compact {{Manifolds}} with {{Special Holonomy}}}.
\newblock Oxford {{Mathematical Monographs}}. Oxford University Press, Oxford,
  New York, July 2000.

\bibitem[Joy07]{joyceRiemannianHolonomyGroups2007}
Dominic Joyce.
\newblock {\em Riemannian Holonomy Groups and Calibrated Geometry}.
\newblock Oxford {{Graduate Texts}} in {{Mathematics}}. Oxford University
  Press, 2007.

\bibitem[KL11]{kovalevK3SurfacesNonsymplectic2011}
Alexei Kovalev and Nam-Hoon Lee.
\newblock K3 surfaces with non-symplectic involution and compact irreducible
  {{G2-manifolds}}.
\newblock {\em Mathematical Proceedings of the Cambridge Philosophical
  Society}, 151(2):193--218, September 2011.

\bibitem[KLL20]{karigiannisLecturesSurveysG2Manifolds2020}
Spiro Karigiannis, Naichung~Conan Leung, and Jason~D. Lotay, editors.
\newblock {\em Lectures and {{Surveys}} on {{G2-Manifolds}} and {{Related
  Topics}}}, volume~84 of {\em Fields {{Institute Communications}}}.
\newblock Springer US, New York, NY, 2020.

\bibitem[Kol13]{kollar2013rational}
J.~Kollar.
\newblock {\em Rational Curves on Algebraic Varieties}.
\newblock Ergebnisse Der Mathematik Und Ihrer Grenzgebiete. 3. {{Folge}} / a
  Series of Modern Surveys in Mathematics. Springer Berlin Heidelberg, 2013.

\bibitem[Kov03]{kovalevTwistedConnectedSums2003}
Alexei Kovalev.
\newblock Twisted connected sums and special {{Riemannian}} holonomy.
\newblock {\em Fuer die reine und angewandte Mathematik}, 2003(565):125--160,
  December 2003.

\bibitem[Kov09]{KovalevFibration}
Alexei Kovalev.
\newblock Coassociative {{K3}} fibrations of compact {{G2}}-manifolds.
\newblock \url{https://doi.org/10.48550/arXiv.math/0511150}, 2009.

\bibitem[Loc87]{lockhartFredholmHodgeLiouville1987}
Robert~B. Lockhart.
\newblock Fredholm, {{Hodge}} and {{Liouville Theorems}} on {{Noncompact
  Manifolds}}.
\newblock {\em Transactions of the American Mathematical Society}, 301:1--35,
  1987.

\bibitem[McL98]{mcleanDeformationsCalibratedSubmanifolds1998}
Robert~C. McLean.
\newblock Deformations of calibrated submanifolds.
\newblock {\em Communications in Analysis and Geometry}, 4:705--747, 1998.

\bibitem[Moor17]{mooreDeformationTheoryCayley2017}
Kim Moore.
\newblock {\em Deformation Theory of {{Cayley}} Submanifolds}.
\newblock PhD thesis, University of Cambridge, 2017.

\bibitem[Moor19]{mooreCayleyDeformationsCompact2019}
Kim Moore.
\newblock Cayley deformations of compact complex surfaces.
\newblock {\em Journal of the London Mathematical Society}, 100(2):668--691,
  2019.

\bibitem[MS94]{mcduffJholomorphicCurvesQuantum1994}
Dusa McDuff and Dietmar Salamon.
\newblock {\em J-Holomorphic {{Curves}} and {{Quantum Coholomogy}}}.
\newblock Number~6 in University {{Lecture Series}}. American Mathematical
  Society, 1994.

\bibitem[{\v S}ok80]{sokurovSmoothnessGeneralAnticanonical1980}
V.~V. {\v S}okurov.
\newblock Smoothness of the general anticanonical divisor on a {{Fano}} 3-fold.
\newblock {\em Mathematics of the USSR-Izvestiya}, 14(2):395, April 1980.

\bibitem[Voi02]{voisinHodgeTheoryComplex2002}
Claire Voisin.
\newblock {\em Hodge {{Theory}} and {{Complex Algebraic Geometry I}}}, volume~1
  of {\em Cambridge {{Studies}} in {{Advanced Mathematics}}}.
\newblock Cambridge University Press, Cambridge, 2002.

\end{thebibliography}

\medskip

\noindent{\small\sc The Mathematical Institute, Radcliffe
Observatory Quarter, Woodstock Road, Oxford, OX2 6GG, U.K.

\noindent E-mail: {\tt gilles.englebert@maths.ox.ac.uk.}
\end{document}